\documentclass[preprint,12pt]{elsarticle}
\usepackage[margin=1.5in]{geometry}

\usepackage{amssymb}
\usepackage{amsmath}
\usepackage{graphicx}
\usepackage{amsfonts}
\usepackage{upgreek}
\usepackage{xcolor}
\usepackage[colorlinks=true, linkcolor=purple, citecolor=cyan, urlcolor=cyan]{hyperref}
\usepackage{enumitem}
\usepackage{graphicx}
\usepackage{float}

\newcommand{\minmod}{\hbox{minmod}}
\newcommand{\sF}{{\cal F}}
\newcommand{\shalf}{{\textstyle \frac{1}{2}}}
\newcommand{\V}[1]{\boldsymbol{#1}}
\newcommand{\sD}{{\cal D}}
\long\def\HIDE#1{}
\newcommand{\vthin}{\mkern 1mu} 

\usepackage{lineno}

\journal{Journal of Computational Physics}

\begin{document}

\begin{frontmatter}



\title{Numerical Study of Reported Multiple Solutions to a Cauchy Problem for the 2d Incompressible Euler Equations}


\author[inst1,inst2]{Allen M. Tesdall\corref{cor1}}
\ead{allen.tesdall83@login.cuny.edu}

\author[inst3]{Richard Sanders}
\ead{rsanders@uh.edu} 

\cortext[cor1]{Corresponding author}

\affiliation[inst1]{
    organization={Department of Mathematics, City University of New York, 
    College of Staten Island},
    city={Staten Island},
    postcode={10314},
    state={NY},
    country={USA}
}

\affiliation[inst2]{
    organization={Physics Program, CUNY Graduate Center},
    city={New York},
    postcode={10016},
    state={NY},
    country={USA}
}

\affiliation[inst3]{
    organization={Department of Mathematics, University of Houston},
    city={Houston},
    postcode={77204},
    state={TX},
    country={USA}
}

\begin{abstract}

Bressan and Shen (2021) jointly with Bressan and Murray (2020)
proposed a simple example of a Cauchy problem on $\mathbb{R}^2$, 
characterized by two key parameters, 
intended to yield nonunique weak solutions of the incompressible Euler equations.
They suggest that a limiting process applied to two specific approximations of their
example initial data yields two distinct solutions.
Here we provide evidence coming from a carefully designed numerical procedure to 
support if, or when, this intriguing claim is valid.
Some of the numerical challenges encountered in the present study are:
(i)~the proposed initial data is unbounded in a neighborhood of the origin,
(ii)~the proposed problem is posed on all of $\mathbb{R}^2$,
(iii)~several essential symmetries present in the true solution must be captured,
(iv)~the proposed example depends on two parameters which may be critically linked.
We demonstrate that an extremely large computational domain is mathematically 
necessary to capture the core symmetries of the exact free-space problem.  
Accordingly we implement a nested-grid strategy specifically designed for its particular
suitability and efficacy for the problem at hand.  
Our strategy also allows us to achieve extreme grid refinement in the region 
where solution accuracy is most critical.
Within our nested-grid framework, we design an underlying discretization scheme 
that is high-order, robust, and fast.  

We consider a broad collection of allowable parameter pairs, but among all of them 
we find only one that may satisfy the suggested claim.

\end{abstract}




\end{frontmatter}



\section{Introduction}
\label{sec:intro}
The motion of an inviscid, incompressible fluid in $\mathbb{R}^2$ is described by the Euler equations,
\begin{align}
\label{euler-primvar}
    \V{u}_t + (\V{u} \cdot \nabla) \V{u} + \nabla p \hskip1ex &= \hskip1ex \V{0}, \\
    \text{div}\,\V{u} \hskip1ex &= \hskip1ex 0. \nonumber
\end{align}
Here $\V{u} = \V{u}(t,\V{x})$ is the fluid velocity, and the scalar function $p=p(t,\V{x})$ is the pressure, with $\V{x}=(x,y).$  
As is well known, (\ref{euler-primvar}) can be rewritten in terms of the stream function and the 
vorticity as the system
\begin{align}
\label{euler-psiomega}
    \omega_{t} + \V{u} \cdot \nabla \omega =0 , \\
    \nabla^2 \psi = \omega, \nonumber
\end{align}
where $\V{u}$, $\psi$ and $\omega$ are related by 
\begin{align*}
    &\V{u} = \nabla^{\perp} \psi \equiv (-\psi_{y}, \psi_{x}), \ \ 
      \nabla^2 \psi = (v_x - u_y) \equiv \text{curl}\, \V{u} = \omega.
\end{align*}
The Cauchy problem is completed by 
initial data\footnote{Throughout, we use a superscript 0 to denote the initial state at $t=0$.}
$\omega^0(\V{x}) = \omega(0,\V{x})$.

Following the work of Elling \cite{ve2013,ve2016-a,ve2016-b} 
on self-similar solutions of the 2d incompressible Euler equations
with potential nonuniqueness for the Cauchy problem, Bressan and Murray \cite{ab2020} 
and Bressan and Shen \cite{ab2021} considered
particular initial data for (\ref{euler-psiomega})
where their initial vorticity was symmetric about
the origin and supported on two infinite wedges.
Expressed in polar coordinates $(r,\theta)$, this took the specific form
\begin{equation}	
\label{init-omega} 
    \V{x} = (r\cos{\theta}, r\sin{\theta}), \quad \omega^0(\V{x}) = r^{-\alpha} \phi(\theta).
\end{equation}
The parameter $\alpha \in (0,2)$ describes the radial dependence of the initial vorticity, 
while the angular dependence is given by $\phi$ where 
\begin{equation}
  \label{phi-def1}
  \phi(\theta + \pi) = \phi(\theta),  \quad  \phi(\theta) = 0 \ \ \text{if} \ \ \theta \in \left[\frac{\pi}{4}, \pi \right].
\end{equation}
Observe that along rays where $\phi \not= 0$ the vorticity $\omega^0(\V{x}) \to \infty$ as $|\V{x}| \to 0$.

Analyses of the problem given above were presented in \cite{ab2021,ab2020}.
The exact initial data $\omega^0$ is given by 
(\ref{init-omega})--(\ref{phi-def1}); however, the authors considered $L^\infty$~regularizations 
of this singular function $\omega^0$ by two families of bounded initial data
$\omega^0_{1,\epsilon}$, $\omega^0_{2,\epsilon}$, taking
\begin{equation}
\label{approx-omega0}
	\omega^0_{1,\epsilon}(\V{x}) = \begin{cases}
  	\omega^0(\V{x}) & \text{if } |\V{x}| > \epsilon \\
  	\epsilon^{-\alpha} & \text{if }  |\V{x}| \leq \epsilon, 
	\end{cases}
\quad
	\omega^0_{2,\epsilon}(\V{x}) = \begin{cases}
  	\omega^0(\V{x}) & \text{if } |\V{x}| > \epsilon \\
  	0 & \text{if }  |\V{x}| \leq \epsilon.
	\end{cases}
\end{equation}
As $\epsilon \to 0$, both families converge to 
$\omega^0$ in $L^p_{loc}$ for any $1 \le p < 2/\alpha$. 
Observe that $\omega^0_{1,\epsilon}$ is supported on two wedges and a small disk connecting them, 
while $\omega^0_{2,\epsilon}$ is supported on the same two wedges, this time separated by a small disk of zero vorticity.

Suppose $\omega_{1,\epsilon}(t,\cdot)$ resp.\ $\omega_{2,\epsilon}(t,\cdot)$ denote
the solution to the Cauchy problem with initial data $\omega^0_{1,\epsilon}(\cdot)$ resp.\
$\omega^0_{2,\epsilon}(\cdot)$.  
The authors of \cite{ab2021,ab2020} suggest that as $\epsilon \to 0$ 
\begin{equation*}
  \omega_{1,\epsilon}(t,\cdot) \to \omega_1(t,\cdot)
  \ \ \hbox{and}\ \ 
  \omega_{2,\epsilon}(t,\cdot) \to \omega_2(t,\cdot),
\end{equation*}
where $\omega_1(t,\cdot)$ resp.\ $\omega_2(t,\cdot)$ contain one resp.\ two spiral centers.
Moreover, they suggest both limits are self-similar.
According to Remark~1 in \cite{ab2021}, this would demonstrate an ``incurable'' ill-posedness
(nonuniqueness) for the $\epsilon$-limit problem.
In \S\ref{sec:scaling} we discuss self-similarity induced by certain scale-symmetries.

\begin{figure}[t]
\centerline{
  \vbox{\hbox{\includegraphics[width=0.50\textwidth]
        {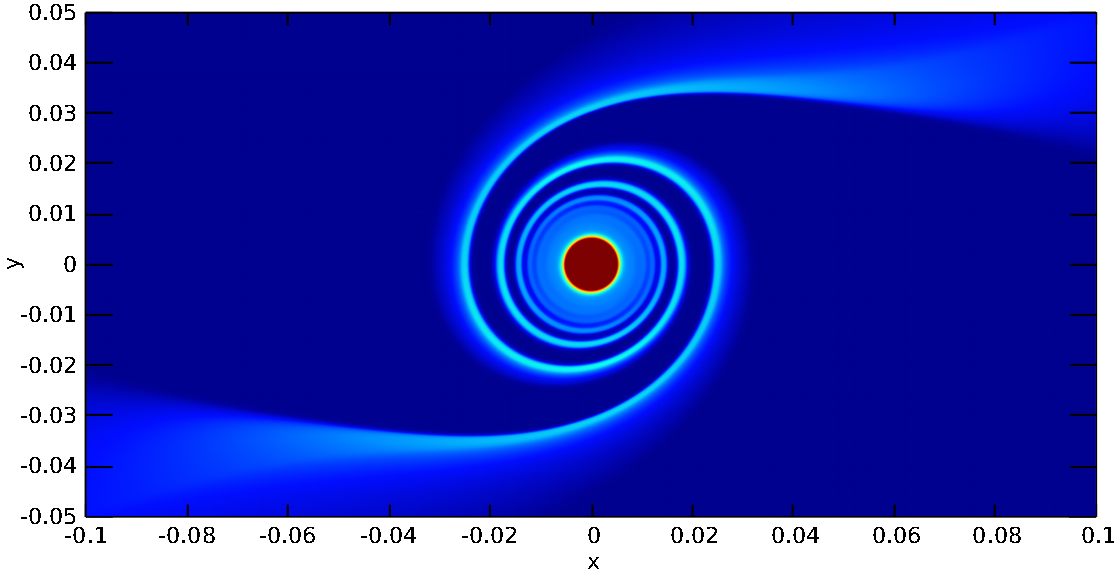}}\hbox{\quad ({\bf a}) One-spiral solution}}
  \hfill
  \vbox{\hbox{\includegraphics[width=0.50\textwidth]
        {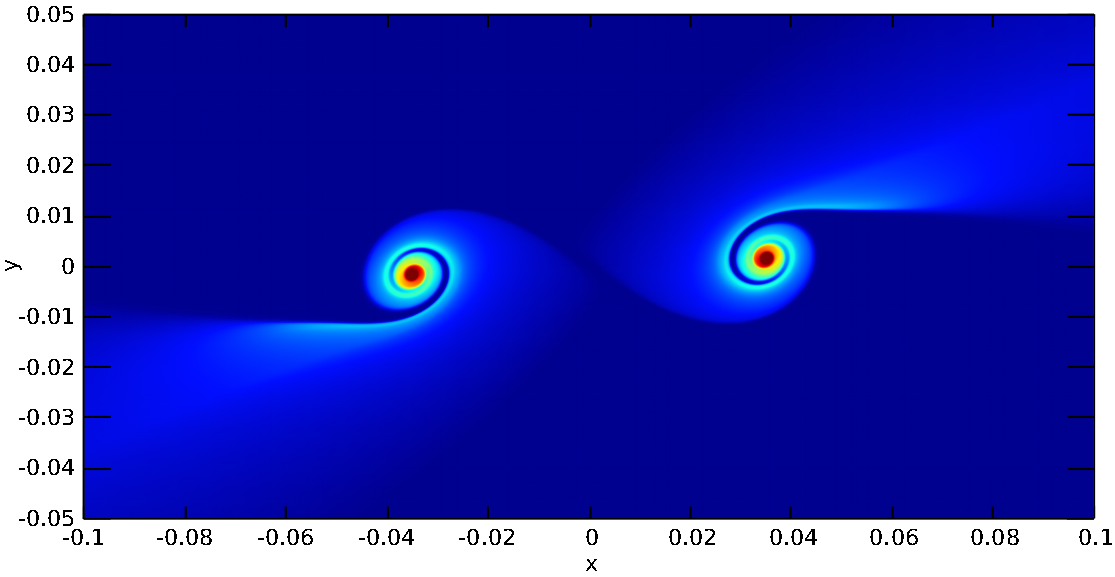}}\hbox{\quad ({\bf b}) Two-spiral solution}}
}
\caption{Our solutions of the two $\epsilon$-regularized Cauchy problems; 
the plots show vorticity contours at $t=1$.  
In\/ {\rm(a)}, solution using initial datum $\omega^0_{1,\epsilon}$, showing a single spiral; 
compare to \cite[Fig. 2]{ab2021}.  
In {\rm(b)}, solution using initial datum $\omega^0_{2,\epsilon}$, 
showing two distinct spirals; compare to \cite[Fig. 3]{ab2021}.
For both {\rm(a)} and {\rm(b)}, $\alpha=0.95$, $\epsilon=0.004$.
The problem domain, $\V{x} \in [-0.2,0.2]\times[-0.2,0.2]$,
is partitioned into a single uniform grid with $4096 \times 4096$ grid cells.
See \S\ref{sec:num-meth} for a full description of the numerical method we employ.}
\label{fig:Our_SingleGrid}
\end{figure}

The analysis in \cite{ab2021,ab2020} employs a sophisticated 
domain decomposition method. 
With this approach they seek to find two distinct self-similar solutions, 
in particular, one containing exactly one spiral and the other containing exactly two spirals.
A critical component of the analysis in \cite{ab2021} relies on the numerical approximation of
the solution of (\ref{euler-psiomega})--(\ref{approx-omega0}) on a bounded intermediate region.
In addition to the analysis, they performed numerical simulations with a Matlab code 
using a finite-difference scheme on a uniform rectangular grid.  
Their code uses
hard wall boundary conditions, $\psi=0$, on a finite domain, set, as of this writing, 
to $[-0.2,0.2] \times [-0.06,0.06]$.
The authors computed solutions to the two $\epsilon$-regularized Cauchy problems, 
and present two vorticity contour plots,
depicting solutions at time $t=1$; see \cite[Figs. 2, 3]{ab2021}.
Their solution with initial datum $\omega^0_{1,\epsilon}$ evolves  
into a single spiral centered at the origin; see \cite[Fig.~2]{ab2021}.  
In their solution with initial datum $\omega^0_{2,\epsilon}$, each wedge in the initial datum 
evolves into its own spiral, with the two spirals symmetric about the origin; 
see \cite[Fig.~3]{ab2021}.
The authors in \cite{ab2021} state that ``numerical simulations indicate that, as $\epsilon \to 0$, 
two distinct limit solutions are obtained''.
In Fig.~\ref{fig:Our_SingleGrid} below we show two specific examples, computed by us, 
illustrating the two generic solutions displayed in \cite{ab2021}.
A single uniform grid of the spatial domain $[-0.2,0.2]\times[-0.2,0.2]$ is 
employed here.  The size of the displayed region is indicated and
parameter values and other problem data are given in the figure caption.  
(These data are not provided in \cite{ab2021} for the solutions presented therein.) 
We note that Bressan and Murray in \cite[Figs. 2, 3]{ab2020}
display contour plots similar to those in \cite{ab2021} and state 
``numerical simulations show the existence of two distinct solutions for the same initial data''.

We introduce the following terminology:
\begin{enumerate}
    [label={Problem \arabic*:}, leftmargin=*, labelindent=1.5\parindent, labelsep=1ex, itemsep=-4pt, 
    topsep=4pt] 
    \item System (\ref{euler-psiomega}) with initial datum $\omega^0_{1,\epsilon}$. 
    \item System (\ref{euler-psiomega}) with initial datum $\omega^0_{2,\epsilon}$. 
\end{enumerate}
The cut-off function $\phi(\theta)$ given in (\ref{phi-def1}) is not yet fully defined; it will be 
in~\S\ref{subsec:init-data}.

In this paper we present high resolution numerical solutions of Problem~1 
and Problem~2 in order to
provide numerical evidence supporting, or not, the assertion of nonuniqueness given in \mbox{\cite{ab2021,ab2020}}.
These problems are posed on all of $\mathbb{R}^2$, while a numerical solution is computed only on some bounded domain.
We review in \S\ref{sec:scaling} certain basic scaling symmetries satisfied by the full free-space problem and self-similarity of $\epsilon$-limits.
We have dealt with self-similar solutions to related problems before; see e.g.\ \cite{t+s+k2008}.  
There we employed self-similar coordinates.
Here it is more natural, because of infinite speed of propagation, 
to employ a very large spatial domain to minimize the effects of far-field numerical boundaries.
To this end, in \S\ref{subsec:num-meth-specifics} we describe a computationally efficient and easy to program nested-grid strategy.  
This strategy is also a means of implementing local grid refinement in an appropriate neighborhood around the origin.  
The basic single grid high-order scheme, which is an integral component of the nested-grid strategy, is derived in \S\ref{subsec:num-meth-general}.

If one attempts to numerically solve Problem~2 on a domain of inadequate size, 
one immediately discovers that the computed movements in time of spiral
center locations (we refer to these as {\em time-trajectories})
do not even approximately satisfy certain key symmetries.
These and other symmetries are derived in \S\ref{sec:scaling} and exemplified in detail in \S\ref{sec:discussion}.  
In \S\ref{sec:numerical}, we present solutions {\em exclusively} computed
using our nested-grid approach on a very large spatial domain.  
We display a certain scaling of spiral time-trajectories for Problem~2 solutions to 
demonstrate that our computed solutions indeed satisfy these symmetries to a high degree of accuracy.  
Having verified this, we use our finest grid results to perform a systematic study designed to 
approximate solution limits as $\epsilon$ tends to zero at a variety of fixed $t > 0$,
maintaining the ratio $\epsilon/\Delta{x}$ at an appropriately large value to accurately resolve the solutions near the origin.
Recall that in Fig.~\ref{fig:Our_SingleGrid} we took a particular example with $\epsilon = 0.004$.  
In \S\ref{sec:numerical} we display numerical solutions of Problem~1 and Problem~2 
for an entire suite of parameter values which are specified at the top of \S\ref{subsec:num-solutions}.
We remark that in our numerical work we are particularly focused on fast execution times so that we 
are able to treat very large spatial domains and integrate to sufficiently large values of $t$.
For our solutions of Problem~2 we carefully monitor spiral center 
locations as functions of time and $\epsilon$.
The results we present in \S\ref{sec:numerical} may be considered somewhat surprising. 

As mentioned, \S\ref{sec:discussion} gives a detailed discussion comparing and
contrasting significant issues and their implications.  
Finally, our conclusions are summarized in \S\ref{sec:conclusion}.


\section{Scaling basics and similarity}
\label{sec:scaling}

Recall from (\ref{init-omega}) and (\ref{approx-omega0}),
the initial data we consider in the present work will take one of two 
specific forms, where, when written in polar coordinates
\begin{align*}
  & \V{x} = (r\cos\theta,r\sin\theta),\\
  &\omega^0_{1,\epsilon}(\V{x}) \equiv 
        \begin{cases}
          r^{-\alpha}\,\phi(\theta) & \text{if } |r| > \epsilon \\
          \epsilon^{-\alpha}        & \text{if } |r| \leq \epsilon
	\end{cases}
   \ \ \hbox{or}\ \
   \omega^0_{2,\epsilon}(\V{x}) \equiv 
        \begin{cases}
          r^{-\alpha}\,\phi(\theta) & \text{if } |r| > \epsilon \\
           0                        & \text{if } |r| \leq \epsilon.
        \end{cases}
\end{align*}
The angular component, $\phi(\theta)$, is a given but not yet fully 
specified nonnegative $\pi$-periodic function.
Suppose $\omega^0_\epsilon(\V{x})$ is defined by either of these two functions.
Observe that since $\V{x}$ is allowed to range over all $\mathbb{R}^2$, 
and $\epsilon$ is positive
\begin{equation}
\label{idata-scaling}
       \omega^0_\epsilon(\epsilon\,\V{x}) 
     = \epsilon^{-\alpha}\omega^0_1(\V{x}). \quad \hbox{(Here the subscript $1$ means $\epsilon=1$.)}
\end{equation}
Next, write the stream function-vorticity system (\ref{euler-psiomega})
in compact scalar form
\begin{align}
\label{unscaled-equation}
  &\omega_t + \nabla^{\perp}_{\V{x}}
             \left( \Delta^{-1}_{\V{x}} \omega \right)
             \cdot 
             \nabla_{\V{x}}\, \omega = 0, \\
  &\omega(0,\V{x}) = \omega^0_\epsilon(\V{x}),\nonumber
\end{align}
and consider scaled variables, $\widehat{t}$, $\widehat{\V{x}}$ and 
$\widehat{\omega}$, defined by
\begin{equation*}
  t = T \, \widehat{t},\ \ 
  \V{x} = X \widehat{\V{x}},\ \ 
  \omega(t,\V{x}) = W\widehat{\omega}(\widehat{t},\widehat{\V{x}}),
\end{equation*}
where the positive constants $T$, $X$ and $W$ are determined as follows.
Introduce these into (\ref{unscaled-equation}) and apply the chain rule 
to find $\widehat{\omega}(\widehat{t},\widehat{\V{x}})$ solves
\begin{align*}
  &\widehat{\omega}_{\widehat{t}} + TW\,\nabla^{\perp}_{\widehat{\V{x}}}
             \left( \Delta^{-1}_{\widehat{\V{x}}} \widehat{\omega} \right)
             \cdot 
             \nabla_{\widehat{\V{x}}}\, \widehat{\omega} = 0, \\
  &W\,\widehat{\omega}(0,\widehat{\V{x}}) 
    = \omega(0,\V{x}) = \omega^0_\epsilon(\V{x}) 
    = \omega^0_\epsilon(X\widehat{\V{x}}).
\end{align*}
Take $X=\epsilon$ and $W=\epsilon^{-\alpha}$ and use (\ref{idata-scaling}) 
on the right side above to reveal that now
$\widehat{\omega}(0,\widehat{\V{x}}) =    \omega^0_1(\widehat{\V{x}})$.
Continue by taking $TW=1$ to finally deduce that
\begin{equation}
  \widehat{\omega}(t,\V{x}) 
  = \epsilon^{\alpha}\omega(\epsilon^\alpha t,\epsilon\vthin\V{x})
\label{solution-to-normalized}
\end{equation}
is a solution of the following {\it normalized} (i.e.\ $\epsilon=1$) 
initial value problem
\begin{align}
\label{scaled-equation}
  &\widehat{\omega}_{t} + \nabla^{\perp}_{\V{x}}
             \left( \Delta^{-1}_{\V{x}} \widehat{\omega} \right)
             \cdot 
             \nabla_{\V{x}}\, \widehat{\omega} = 0, \\
  &\widehat{\omega}(0,\V{x}) = \omega^0_1(\V{x}). \nonumber
\end{align}
The hats over $\V{x}$ and $t$ are suppressed here simply to tidy-up the notation.
To~summarize, if we take
\begin{equation}
\label{normalized-epsilon-scaling}
  t = \epsilon^\alpha \vthin \widehat{t},\ \ 
  \V{x} = \epsilon\vthin \widehat{\V{x}},\ \ 
  \omega(t,\V{x}) 
        = \epsilon^{-\alpha}\vthin \widehat{\omega}(\widehat{t},\widehat{\V{x}}), 
\end{equation}
where $\omega(t,\V{x})$ denotes a solution of the unscaled initial value 
problem (\ref{unscaled-equation}), then $\widehat{\omega}(t,\V{x})$ given 
in (\ref{solution-to-normalized}) solves the normalized problem 
(\ref{scaled-equation}).
More generally, it follows that if $\omega_{\epsilon_0}(t,\V{x})$ 
denotes a solution of (\ref{unscaled-equation}) for one fixed 
$\epsilon=\epsilon_0>0$, then $\omega_\epsilon(t,\V{x})$ given by
\begin{equation}
\label{basic-scaling}
  \omega_{\epsilon}(t,\V{x}) 
  =  \rho^\alpha 
     \omega_{\epsilon_0}(\rho^\alpha t, \rho\,\V{x}),
     \ \hbox{where $\rho = \epsilon_0/\epsilon$,}
\end{equation}
solves (\ref{unscaled-equation}) for arbitrary $\epsilon>0$.

We continue our discussion concerning observable consequences due to scaling symmetry. 
Let $\V{x}_{\epsilon}(t)\in\mathbb{R}^2$ denote some clearly identifiable 
$t$-parameterized curve associated to the vorticity field $\omega_\epsilon(t,\V{x})$. 
For example, this curve could represent the spatial location of $\omega_\epsilon$'s 
locally peak value with varying time.
For the $\epsilon$-family of such curves, the scaling identity (\ref{basic-scaling}) 
just given shows
\begin{align}
\label{traj-symmetry}
  &\omega_\epsilon(t,\V{x}_\epsilon(t)) 
  = 
  \rho^\alpha\omega_{\epsilon_0}(\rho^\alpha t,\rho\vthin\V{x}_\epsilon(t))
  \equiv
  \rho^\alpha\omega_{\epsilon_0}(\rho^\alpha t,\V{x}_{\epsilon_0}(\rho^\alpha t)), \\
  &\Rightarrow\quad
   \rho\vthin \V{x}_\epsilon(t) = \V{x}_{\epsilon_0}(\rho^\alpha t),\ 
                  \hbox{where again $\rho=\epsilon_0/\epsilon$.} \nonumber
\end{align}
Restated, given a curve $\V{x}_{\epsilon_0}(t)$ for one $\epsilon_0>0$, curves
$\V{x}_\epsilon(t)$ can be generated via (\ref{traj-symmetry}) for 
arbitrary $\epsilon>0$.
We use this in \S\ref{sec:numerical} to assess in part our numerical results.
Approximations of $\V{x}_\epsilon(t)$ are extracted from our computed 
numerical data, and the degree to which they satisfy this ``exact solution'' 
symmetry constraint is examined.

A classical technique used in the study of fluid mechanics is to seek 
special, what are called, {\it self-similar} solutions;
see e.g.\ \cite{kgg1947}, \cite{t+h2002}, \cite{t+s+k2008}.  
From~(\ref{normalized-epsilon-scaling}) see that $\V{y}$ and $\Psi$
given by
\[
   \frac{\widehat{\V{x}}}{\widehat{t}^{\, 1/\alpha}} 
 = \frac{\V{x}}{t^{1/\alpha}} \equiv \V{y}
   \ \ \hbox{and}\ \ 
   \widehat{t}\,\widehat{\omega}(\widehat{t},\widehat{\V{x}}) 
 = t\,\omega(t,\V{x}) \equiv  \Psi
\]
are invariant under scaling.
These serve here as self-similar independent and dependent variables.
A self-similar solution, say $\Omega$, is one that only depends on 
variable $\V{y}$.
Let $t>0$ and $\epsilon_0>0$ be fixed, and suppose $\tau \ge t$ but
otherwise arbitrary.
Set $\epsilon/\epsilon_0 = \left( t/\tau \right)^{1/\alpha}$
in (\ref{basic-scaling}) to see
\[
  t\,\omega_\epsilon(t,\V{x})
  =
  \tau\,\omega_{\epsilon_0}(\tau, \tau^{1/\alpha}(\V{x}/t^{1/\alpha}) )
  =
  \tau\,\omega_{\epsilon_0}(\tau, \tau^{1/\alpha}\V{y} ).
\]
Note that $\epsilon \to 0 \iff \tau \to \infty$.
Therefore, at least formally, 
\begin{equation}
  \label{eps-limit-of-weps}
  t\,\omega(t,\V{x})
  \equiv
  \lim_{\epsilon \to 0}
    \Big( t\,\omega_\epsilon(t,\V{x}) \Big)
  = 
  \lim_{\tau \to \infty}
    \Big(\tau\,\omega_{\epsilon_0}(\tau, \tau^{1/\alpha}\V{y} )\Big)
  \equiv
  \Omega(\V{y}).
\end{equation}
Using $\omega(t,\V{x}) = t^{-1}\Omega(\V{y})$ in 
(\ref{unscaled-equation}) gives the following self-similar formulation
\begin{align*}
 \hbox{$t>0$}\ \ &\Rightarrow\ \  
    \left( \nabla^\perp_{\V{y}} \Delta^{-1}_{\V{y}} \Omega
        -\frac{1}{\alpha}\,\V{y} \right) \cdot \nabla_{\V{y}}\,\Omega
    = \Omega,\\
 \hbox{$t\to0$}\ \ &\Rightarrow\ \  
    \lim_{|\V{y}| \to \infty} |\V{y}|^\alpha\, \Omega(\V{y}) = \phi(\theta).
\end{align*}
Bressan and Shen in \cite{ab2021} and Bressan and Murray in \cite{ab2020}
consider precisely this self-similar free-space boundary value problem.

For a certain range of parameter $\alpha$ and 
angular cut-off function $\phi(\theta)$,
the authors of \cite{ab2021} and \cite{ab2020} strongly suggest the 
self-similar equation above admits two completely distinct solutions,
one containing a single spiral, say $\Omega_1(\V{y})$, 
and the other containing exactly two, say $\Omega_2(\V{y})$.
(By spiral we mean an isolated point of maximum/infinite vorticity.)
These were obtained via a sophisticated construction.
Specifically, an approximate solution of Cauchy problem~(\ref{unscaled-equation})
with Problem~1's initial datum, i.e.\ $\omega^0_{1,\epsilon_0}(\V{x})$, 
was employed to build $\Omega_1(\V{y})$.
Likewise, an approximate solution with Problem~2's initial datum, 
$\omega^0_{2,\epsilon_0}(\V{x})$, was used to build $\Omega_2(\V{y})$.
(Here we use subscripts $1$~resp.~$2$ to refer to solution with
Problem~1~resp.~Problem~2 initial data.)\ \ 
In \cite{ab2020} it is further suggested that, at an ``intuitive level'' 
in some strong a.e.~sense, for $t \ge t_0 > 0$
\begin{equation}
  \label{prob1-prob2-limits}
  (i)\ \ 
  t\,\omega_{1,\epsilon}(t,\V{x}) \to \Omega_1(\V{y}) 
  \ \ \hbox{and}\ \ 
  (ii)\ \ 
  t\,\omega_{2,\epsilon}(t,\V{x}) \to \Omega_2(\V{y})
  \ \ \hbox{as $\epsilon\to 0$.}
\end{equation}
Whether $(i)$ and/or $(ii)$ can be confirmed or not by numerical experiment is
the principal topic of this present work.

We consider $(ii)$ in (\ref{prob1-prob2-limits}) the most 
interesting, and in fact a critical so-called intuitive aspect of the
present problem.
The original time dependent Cauchy problem as well as its self-similar 
analogue both possess antipodal symmetry.
Therefore, $\Omega_1$'s single spiral location must lie at the origin,
whereas $\Omega_2$'s two distinct spiral locations are symmetric about 
and therefore do not lie at the origin. 
We can easily investigate numerically if $(ii)$ seems likely as follows.
From~(\ref{eps-limit-of-weps}), we expect in a general sense 
$t\,\omega_{2,\epsilon}(t,\V{x}) \sim \Omega(\V{y})$ as $\epsilon \to 0$.
To determine whether $\Omega(\V{y})$ contains a spiral center
$\V{y}_* \not = \V{0}$, for $\epsilon>0$ we track the location 
of peak vorticity in time found during numerical time integration.
Let $\V{x}_\epsilon(t)$ denote this path.
From
\[
  t\,\omega_{2,\epsilon}(t,\V{x}_\epsilon(t)) 
    \sim \Omega(\V{x}_\epsilon(t)/t^{1/\alpha})
  \ \ \hbox{we expect to see}\ \
  \V{x}_\epsilon(t)/t^{1/\alpha} \sim \V{y}_*.
\]
That is, when $\V{y}_* \not= \V{0}$ and for adequately large $t>0$, 
we expect
\begin{equation}
  \label{what-we-expect}
  \V{x}_\epsilon(t) \sim t^{1/\alpha}\,\V{y}_*,
\end{equation}
where the path on the right hand side above is a straight line emanating
from the origin and has easily observable speed.


\section{The numerical methods}
\label{sec:num-meth}

\subsection{Discretization generalities.}
\label{subsec:num-meth-general}

The spatial domain $\Omega =  [-L,L] \times [-L,L]$ is partitioned into a 
uniform grid of cells, $\bigcup_{i,j} (x_i, x_{i+1}) \times (y_j,y_{j+1})$, 
where $\Delta x = x_{i+1}-x_i$ and $\Delta y = y_{j+1}-y_j$ are constant.
Regard $\omega_{i,j}$ as a cellwise approximation of the
vorticity $\omega$ from (\ref{euler-psiomega}) at cell centers 
$(x_{i+1/2},y_{j+1/2})$.
Regard $\psi_{i,j}$ as a pointwise approximation of the stream function
$\psi$ living at cell vertices $(x_i,y_j)$ determined by inverting 
the second centered difference scheme for $\nabla^2 \psi = \omega$
\begin{equation}
\label{stream-function-discretization}
\begin{aligned}
  \sD^2(\psi)_{i,j} &:= 
    \frac{\left(\psi_{i+1,j}-2\psi_{i,j}+\psi_{i-1,j}\right)}{\Delta x^2}
  + \frac{\left(\psi_{i,j+1}-2\psi_{i,j}+\psi_{i,j-1}\right)}{\Delta y^2}\\
  &= \nu(\omega)_{i,j}, 
\end{aligned}
\end{equation}
where $\nu(\omega)\approx\omega$ denotes the mass lumping operator taking
cell-centered quantities to vertex-centered; specifically
\begin{equation*}
  \nu(\omega)_{i,j} := \frac{1}{4} 
   \left( \omega_{i,j} + \omega_{i-1,j} 
                       + \omega_{i,j-1} + \omega_{i-1,j-1} \right).
\end{equation*}
Boundary conditions for (\ref{stream-function-discretization}) 
are discussed in the next subsection.
This discrete problem is inverted by a multigrid \hbox{V-cycle} algorithm 
where pre-\ and post-smoothing is accomplished by three iterations of 
artificial time relaxation applied to the residual
(i.e.~$\psi^{k+1} = \psi^k + \rho\,(\sD^2(\psi^k)-\nu(\omega))$ with relaxation
parameter $\rho = 0.25/(1/\Delta x^2 + 1/\Delta y^2)$) together with 
standard intergrid transfer operators.
On cell edge $(x_{i},y_j)~\to~(x_i,y_{j+1})$ (i.e.~left cell edge) 
the discrete velocity's $x$-component is taken as
\begin{equation*}
  u_{i,j} = -\frac{1}{\Delta y} \left( \psi_{i,j+1} - \psi_{i,j}\right)
  \approx -\frac{\partial}{\partial y}\psi(x_i,y_{j+1/2}),
\end{equation*}
and on edge $(x_{i},y_j)~\to~(x_{i+1},y_j)$ (i.e.~bottom cell edge)
its $y$-component is taken as
\begin{equation*}
  v_{i,j} = \frac{1}{\Delta x} \left( \psi_{i+1,j} - \psi_{i,j}\right)
  \approx \frac{\partial}{\partial x}\psi(x_{i+1/2},y_j).
\end{equation*}
Both are second order approximations at edge centers.
Of course, the stream function formulation for velocity in 
(\ref{euler-psiomega}) automatically leads to $\nabla \cdot \V{u} = 0$,
and here the discrete analogue reads
\begin{equation}
\label{discrete-div-free-vel}
  \frac{u_{i+1,j}-u_{i,j}}{\Delta x}
 +\frac{v_{i,j+1}-v_{i,j}}{\Delta y} = 0.
\end{equation}
It is natural therefore to rewrite the vorticity equation in conservation form
\begin{equation}
\label{consform-advection}
  \frac{\partial\, \omega}{\partial t} + \V{u} \cdot \nabla\,\omega
  = 0 
  \quad\Rightarrow\quad
  \frac{\partial\, \omega}{\partial t}
  + \frac{\partial}{\partial x} \left(u\,\omega\right)
  + \frac{\partial}{\partial y} \left(v\,\omega\right) = 0,
\end{equation}
and to then discretize the resulting conservation law via finite volume.

A robust first order accurate difference scheme for the 
vorticity conservation law on the right in (\ref{consform-advection})
is explicit upwind differencing
\begin{equation*}
   \frac{ \omega^{n+1}_{i,j} -  \omega^{n}_{i,j} }{\Delta t}
 + \frac{ f^n_{i+1,j} - f^n_{i,j}}{\Delta x}
 + \frac{ g^n_{i,j+1} - g^n_{i,j}}{\Delta y} = 0,\\
\end{equation*}
where the numerical fluxes $f$ and $g$ are given by the following upwind
formulae
\begin{align*}
  f^n_{i,j} &= \min(u^n_{i,j},0)\,\omega^n_{i,j} 
             + \max(u^n_{i,j},0)\,\omega^n_{i-1,j}, \\ 
  g^n_{i,j} &= \min(v^n_{i,j},0)\,\omega^n_{i,j} 
             + \max(v^n_{i,j},0)\,\omega^n_{i,j-1}.
\end{align*}
The velocity field, $(u^n_{i,j},v^n_{i,j})$ above, should be updated at every 
time step as outlined in the previous paragraph.

The basic first order upwind scheme above is extended to a 
second order method, see e.g.\ \cite{bvl1974}, \cite{pks1984}, in 
the following manner.
The minmod slope limiter 
\begin{equation*}
  \minmod(a,b) = 
    \begin{cases}
       a & \text{if } |a| < |b|,\ ab > 0 \\
       b & \text{if } |b| < |a|,\ ab > 0 \\
       0 & \text{otherwise,}
    \end{cases}
\end{equation*}
is used to define four edge values (Left, Right, Bottom, Top) per cell
\begin{equation}
\label{limited-edge-points}
\begin{aligned}
  s^{L|R}_{i,j} &= \omega_{i,j} \mp \frac{1}{2}\,
    \minmod (\ \omega_{i+1,j}-\omega_{i,j},\ \omega_{i,j}-\omega_{i-1,j}\ ),\\
  s^{B|T}_{i,j} &= \omega_{i,j} \mp \frac{1}{2}\,
    \minmod (\ \omega_{i,j+1}-\omega_{i,j},\ \omega_{i,j}-\omega_{i,j-1}\ ).
\end{aligned}
\end{equation}
We employ it here due to its fast execution speed and programming 
simplicity.
Note the following obvious facts for later reference,
\begin{equation}
\label{limited-points-in-I}
\begin{aligned}
  s^{L}_{i,j} \in I(\omega_{i-1,j},\,\omega_{i,j}),&\quad
  s^{R}_{i,j} \in I(\omega_{i+1,j},\,\omega_{i,j}),\quad
  \shalf\left(s^L_{i,j} + s^R_{i,j}\right) = \omega_{i,j}, \\
  s^{B}_{i,j} \in I(\omega_{i,j-1},\,\omega_{i,j}),&\quad
  s^{T}_{i,j} \in I(\omega_{i,j+1},\,\omega_{i,j}),\quad
  \shalf\left(s^B_{i,j} + s^T_{i,j}\right) = \omega_{i,j}, \\
\end{aligned}
\end{equation}
where $I(a,b)$ denotes the interval $[\,\min(a,b),\,\max(a,b)\,]$.

Now, consider the following spatial discretization where here we insert the 
limited cell edge points given in (\ref{limited-edge-points}) into the basic 
first order upwind scheme discussed earlier
\begin{equation}
\label{second-order-spatial}
\begin{aligned}
  \sF(\V{u},\omega&)_{i,j} := 
    \frac{ f_{i+1,j} - f_{i,j}}{\Delta x}
   +\frac{ g_{i,j+1} - g_{i,j}}{\Delta y}, \ \ \text{where}\\
  f_{i,j} &= \min(u_{i,j},0)\,s^L_{i,j} 
             + \max(u_{i,j},0)\,s^R_{i-1,j}, \\ 
  g_{i,j} &= \min(v_{i,j},0)\,s^B_{i,j} 
             + \max(v_{i,j},0)\,s^T_{i,j-1}.
\end{aligned}
\end{equation}
This is a generically second order accurate 
difference method based on $\minmod$ reconstruction from cell averages
coming from the finite volume paradigm.
To march forward in time, 
we are naturally led to employing Heun's two-step method, see e.g.\ \cite{shu1988},
\begin{equation}
\label{full-second-order}
\begin{aligned}
  \hbox{(p)}\qquad
  &\widetilde\omega = \omega^n - \Delta t\, \sF(\V{u}^n,\omega^n), \\
  \hbox{(c)}\qquad
  &\omega^{n+1}     = \omega^n - \shalf\Delta t\,
            \big(
            \sF(\V{u}^n,\omega^n)
          + \sF(\widetilde{\V{u}},\widetilde \omega)
            \big). \\
\end{aligned}
\end{equation}
At each time step $n$, Heun's method requires two numerical inversions of the 
Laplacian: one to determine $\V{u}^n$ from $\omega^n$ and a second to determine 
$\widetilde{\V{u}}$ from $\widetilde \omega$.
Clearly, Heun's method is second order in time.
The predictor step (p) is first order Euler, while the corrector step (c)
is based on an explicit form of the second order trapezoidal rule.

It is interesting to observe that the second order discretization we employ,
that is (\ref{full-second-order}) together with (\ref{second-order-spatial}), 
satisfies the maximum principle
\begin{equation*}
  \min_{i,j}(\omega^n_{i,j}) \le \omega^{n+1}_{i,j} \le \max_{i,j}(\omega^n_{i,j}),
\end{equation*}
subject to a certain CFL-like time step size restriction.
To see this is true, first consider
\begin{equation*}
    \omega_{i,j} - \Delta t\, \sF(\V{u},\omega)_{i,j},
\end{equation*}
which by way of (\ref{limited-points-in-I}) and (\ref{second-order-spatial})
can be written as
\begin{align*}
    &\frac{1}{4} \left(s^L_{i,j}+s^R_{i,j}+s^B_{i,j}+s^T_{i,j}\right) \\
    & - \frac{\Delta t}{\Delta x}
        \big(
        (\min(u_{i+1,j},0)\,s^L_{i+1,j} + \max(u_{i+1,j},0)\,s^R_{i,j})\\
    &\qquad
       -(\min(u_{i,j},0)\,s^L_{i,j} + \max(u_{i,j},0)\,s^R_{i-1,j})
        \big)\\
    & - \frac{\Delta t}{\Delta y}
        \big(
        (\min(v_{i,j+1},0)\,s^B_{i,j+1} + \max(v_{i,j+1},0)\,s^T_{i,j})\\
    &\qquad
      -(\min(v_{i,j},0)\,s^B_{i,j} + \max(v_{i,j},0)\,s^T_{i,j-1})
        \big).
\end{align*}
This can be further rearranged into eight terms 
\begin{align*}
   &\left(
     \frac{1}{4} + \frac{\Delta t}{\Delta x}\min(u_{i,j},0)
    \right) s^L_{i,j}
  + \left(
    -\frac{\Delta t}{\Delta x}\min(u_{i+1,j},0)
    \right)s^L_{i+1,j}\\
 + &\left( 
     \frac{1}{4} - \frac{\Delta t}{\Delta x}\max(u_{i+1,j},0)
    \right) s^R_{i,j}
  + \left(
     \frac{\Delta t}{\Delta x}\max(u_{i,j},0)
    \right)s^R_{i-1,j} \\
 + &\left( 
     \frac{1}{4} + \frac{\Delta t}{\Delta y}\min(v_{i,j},0)
    \right) s^B_{i,j}
 + \left(
    -\frac{\Delta t}{\Delta y}\min(v_{i,j+1},0)\right)s^B_{i,j+1} \\
 + &\left( 
     \frac{1}{4} - \frac{\Delta t}{\Delta y}\max(v_{i,j+1},0)
    \right)s^T_{i,j}
 + \left( 
     \frac{\Delta t}{\Delta y}\max(v_{i,j},0)
   \right)s^T_{i,j-1}.
\end{align*}
Notice that for $\Delta t$ taken small enough, 
the bracketed coefficients above are all nonnegative,
and by (\ref{discrete-div-free-vel}) they sum as follows
\begin{align*}
    1 & + \frac{\Delta t}{\Delta x}
        \big( 
          \min(u_{i,j},0)-\min(u_{i+1,j},0)-\max(u_{i+1,j},0)+\max(u_{i,j},0)
        \big)\\
      & + \frac{\Delta t}{\Delta y}
        \big(
          \min(v_{i,j},0)-\min(v_{i,j+1},0)-\max(v_{i,j+1},0)+\max(v_{i,j},0)
        \big)\\
  = 1 & - \Delta t \left( 
           \frac{u_{i+1,j}-u_{i,j}}{\Delta x} 
         + \frac{v_{i,j+1}-v_{i,j}}{\Delta y}
        \right) = 1.
\end{align*}
Therefore, by convexity, and since the facts referenced
in (\ref{limited-points-in-I}) imply 
$\min(\omega) \le \min(s^L, s^R, s^B, s^T)$ and
$\max(s^L, s^R, s^B, s^T) \le \max(\omega)$, we conclude
\begin{equation*}
  \min(\omega) 
  \le 
  \omega_{i,j} - \Delta t\, \sF(\V{u},\omega)_{i,j}
  \le
  \max(\omega).
\end{equation*}
Finally note that the corrector step in Heun's method, 
(\ref{full-second-order}), can also be rewritten as a convex combination
\begin{equation*}
 \omega^{n+1} = \frac{1}{2}\,\omega^n 
              + \frac{1}{2} \left(
                \widetilde\omega 
              - \Delta t\,\sF(\widetilde{\V{u}},\widetilde{\omega})
                            \right),
\end{equation*}
from which the maximum principle for our second order approach follows.

\subsection{Implementation specifics.}
\label{subsec:num-meth-specifics}
We programmed algorithms (\ref{stream-function-discretization}) and
(\ref{second-order-spatial}), (\ref{full-second-order}) in
the C-language and compiled for an Intel i9 
processor using Intel's 
{\tt icx} compiler with the flag {\tt -Ofast} and 64-bit floats.
Throughout we take $\Delta x = \Delta y$.
Since multigrid is employed to solve (\ref{stream-function-discretization})
we use $2^l$ number of grid points in each coordinate direction.
\hbox{V-cycles} are terminated when the residual max-norm drops below
$5{\cdot}10^{-7}$ which requires at most six cycles independent of grid 
size refinement.

Two of our results presented in the present work are performed
on a single uniform grid of the problem domain $ [-L,L] \times [-L,L]$.
In these cases we take $\psi_{i,j} = 0$ on the domain boundary.
Doing so guarantees $u_{i,j} = 0$ on the left and right and 
$v_{i,j}=0$ on the bottom and top boundary segments.
Therefore, the vorticity boundary flux is zero.
In these single grid cases, we take $L=0.2$.

All other presented results use a sequence of half spatially divided nested 
grids, each using the same number of grid points but with half the spatial 
extent.
That is, for $l=0, 1, \ldots$, sequentially bisect the outermost domain to the final desired grid extent:
\begin{equation}
\label{nested-domains}
  \begin{aligned}
      \Omega_l &= [-L_0/2^{l},L_0/2^{l}]\times[-L_0/2^{l},L_0/2^{l}] \\
                       &\subseteq [-L_0,L_0] \times [-L_0,L_0].
  \end{aligned}
\end{equation}
For related work see e.g.\ \cite{o+s1983}.  
Boundary conditions for $\Omega_0$ are treated 
exactly the same as done in the single grid cases, i.e.~$\psi_{i,j}=0$
for $(x_i,y_j)\in\partial\Omega_0$.
Space/time boundary conditions for $\psi$ and $\omega$ on grid 
$\Omega_l$ with $l\ge 1$ are obtained sequentially by storing 
accurately interpolated solution data from appropriate interior points given by 
the previous grid's (i.e.~the double sized $\Omega_{l-1}$) computed solution.
This grid nesting strategy is a computationally inexpensive and simple to implement
form of local grid refinement.  
We always take $L_0 = 64.0$ and generally stop at $L_0/2^9 = 0.125$.
We exploit this strategy for two additional reasons.
First, we test how taking $\psi = 0$ on the so-called {\it far field} 
$\partial\Omega_0$ affects our computed solutions.
To this end we take $L_0=64.0$ and, specifically for this test, 
we refine only $8$ levels to $L_0/2^8 = 0.25$.
Results are compared to the single grid output discussed in the previous paragraph with $L=0.2$.
Second, we use this strategy to {\it cheaply} test for discretization accuracy.  


\section{Numerical results}
\label{sec:numerical}

\subsection{Defining the initial data.}
\label{subsec:init-data}

The initial data $\omega^0(\V{x})$ as formulated by Bressan and Shen in \cite{ab2021} 
(given here by equations (\ref{init-omega})--(\ref{phi-def1})) 
do not specify what the function $\phi(\theta)$ is, only that it be $\pi$-periodic and supported on two wedges.  
For convenience, as done in \cite{ws2020}, we rotate $\phi(\theta)$ so that it is symmetric about the $x$-axis.
Specifically, we take $\phi(\theta)$ to be the $\pi$-periodic extension of the tent function
\begin{equation}
  \label{phi-def2}
  \phi(\theta) = \begin{cases}
      \theta_w - |\theta|, & \text{if}\ |\theta| \le \theta_w \\
                                0, & \text{if}\ \theta_w < |\theta| \le \pi/2, 
  \end{cases}
\end{equation}
where $\theta_w$ is some given angle, with  $0 < \theta_w < \pi/2$.  
Hence, the initial data $\omega^0_{1,\epsilon}$ and $\omega^0_{2,\epsilon}$, defined in (\ref{approx-omega0}),
depend on three parameters: $\alpha$, $\epsilon$ and also $\theta_w$.

\begin{figure}
  \centerline{
    \vbox{\hbox{\includegraphics[width=0.7\textwidth]{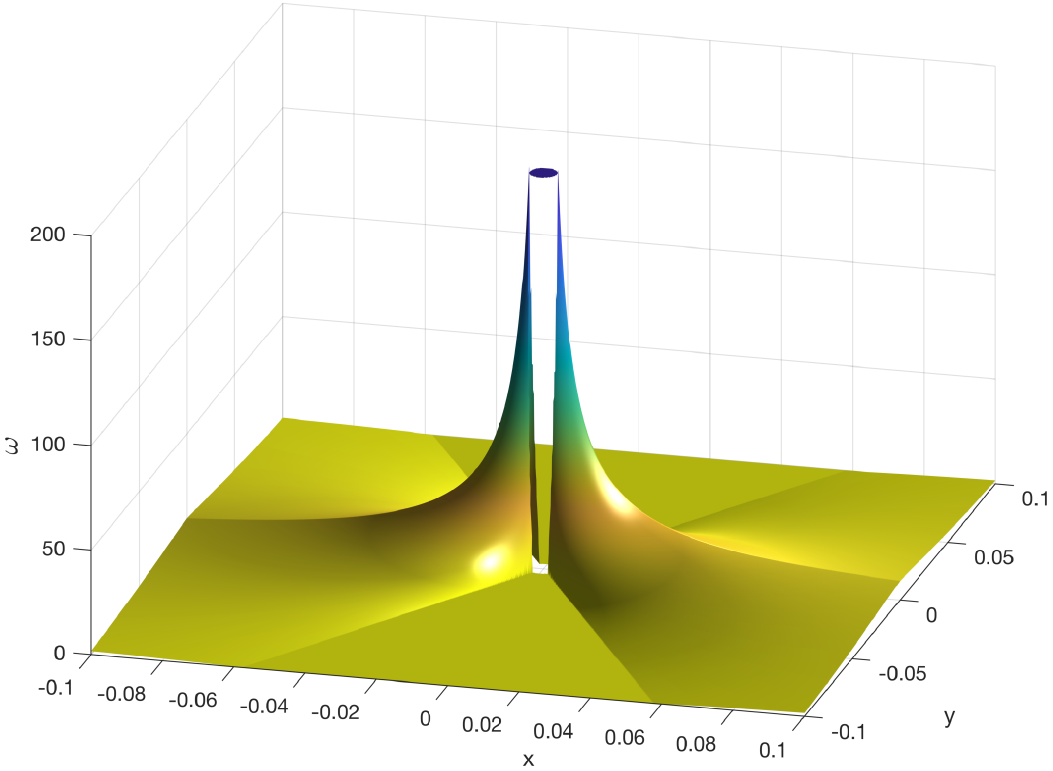}} \hbox{(a) \hskip1em $\omega^0_{1,\epsilon},$
    the initial data for Problem 1.}}
  }
  \centerline{
    \vbox{\hbox{\includegraphics[width=0.7\textwidth]{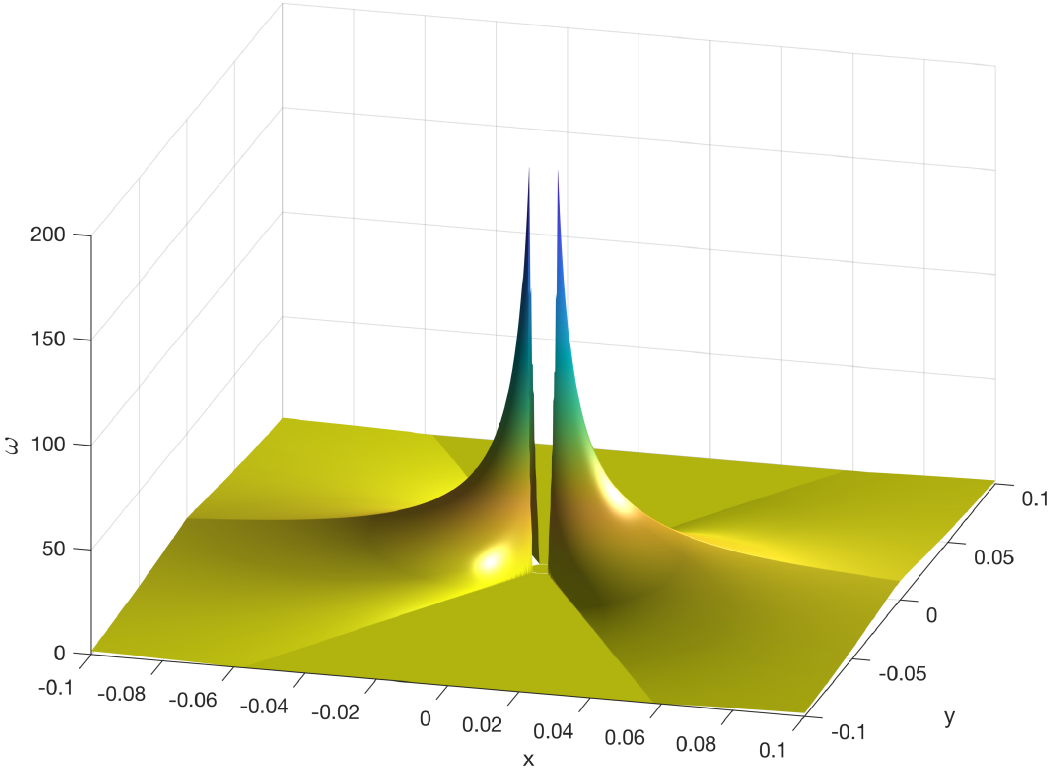}} \hbox{(b) \hskip1em $\omega^0_{2,\epsilon},$
    the initial data for Problem 2.}}
  }
\caption{Surface plots of the two initial vorticity distributions, (\ref{approx-omega0}), with $\phi(\theta)$ as defined in 
	(\ref{phi-def2}).  
	Parameters $\alpha=0.95$, $\epsilon=0.004$, $\theta_{w}=\pi/3$.
	The value of vorticity $\omega$ in the disk at the ``top'' of the surface in 
	(a) is $\epsilon^{-0.95} \approx 189.7$.
	The supremum of $\omega$ in both plots is $\epsilon^{-0.95} \frac{\pi}{3} \approx 198.6$. 
	The graphs in (a) and (b) are displayed only over the indicated portion of the computational domain.
}
\label{fig:InitVort}
\end{figure}

In Fig.~\ref{fig:InitVort} we illustrate these bounded, discontinuous initial data for the example with 
$(\alpha, \epsilon, \theta_{w}) = (0.95, 0.004, \pi/3)$.  
The discontinuity on the circle $r = \epsilon$ and the support of the data are clearly visible in both plots.  

\subsection{The numerical solutions.}
\label{subsec:num-solutions}

In view of the discussion in \S\ref{subsec:init-data} above, we begin by restating:
\begin{enumerate}
    [label={Problem \arabic*:}, leftmargin=*, labelindent=1.5\parindent, labelsep=1ex, itemsep=-4pt, 
    topsep=4pt] 
    \item System (\ref{euler-psiomega}) with initial datum $\omega^0_{1,\epsilon}$ completed by (\ref{phi-def2}).
    \item System (\ref{euler-psiomega}) with initial datum $\omega^0_{2,\epsilon}$ completed by (\ref{phi-def2}).
\end{enumerate}
We solve Problem~1 and Problem~2 for certain combinations of the following:
\begin{align}
\label{eq:parameter-set}
	&\alpha = 0.95,\ 1.05,\ 1.4. \\
	&\epsilon = 0.004,\ 0.002,\ 0.001,\ 0.0005,\ 0.00025,\ 0.000125. \nonumber \\ 
	&\theta_w = \frac{\pi}{8},\ \frac{3\pi}{16},\ \frac{\pi}{4},\ \frac{\pi}{16}. \nonumber
\end{align}
Some of these parameter values are motivated by a related study, \cite{bjl2021}, on 
a small Mach number compressible analogue of our incompressible problem; 
in particular, their Example~3 is most aligned with our incompressible work. 
The authors of \cite{bjl2021} define the angular cut-off function $\phi(\theta)$
and give specific values for initial data parameters $\alpha$ and $\theta_w$.
In our baseline calculations, see \S\ref{subsubsec:baseline-param-values},
we match the values of $\alpha$ and $\theta_w$ used in their Example~3.
The ratios of successive $\epsilon$'s we give above were chosen to be $1/2$
for convenience due to our scaling arguments in \S\ref{sec:scaling}.

Throughout, we use a sequence of nested grids as described in (\ref{nested-domains}).  
We take the outermost extent $L_0=64.0$ in order to reasonably approximate the free-space
problem on which the scaling discussed in \S\ref{sec:scaling} is valid.  
Solutions are computed using a sequence of either 10 or 11 nested subdomains, $\Omega_l$, where
\begin{align}
\label{eq:nesting}
    \Omega_0      &= [-64.0,64.0] \times [-64.0,64.0], \ \cdots, \\
    \Omega_{9}   &= [-0.125,0.125] \times [-0.125,0.125], \nonumber \\  
    \Omega_{10} &= [-0.0625,0.0625] \times [-0.0625,0.0625]. \nonumber
\end{align}
Each $\Omega_l$ is covered by a computational grid comprised of $4096^2$ uniform cells.  
Only solutions on the innermost/most-refined grids are presented in what follows.

\subsubsection{Baseline parameter calculations}
\label{subsubsec:baseline-param-values}

Here and throughout this entire subsection we use the following set of parameter values:  
$\alpha=0.95$, $\theta_w=\pi/8$, and selected values of $\epsilon$ given in (\ref{eq:parameter-set}).

\begin{figure}
  \centerline{
    \vbox{\hbox{\includegraphics[height=0.22\textheight]{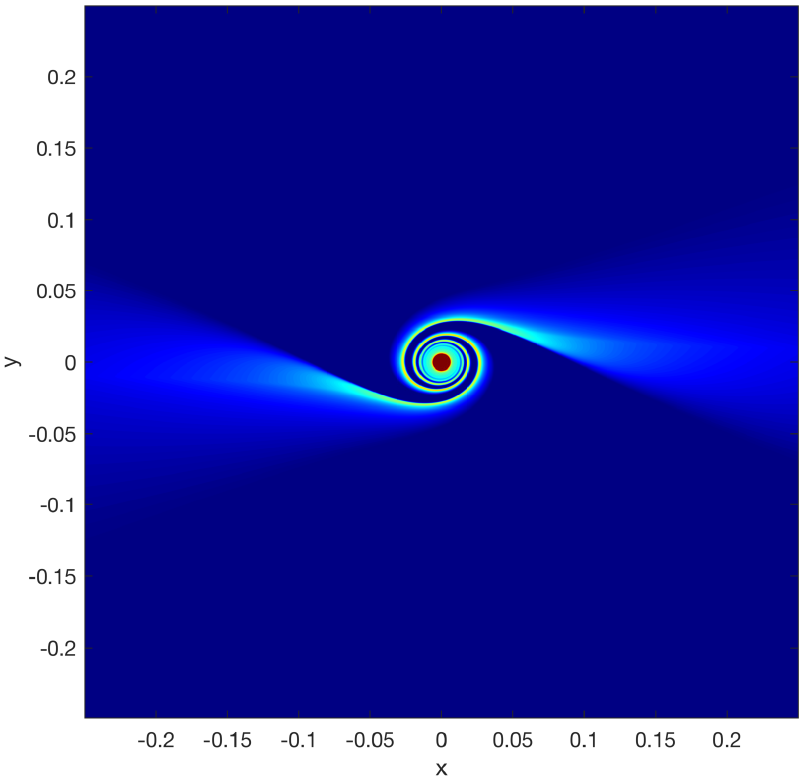}} 
    \hbox{\quad (a) {\it $\omega$}-contours, $t=1.0$}} 
    \qquad
    \vbox{\hbox{\includegraphics[height=0.25\textheight]{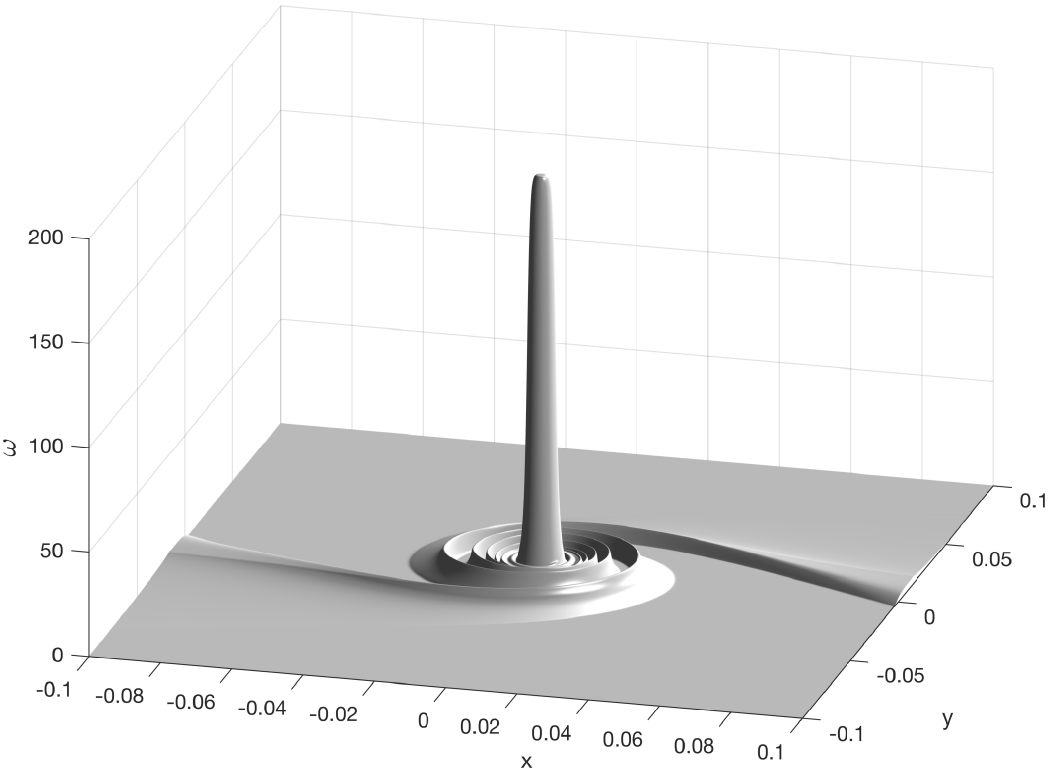}}
    \hbox{\quad (b) {\it $\omega$}-surface, $t=1.0$}}
  } 
  \medskip
  \centerline{
    \vbox{\hbox{\includegraphics[height=0.22\textheight]{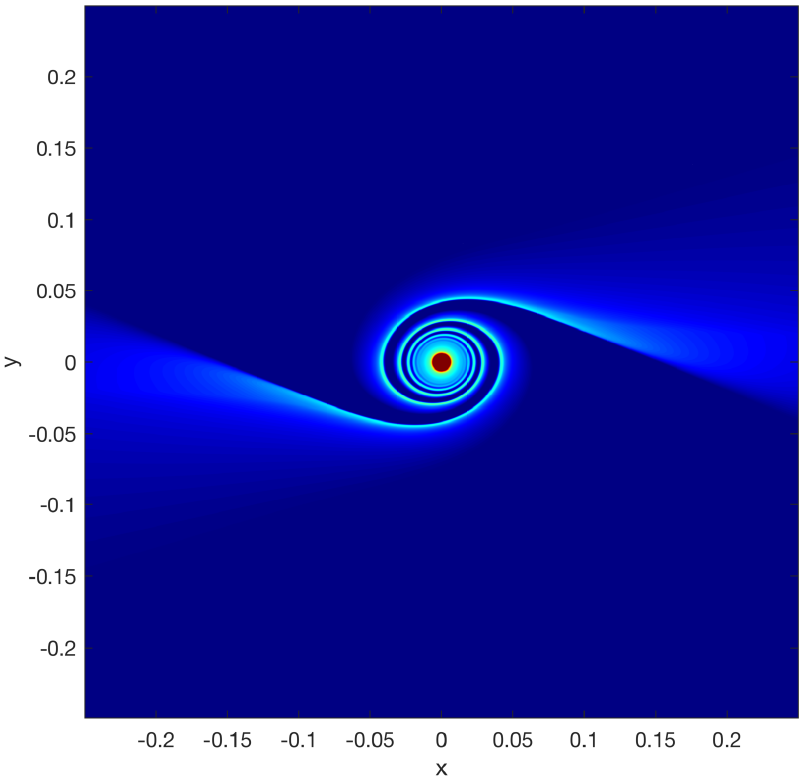}} 
    \hbox{\quad (c) {\it $\omega$}-contours, $t=2.0$}} 
    \qquad
    \vbox{\hbox{\includegraphics[height=0.25\textheight]{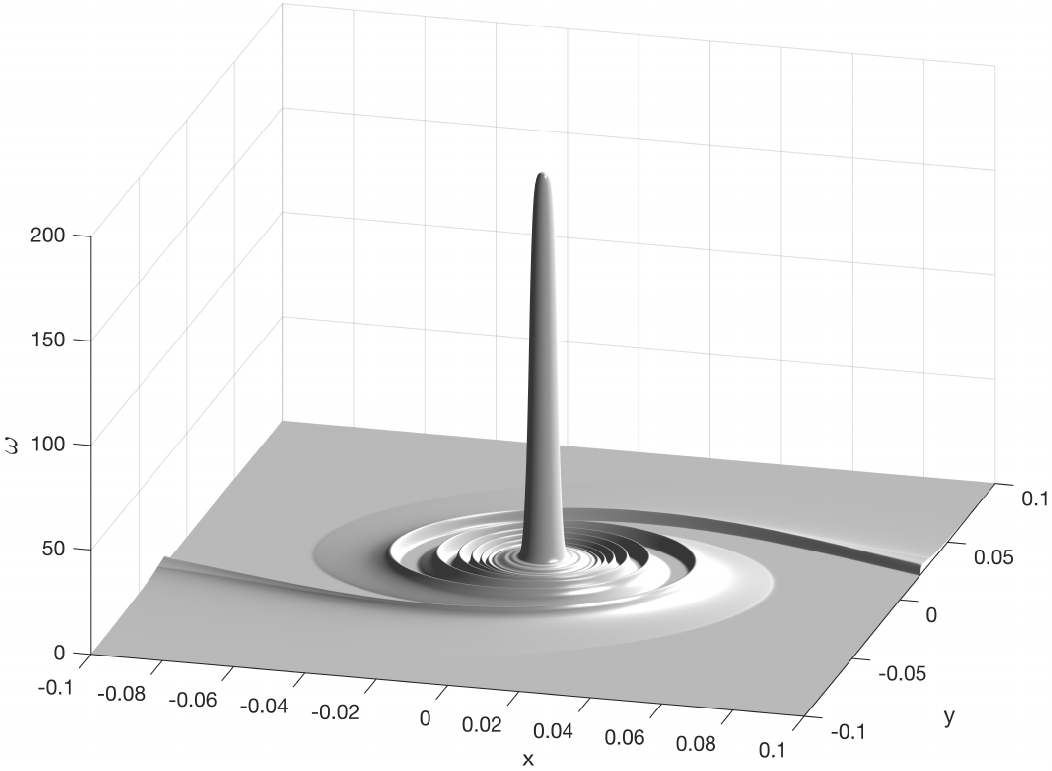}}
    \hbox{\quad (d) {\it $\omega$}-surface, $t=2.0$}}
  }
  \medskip
  \centerline{
    \vbox{\hbox{\includegraphics[height=0.22\textheight]{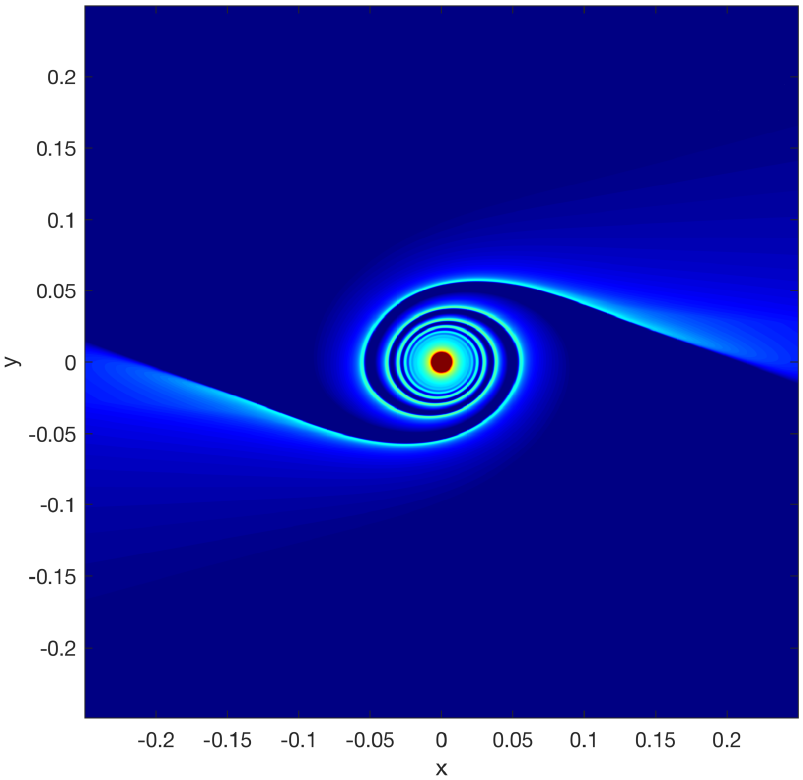}} 
    \hbox{\quad (e) {\it $\omega$}-contours, $t=3.0$}} 
    \qquad
    \vbox{\hbox{\includegraphics[height=0.25\textheight]{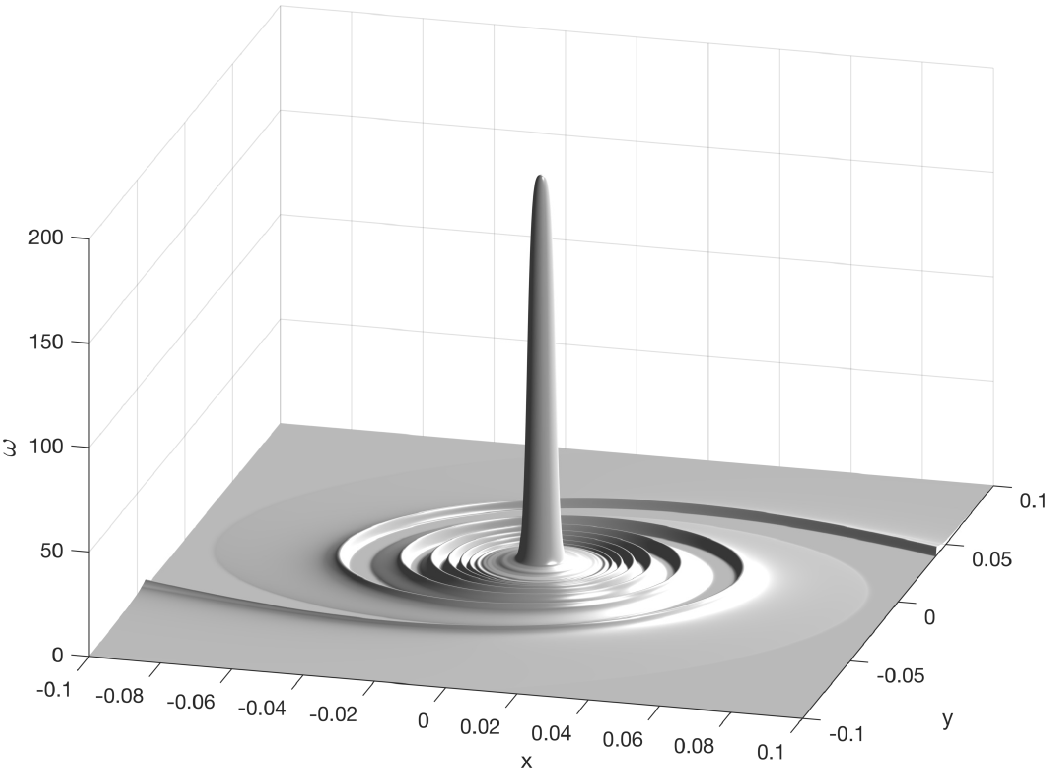}}
    \hbox{\quad (f) {\it $\omega$}-surface, $t=3.0$}}
  }
  \caption{Contour plots (left) and surface plots (right) of vorticity for solutions of Problem~1
           at increasing times $t$, showing the evolution of the one-spiral solution.  
           The contour plots display the grid whose extent is $[-0.25,0.25] \times [-0.25, 0.25]$ 
           containing $4096^2$ grid cells.  The surface plots illustrate details of the flow field
           including the vorticity magnitude; only the indicated portion of the innermost grid, 
           i.e.\ $\Omega_9$ in (\ref{eq:nesting}), is displayed.
           The solutions presented are for the case $(\alpha, \epsilon, \theta_w) = (0.95, 0.004, \pi/8)$.}
  \label{fig:ContSurf-P1}
\end{figure}

We first present a general impression of the evolution in time of solutions of Problem~1 and Problem~2,
taking $\epsilon=0.004$.  
In Fig.~\ref{fig:ContSurf-P1} below we display solutions of the first problem. 
Each row of this figure depicts the contour plot (left) and the graph (right) of the vorticity $\omega(x,y)$ at time $t=1, 2,$ or $3.$  
Graphs are displayed over only the indicated portion of the most-refined grid.
(Note that all plots presented in this work can be greatly enlarged 
without appreciable loss of resolution.)\ \ 
As shown, a single spiral winds itself in a counter-clockwise direction around the origin.  
As time progresses, the spiral grows outwards and remains centered at the origin.  

Figure~\ref{fig:ContSurf-P2} below depicts the time evolution of solutions of Problem 2, 
and is exactly analogous to Fig.~\ref{fig:ContSurf-P1}.  
The contour plots and graphs are displayed over the same regions, and all details of the computations, 
except for the initial data, are the same as those in Fig.~\ref{fig:ContSurf-P1}.  
The plots show that now two spirals form; they are symmetric about the origin and move apart along 
trajectories from $t=1$ to $t=3$.  More will be said about these trajectories later.
Each spiral consists of a counter-clockwise winding about a moving center.  
The fine structure of the moving spirals can be viewed in greater detail by enlarging the image.

\begin{figure}
  \centerline{
    \vbox{\hbox{\includegraphics[height=0.22\textheight]{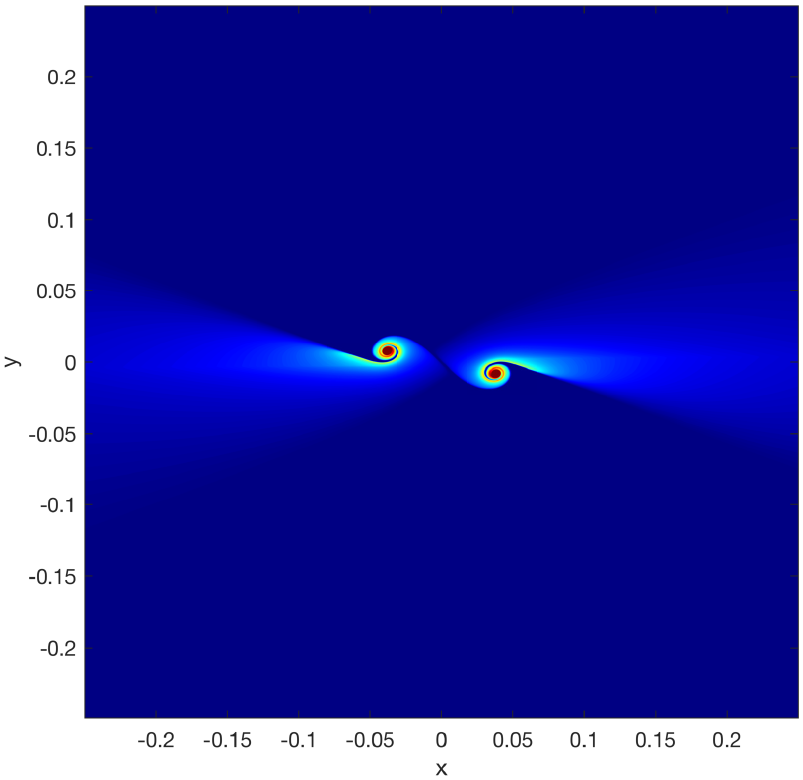}} 
    \hbox{\quad (a) {\it $\omega$}-contours, $t=1.0$}}
    \qquad
    \vbox{\hbox{\includegraphics[height=0.25\textheight]{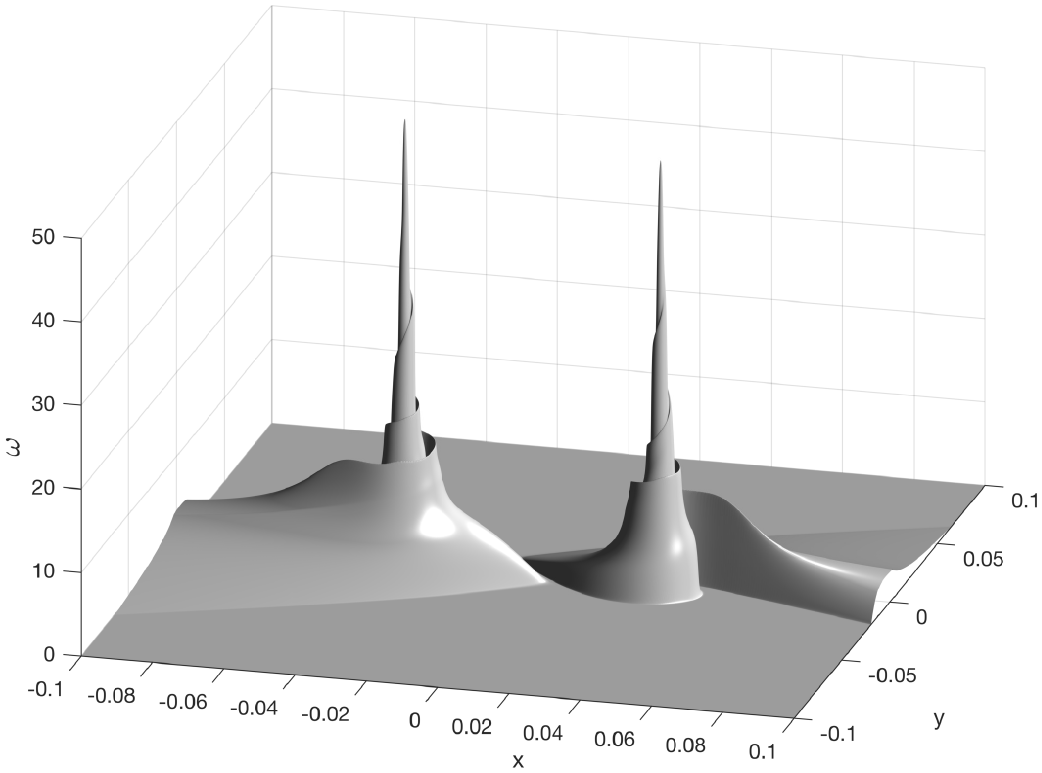}}
    \hbox{\quad (b) {\it $\omega$}-surface, $t=1.0$}}
  } 
  \medskip
  \centerline{
    \vbox{\hbox{\includegraphics[height=0.22\textheight]{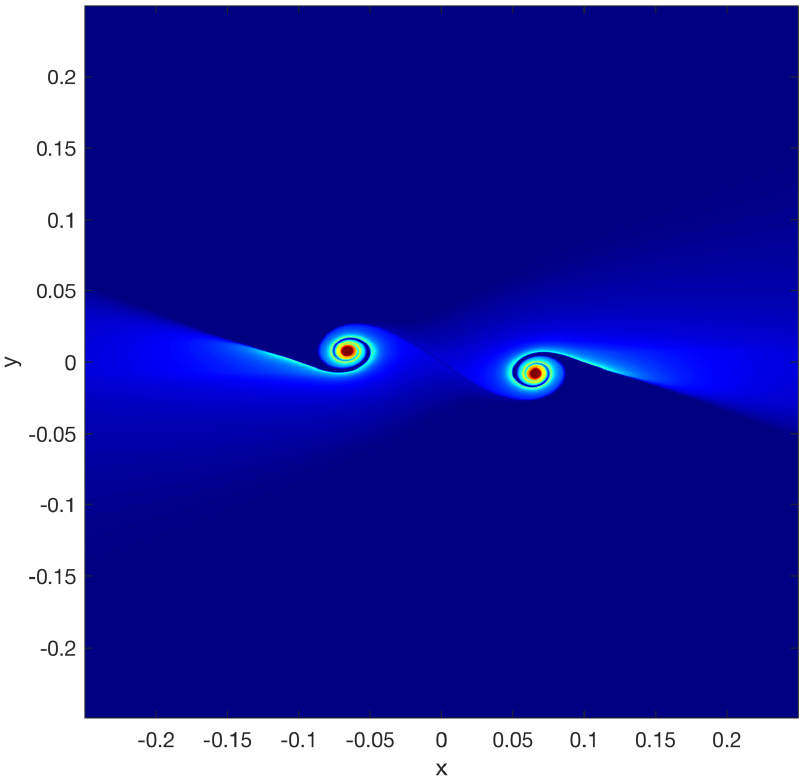}} 
    \hbox{\quad (c) {\it $\omega$}-contours, $t=2.0$}}
    \qquad
    \vbox{\hbox{\includegraphics[height=0.25\textheight]{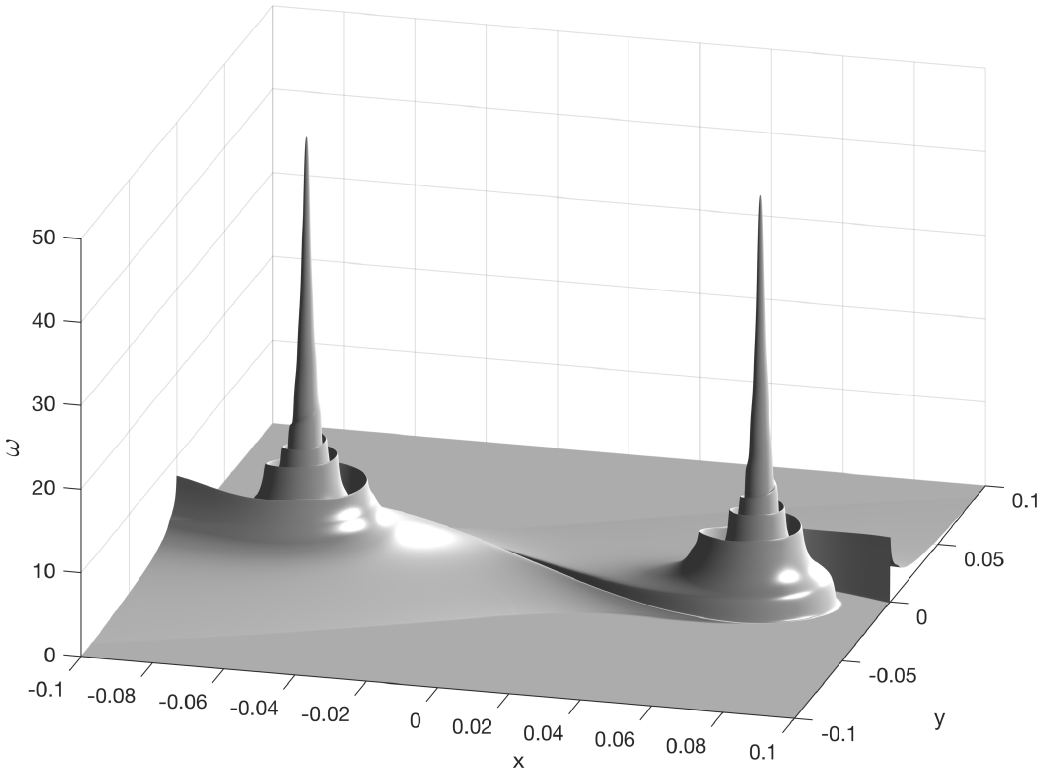}}
    \hbox{\quad (d) {\it $\omega$}-surface, $t=2.0$}}
  }
  \medskip
  \centerline{
    \vbox{\hbox{\includegraphics[height=0.22\textheight]{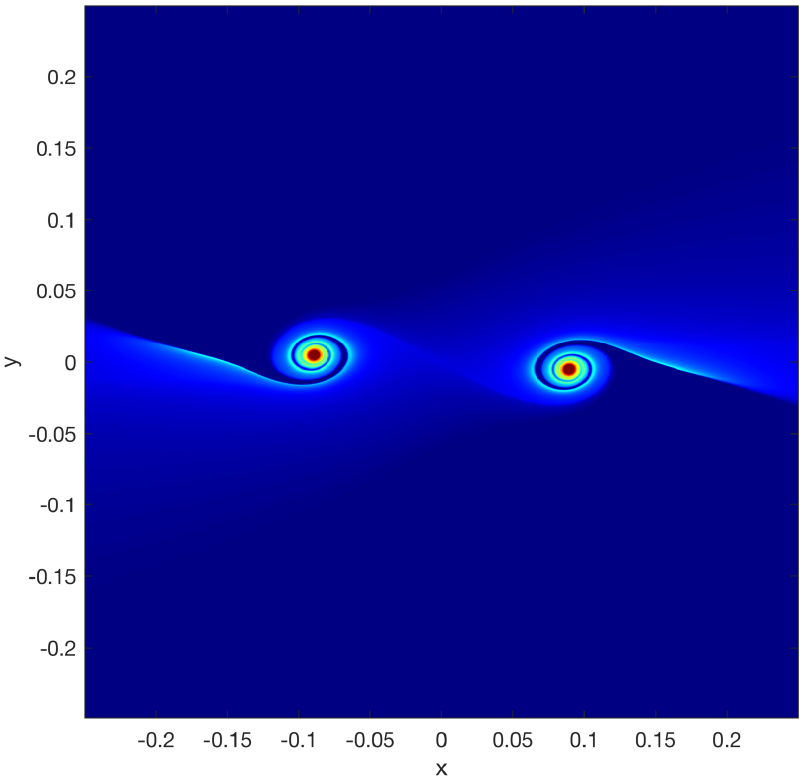}} 
    \hbox{\quad (e) {\it $\omega$}-contours, $t=3.0$}}
    \qquad
    \vbox{\hbox{\includegraphics[height=0.25\textheight]{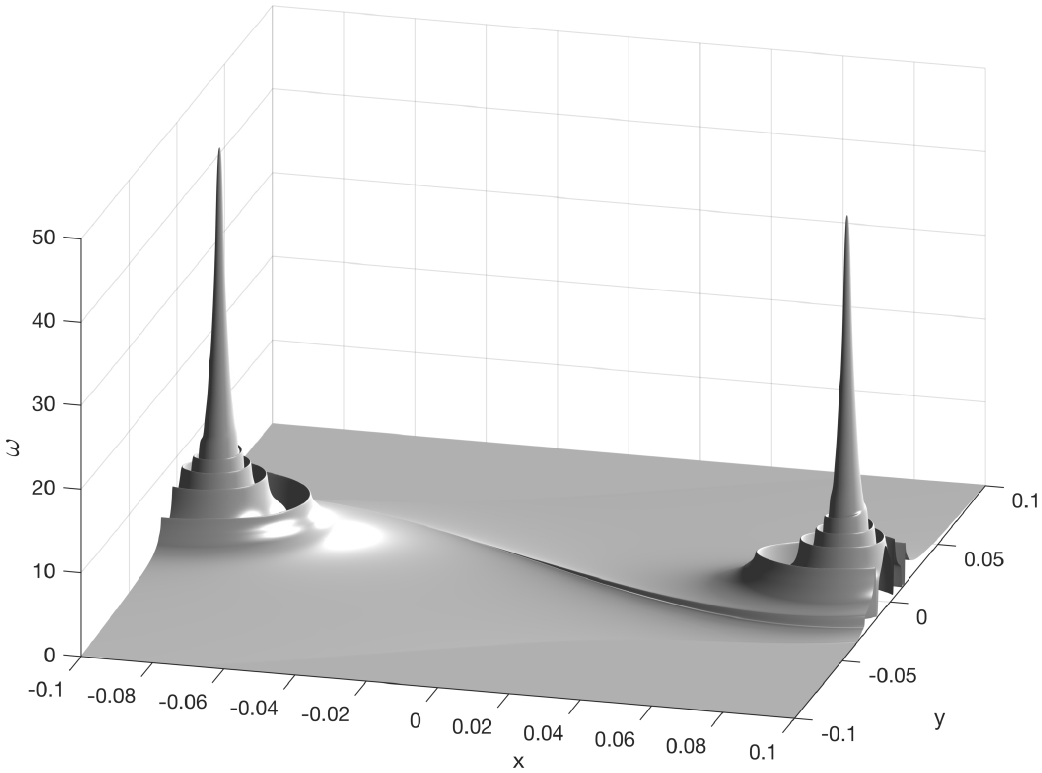}}
    \hbox{\quad (f) {\it $\omega$}-surface, $t=3.0$}}
  }
  \caption{Contour plots (left) and surface plots (right) of vorticity for solutions of Problem~2
           at increasing times $t$, showing the evolution of the two-spiral solution.  
           The contour plots display the grid whose extent is $[-0.25,0.25] \times [-0.25, 0.25]$ 
           containing $4096^2$ grid cells.  The surface plots illustrate details of the flow field
           including the vorticity magnitude; only the indicated portion of the innermost grid,
           i.e.\ $\Omega_9$ in (\ref{eq:nesting}), is displayed.
           The solutions presented are for the case $(\alpha, \epsilon, \theta_w) = (0.95, 0.004, \pi/8)$.}
  \label{fig:ContSurf-P2}
\end{figure}

\begin{figure}[!hp]
  \centerline{
    \vbox{\hbox to 6.5em{\vspace{0.065\textheight} \hbox{$\epsilon=0.004$\quad}} }   
    \vbox{\hbox{\includegraphics[height=0.12\textheight]{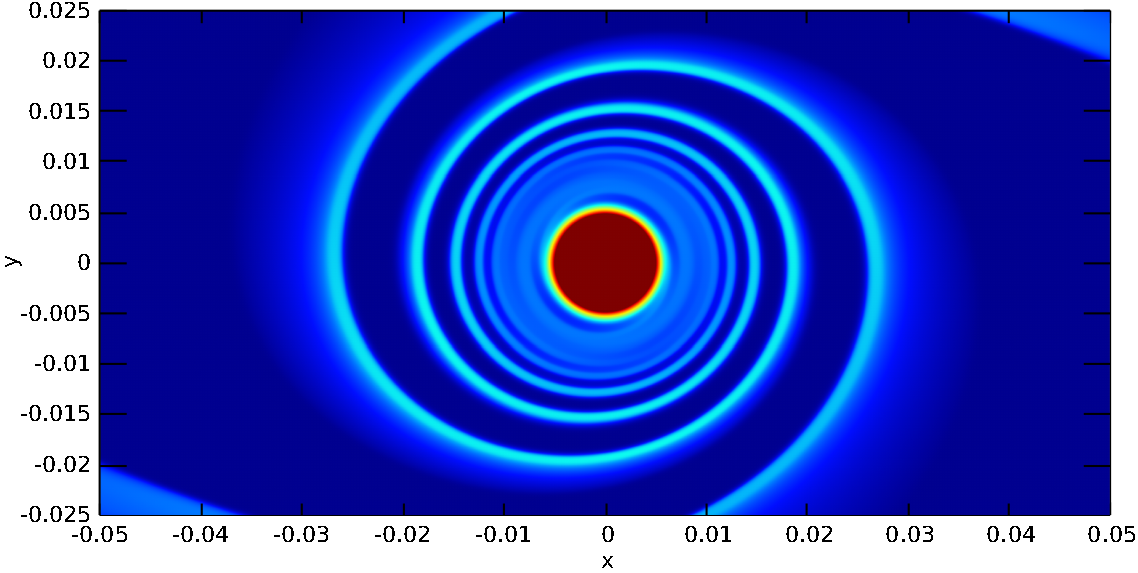}} 
    \vspace{-2ex} \hbox{\hspace{-2ex} (a)}}
    \quad
    \vbox{\hbox{\includegraphics[height=0.12\textheight]{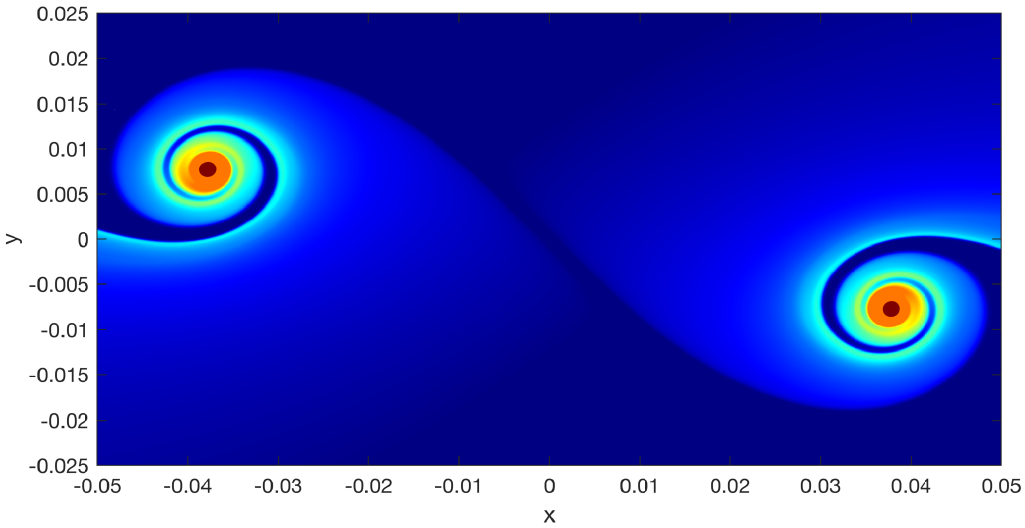}} 
    \vspace{-2ex} \hbox{\hspace{-2ex} (b)}}
  } \vspace{1ex}
  \centerline{
    \vbox{\hbox to 6.5em{\vspace{0.065\textheight} \hbox{$\epsilon=0.002$\quad}} }   
    \vbox{\hbox{\includegraphics[height=0.12\textheight]{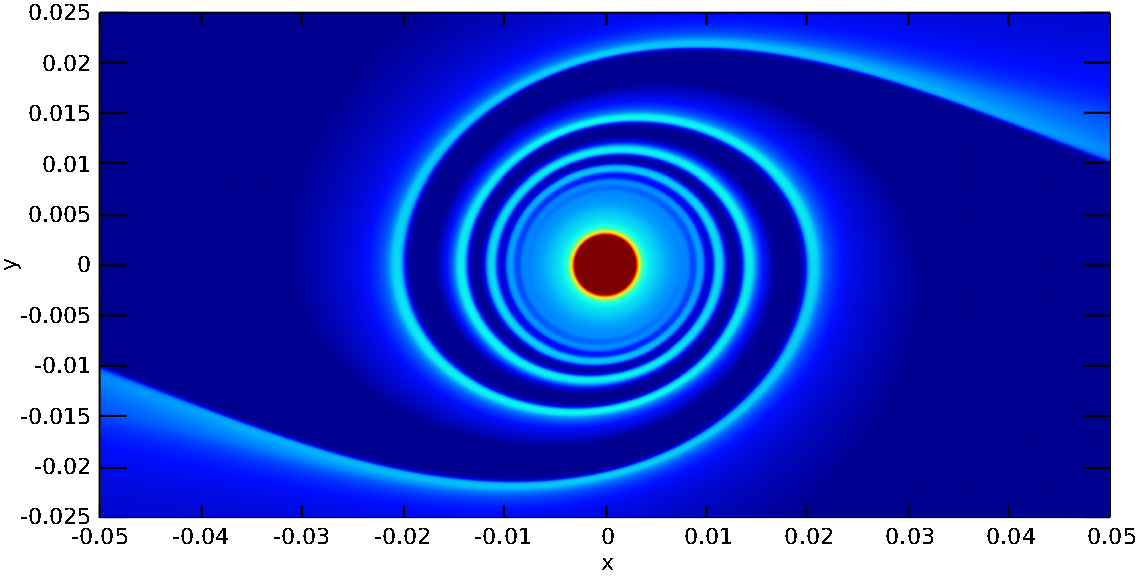}} 
    \vspace{-2ex} \hbox{\hspace{-2ex} (c)}}
    \quad
    \vbox{\hbox{\includegraphics[height=0.12\textheight]{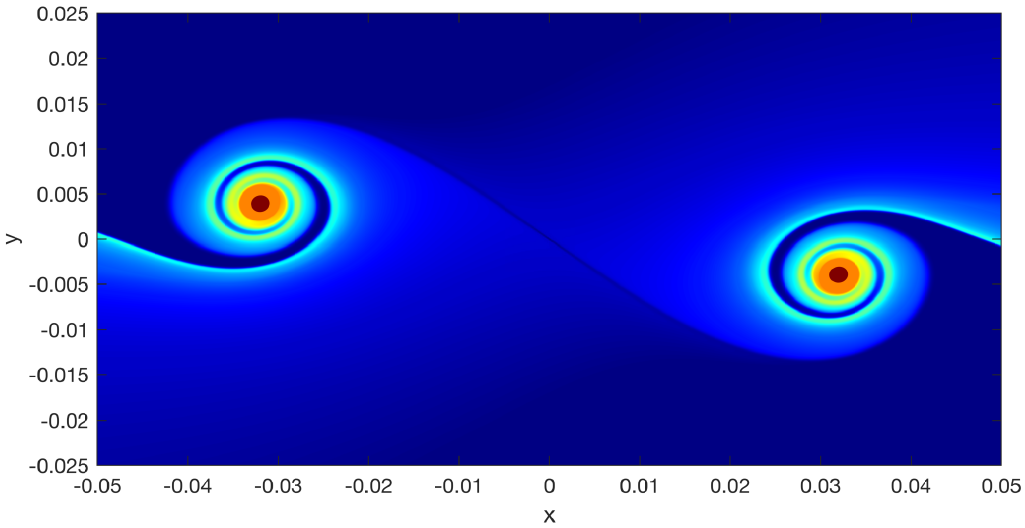}} 
    \vspace{-2ex} \hbox{\hspace{-2ex} (d)}}
  } \vspace{1ex}
  \centerline{
    \vbox{\hbox to 6.5em{\vspace{0.065\textheight} \hbox{$\epsilon=0.001$\quad}} }   
    \vbox{\hbox{\includegraphics[height=0.12\textheight]{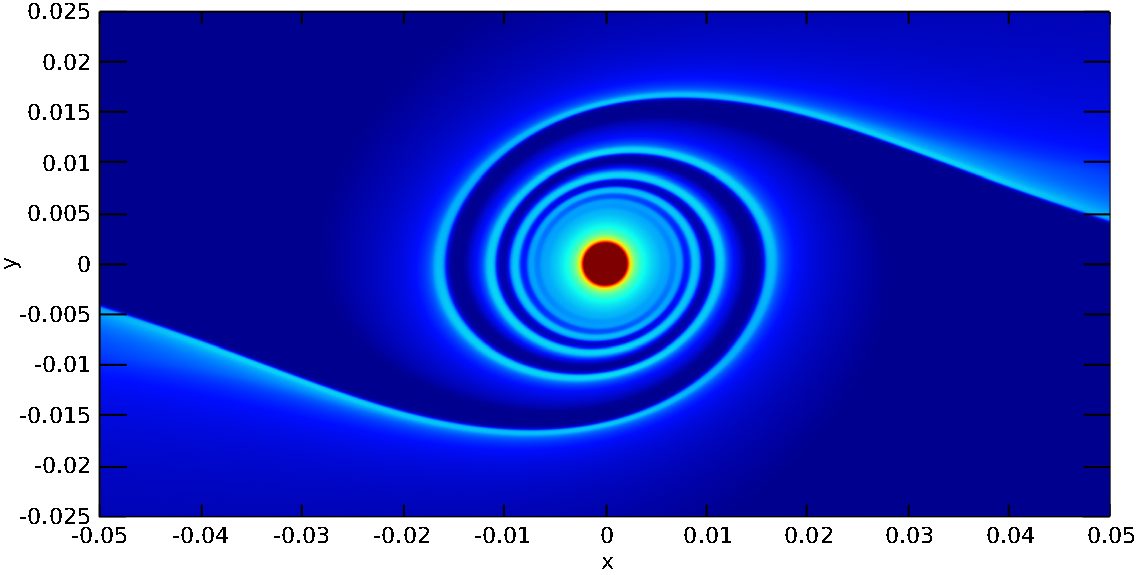}} 
    \vspace{-2ex} \hbox{\hspace{-2ex} (e)}}
    \quad
    \vbox{\hbox{\includegraphics[height=0.12\textheight]{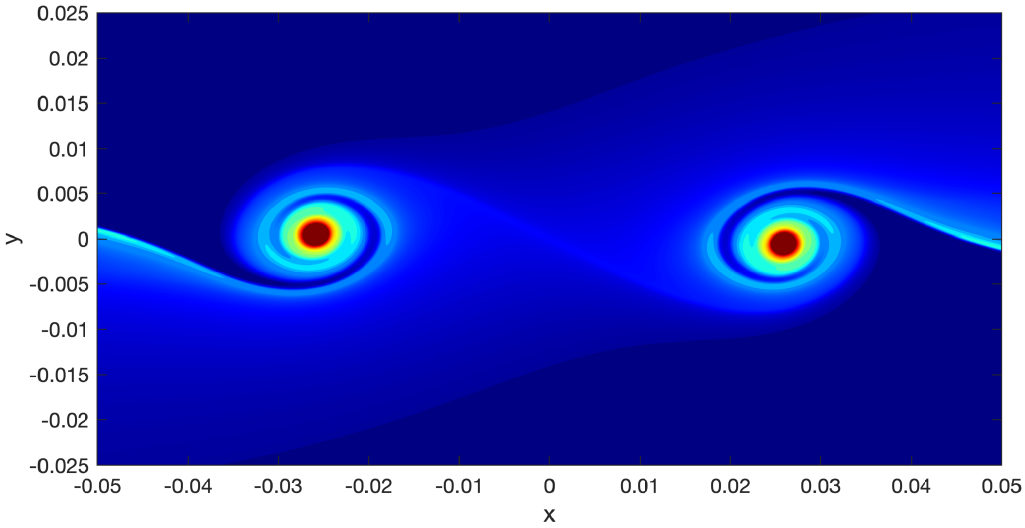}} 
    \vspace{-2ex} \hbox{\hspace{-2ex} (f)}}
  } \vspace{1ex}
  \centerline{
    \vbox{\hbox to 6.5em{\vspace{0.065\textheight} \hbox{$\epsilon=0.0005$\quad}} }   
    \vbox{\hbox{\includegraphics[height=0.12\textheight]{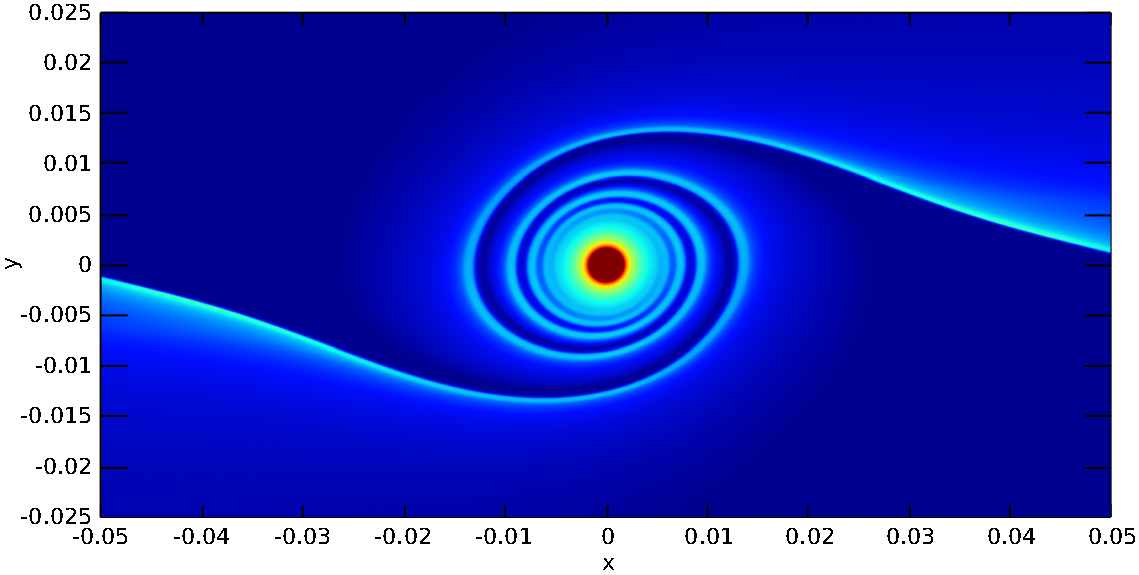}} 
    \vspace{-2ex} \hbox{\hspace{-2ex} (g)}}
    \quad
    \vbox{\hbox{\includegraphics[height=0.12\textheight]{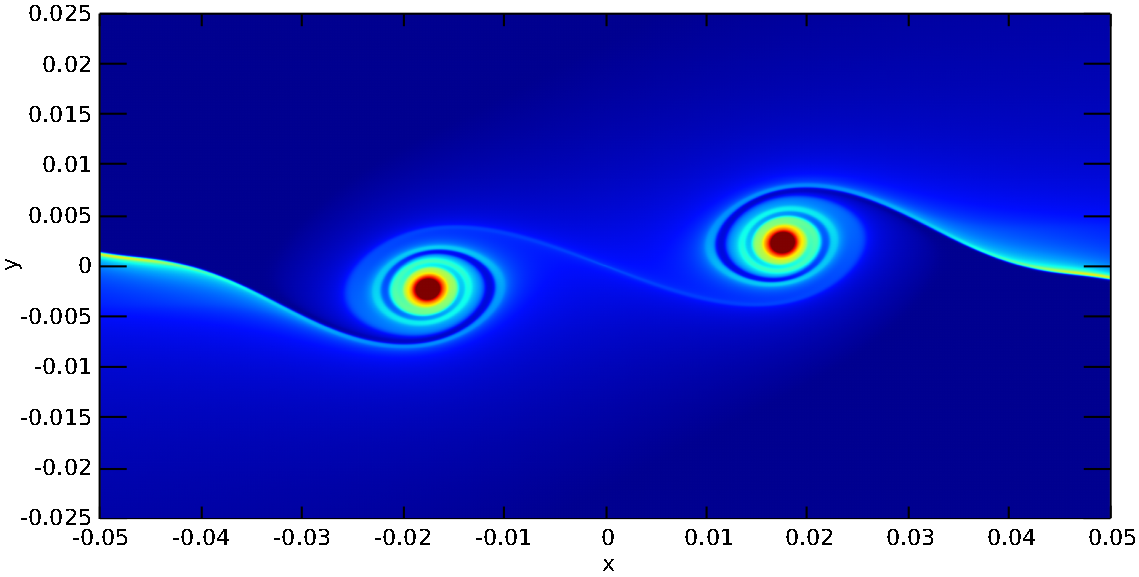}} 
    \vspace{-2ex} \hbox{\hspace{-2ex} (h)}}
  } \vspace{1ex}
  \centerline{
    \vbox{\hbox to 6.5em{\vspace{0.065\textheight} \hbox{$\epsilon=0.00025$\quad}} }   
    \vbox{\hbox{\includegraphics[height=0.12\textheight]{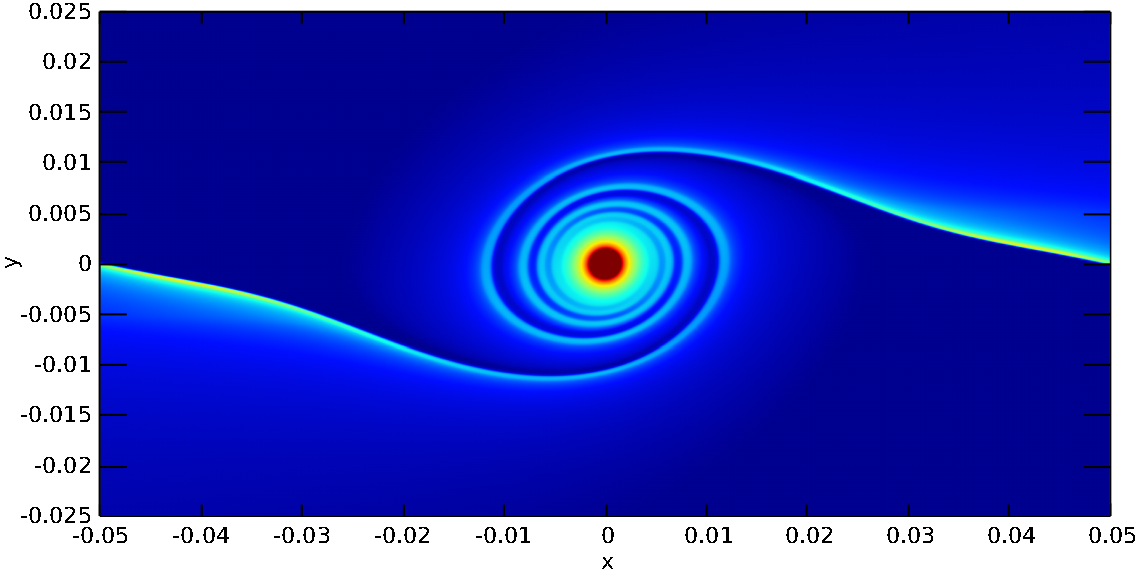}} 
    \vspace{-2ex} \hbox{\hspace{-2ex} (i)}}
    \quad
    \vbox{\hbox{\includegraphics[height=0.12\textheight]{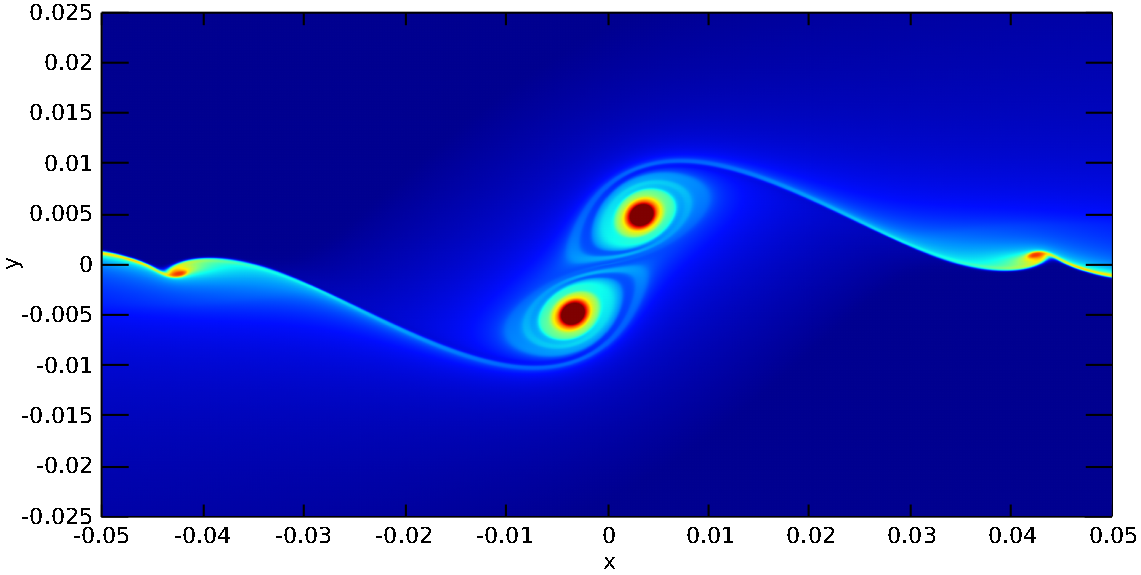}} 
    \vspace{-2ex} \hbox{\hspace{-2ex} (j)}}
  } \vspace{1ex}
  \centerline{
    \vbox{\hbox to 6.5em{\vspace{0.065\textheight} \hbox{$\epsilon=0.000125$\quad}} }   
    \vbox{\hbox{\includegraphics[height=0.12\textheight]{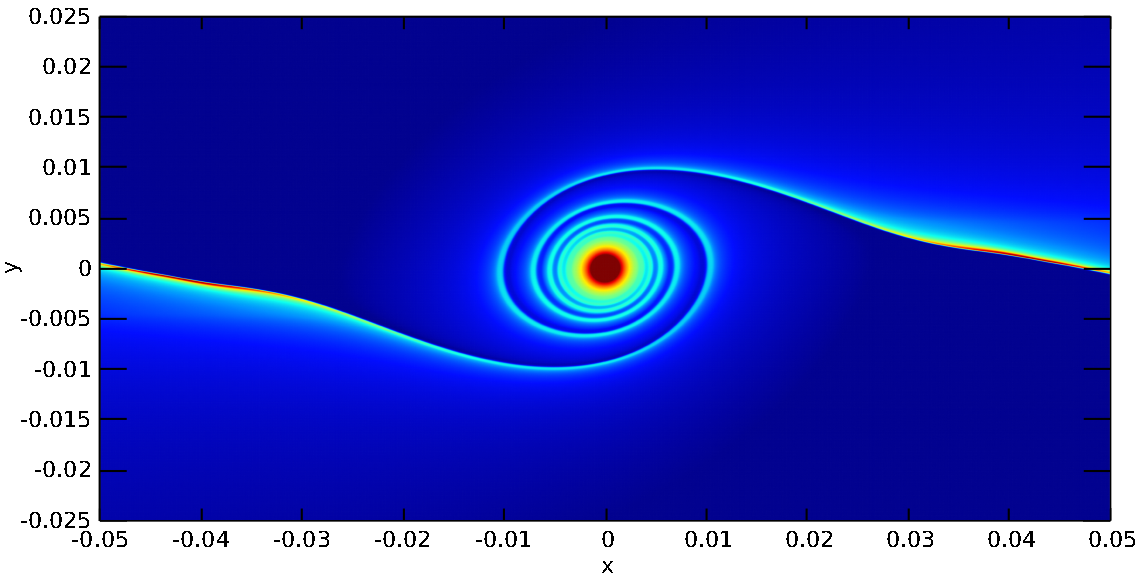}} 
    \vspace{-2ex} \hbox{\hspace{-2ex} (k)}}
    \quad
    \vbox{\hbox{\includegraphics[height=0.12\textheight]{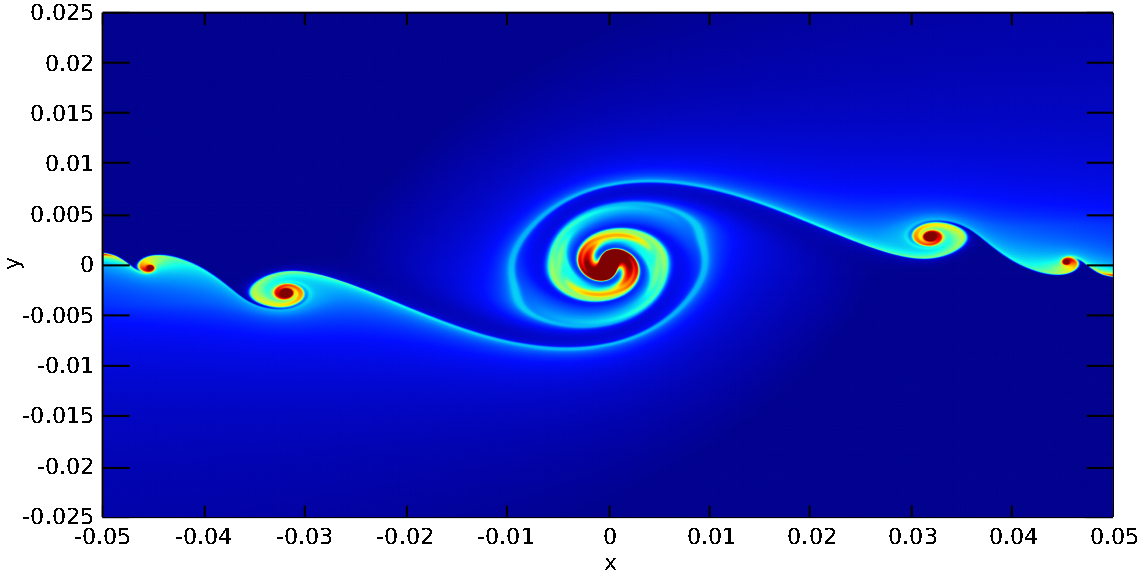}} 
    \vspace{-2ex} \hbox{\hspace{-2ex} (l)}}
  }
\caption{An illustration of the solution limit as $\epsilon$ tends to zero with $t$ fixed, by a sequence of vorticity contour plots; here $t=1$.  Plots on the left of each row are solutions of Problem~1, plots on the right are solutions of Problem~2, 
and both solutions are computed using the value of $\epsilon$ displayed next to the row.    
All plots display a region of the same physical size.
Each subgrid in the computations contains $4096^2$ cells.  
For $0.00025 \le \epsilon \le 0.004$ the 
innermost nested domain terminates at $\Omega_9$; see (\ref{eq:nesting}).  However for $\epsilon = 0.000125$ 
the innermost nested domain terminates at $\Omega_{10}$.
}
  \label{fig:EpsTo0-t=1}
\end{figure}

\begin{figure}[!hp]
  \centerline{
    \vbox{\hbox to 6.5em{\vspace{0.065\textheight} \hbox{$\epsilon=0.004$\quad}} }   
    \vbox{\hbox{\includegraphics[height=0.12\textheight]{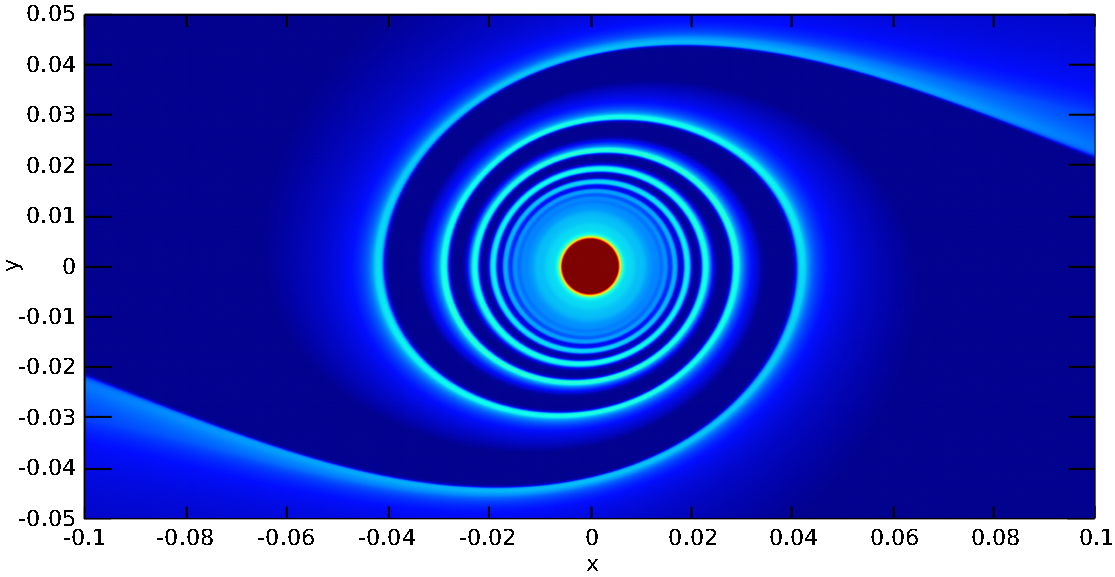}} 
    \vspace{-2ex} \hbox{\hspace{-2ex} (a)}}
    \quad
    \vbox{\hbox{\includegraphics[height=0.12\textheight]{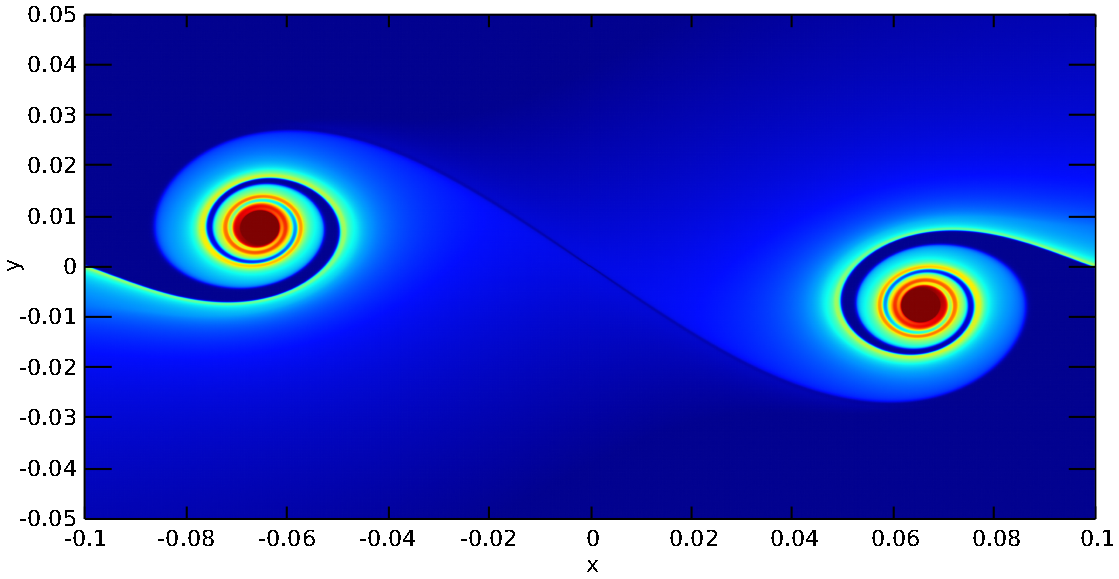}} 
    \vspace{-2ex} \hbox{\hspace{-2ex} (b)}}
  } \vspace{1ex}
  \centerline{
    \vbox{\hbox to 6.5em{\vspace{0.065\textheight} \hbox{$\epsilon=0.002$\quad}} }   
    \vbox{\hbox{\includegraphics[height=0.12\textheight]{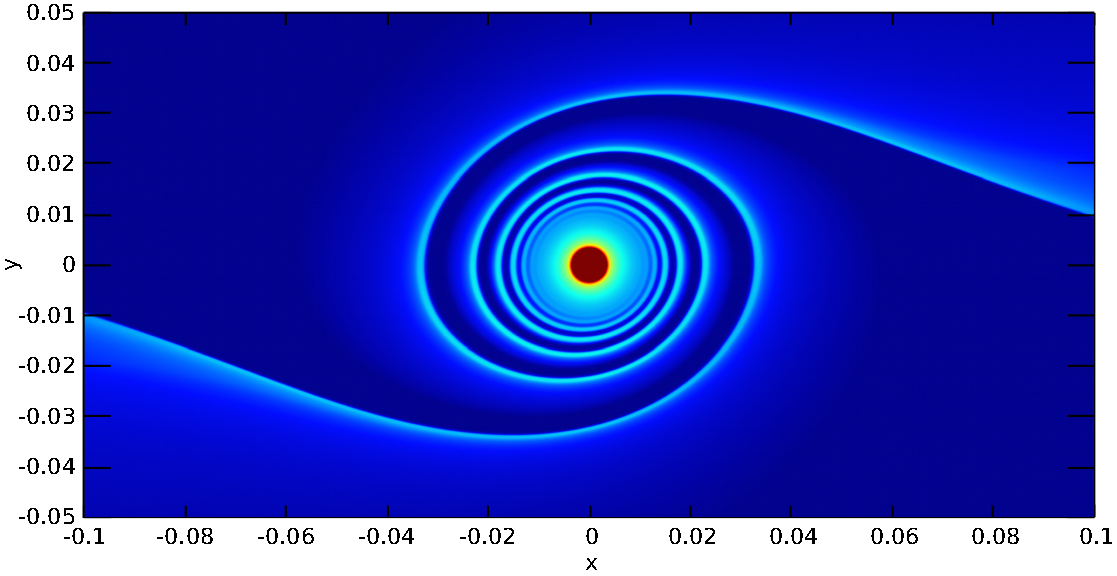}} 
    \vspace{-2ex} \hbox{\hspace{-2ex} (c)}}
    \quad
    \vbox{\hbox{\includegraphics[height=0.12\textheight]{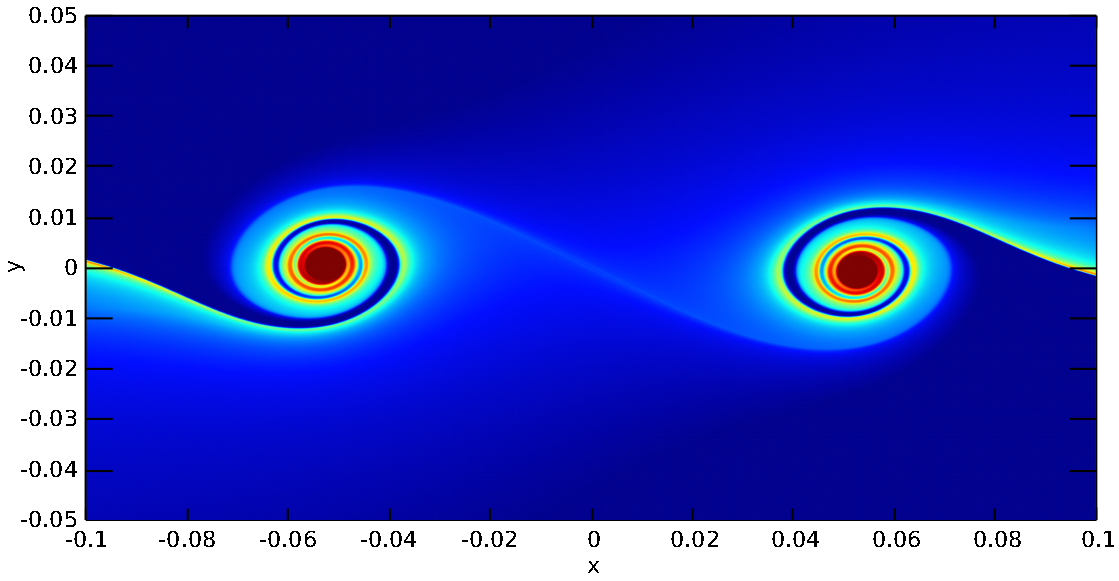}} 
    \vspace{-2ex} \hbox{\hspace{-2ex} (d)}}
  } \vspace{1ex}
  \centerline{
    \vbox{\hbox to 6.5em{\vspace{0.065\textheight} \hbox{$\epsilon=0.001$\quad}} }   
    \vbox{\hbox{\includegraphics[height=0.12\textheight]{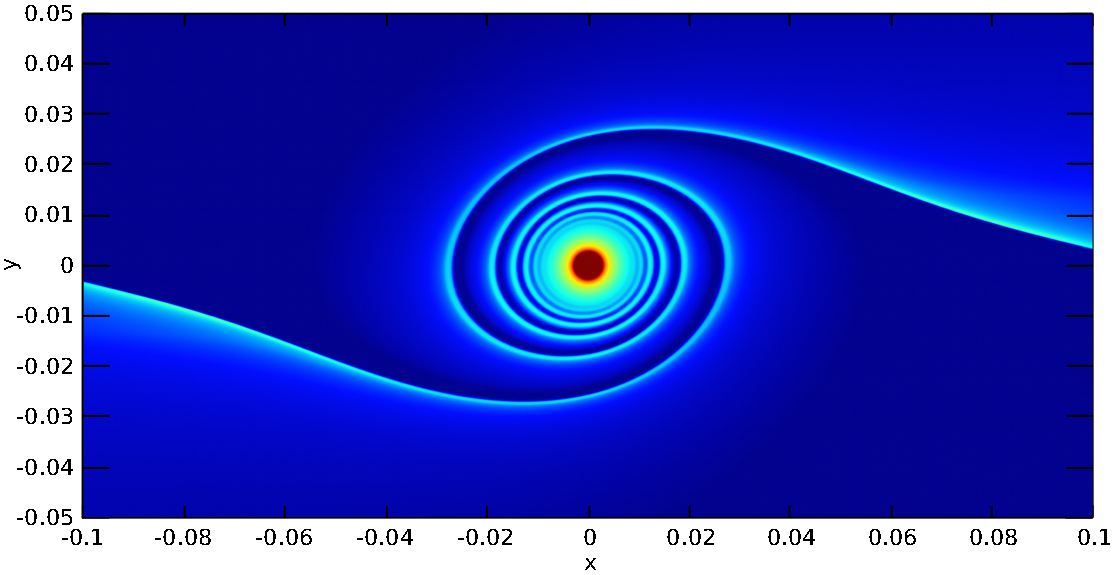}} 
    \vspace{-2ex} \hbox{\hspace{-2ex} (e)}}
    \quad
    \vbox{\hbox{\includegraphics[height=0.12\textheight]{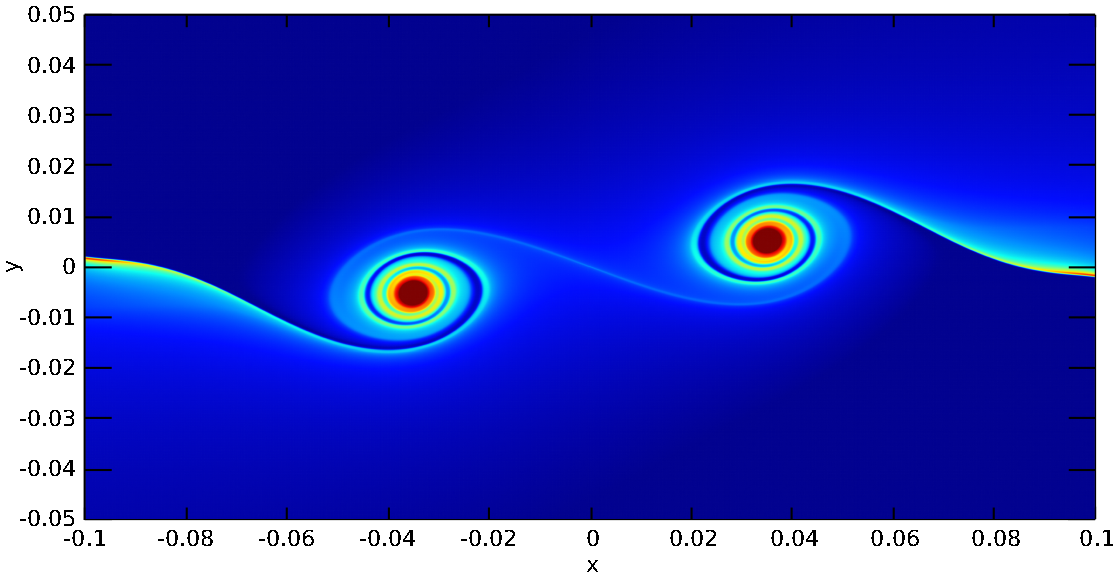}} 
    \vspace{-2ex} \hbox{\hspace{-2ex} (f)}}
  } \vspace{1ex}
  \centerline{
    \vbox{\hbox to 6.5em{\vspace{0.065\textheight} \hbox{$\epsilon=0.0005$\quad}} }   
    \vbox{\hbox{\includegraphics[height=0.12\textheight]{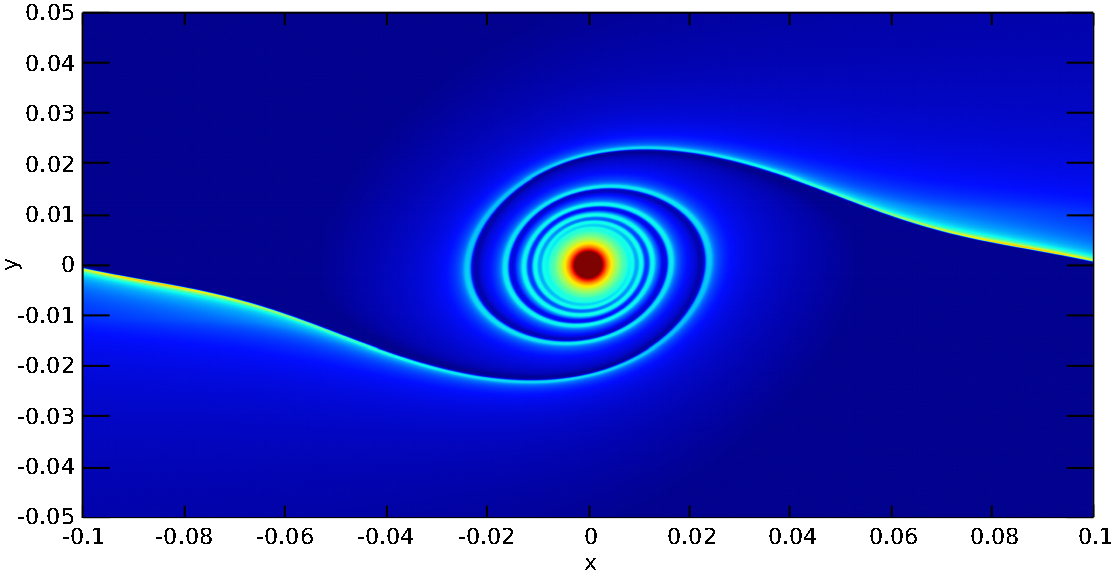}} 
    \vspace{-2ex} \hbox{\hspace{-2ex} (g)}}
    \quad
    \vbox{\hbox{\includegraphics[height=0.12\textheight]{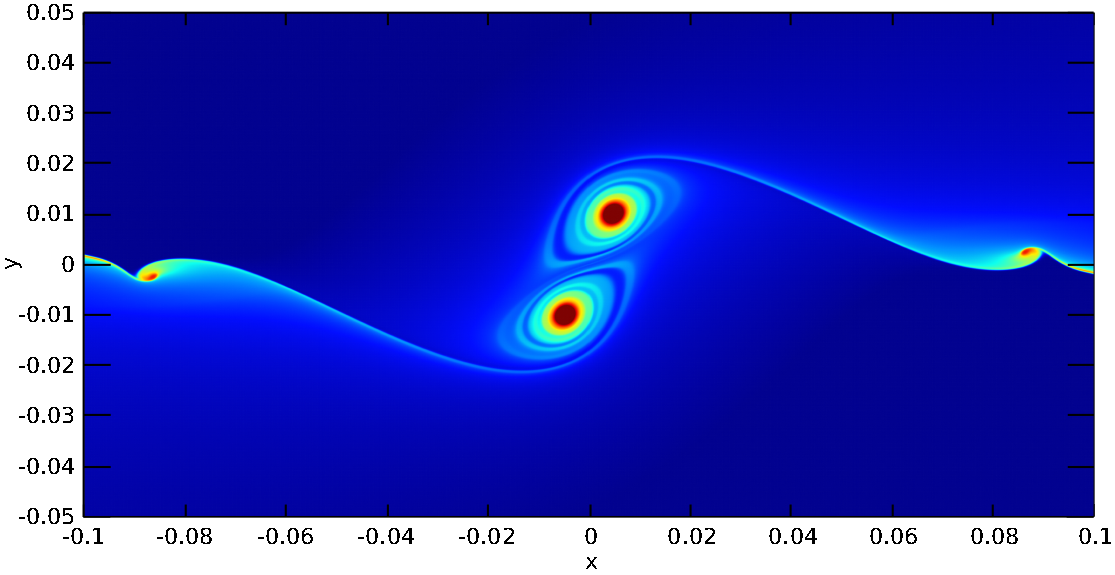}} 
    \vspace{-2ex} \hbox{\hspace{-2ex} (h)}}
  } \vspace{1ex}
  \centerline{
    \vbox{\hbox to 6.5em{\vspace{0.065\textheight} \hbox{$\epsilon=0.00025$\quad}} }   
    \vbox{\hbox{\includegraphics[height=0.12\textheight]{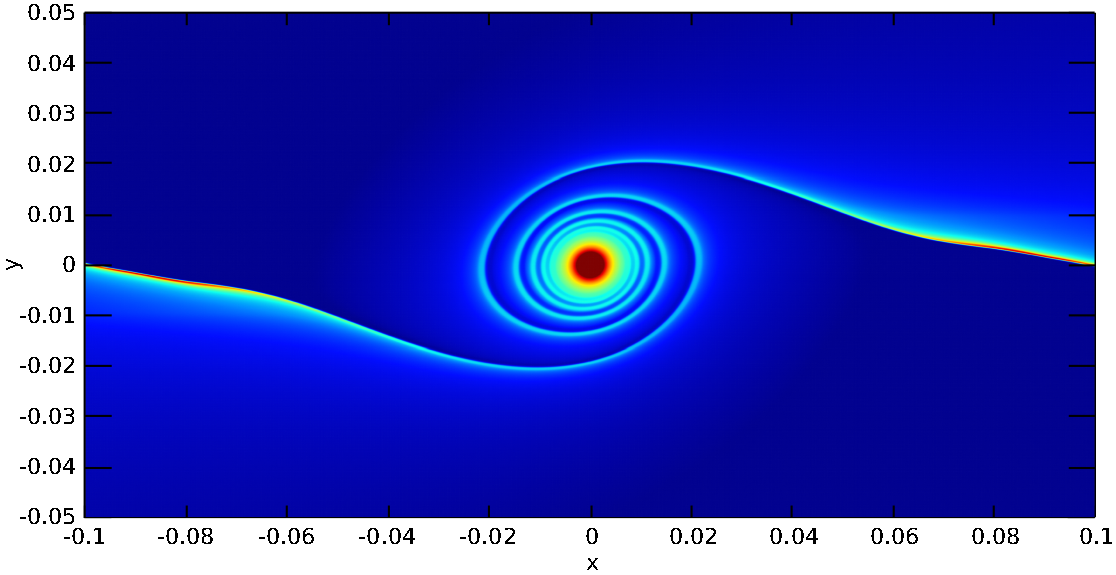}} 
    \vspace{-2ex} \hbox{\hspace{-2ex} (i)}}
    \quad
    \vbox{\hbox{\includegraphics[height=0.12\textheight]{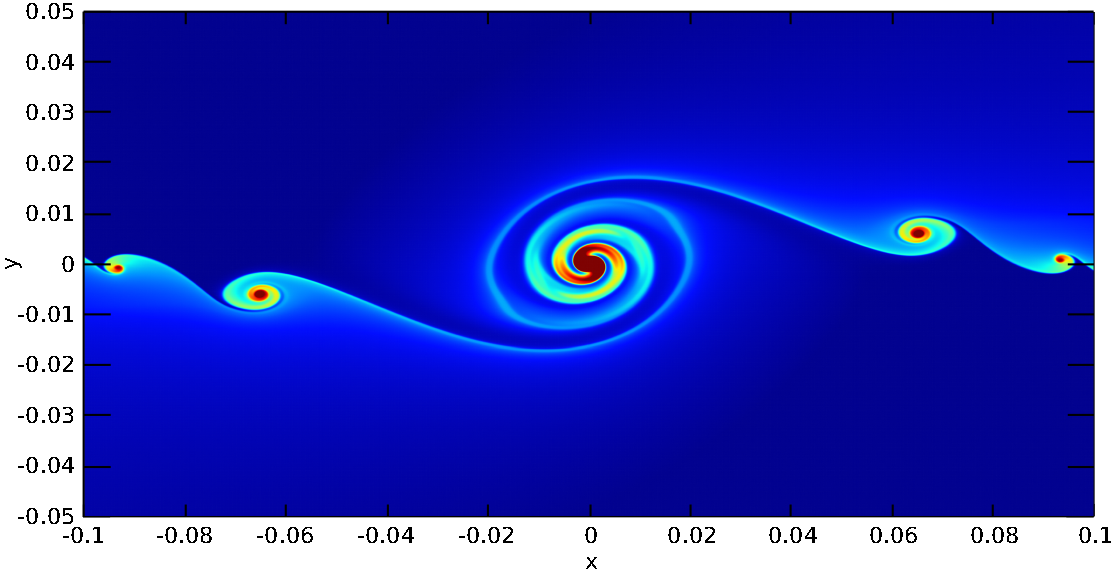}} 
    \vspace{-2ex} \hbox{\hspace{-2ex} (j)}}
  } \vspace{1ex}
  \centerline{
    \vbox{\hbox to 6.5em{\vspace{0.065\textheight} \hbox{$\epsilon=0.000125$\quad}} }   
    \vbox{\hbox{\includegraphics[height=0.12\textheight]{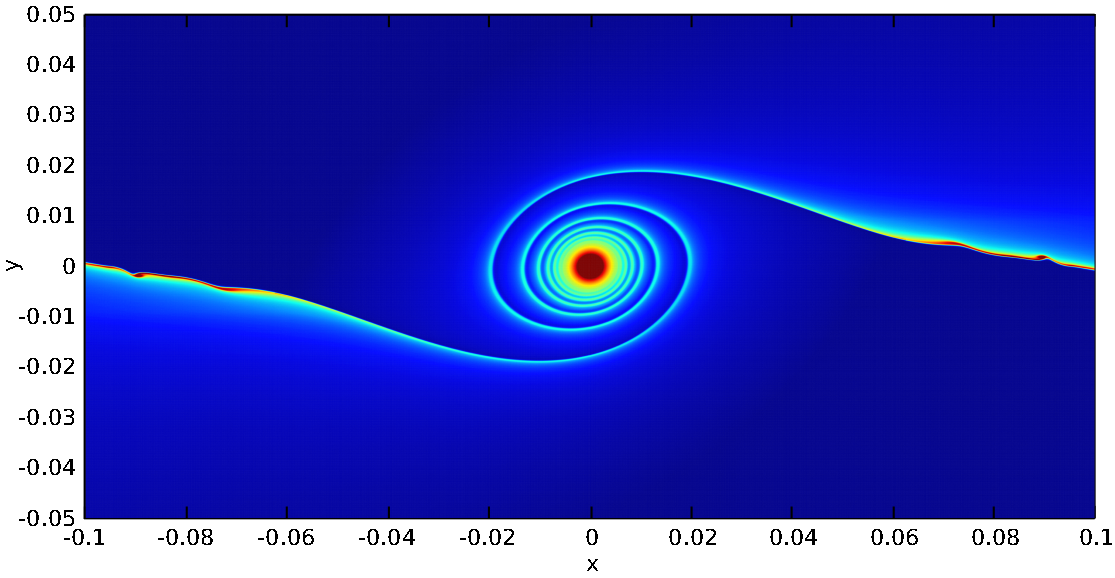}} 
    \vspace{-2ex} \hbox{\hspace{-2ex} (k)}}
    \quad
    \vbox{\hbox{\includegraphics[height=0.12\textheight]{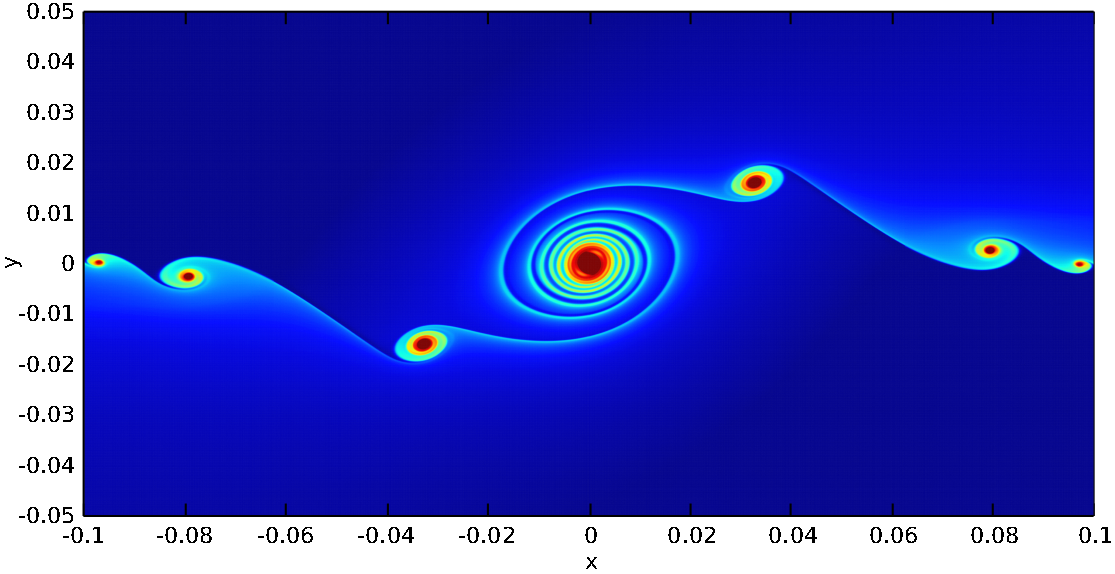}} 
    \vspace{-2ex} 
    \hbox{\hspace{-2ex} (l)}}
  }
\caption{An illustration of the solution limit as $\epsilon$ tends to zero with $t$ fixed, by a sequence of vorticity contour plots; here $t=2$.  Plots on the left of each row are solutions of Problem~1, plots on the right are solutions of Problem~2, 
and both solutions are computed using the value of $\epsilon$ displayed next to the row.    
All plots display a region of the same physical extent; 
note, however, now the region displayed is twice as large as given in Fig.~\ref{fig:EpsTo0-t=1} due
to the temporal movement of the twin spiral centers.
Each subgrid in the computations contains $4096^2$ cells.  
For $\epsilon$ in the range $0.00025$ through $0.004$ the 
innermost nested domain terminates at $\Omega_9$; see (\ref{eq:nesting}).  However for $\epsilon = 0.000125$ 
the innermost nested domain terminates at $\Omega_{10}$.
}
  \label{fig:EpsTo0-t=2}
\end{figure}

In Figs.~\ref{fig:EpsTo0-t=1} and \ref{fig:EpsTo0-t=2} we now come to the main question:  
what happens when $t$ is fixed and $\epsilon$ tends to zero?  
We will take the decreasing sequence of $\epsilon$ given in (\ref{eq:parameter-set}).
In determining this behavior numerically, we are careful to maintain an appropriate ratio 
of $\epsilon/\Delta{x}$, so that the disk of radius $\epsilon$ which is, initially, either filled with positive vorticity (Problem~1) 
or has identically zero vorticity (Problem~2) is well captured.  
In Figs.~\ref{fig:EpsTo0-t=1} and \ref{fig:EpsTo0-t=2} which follow, we maintain a ratio $\epsilon/\Delta{x}$ of at least 4.096 
(providing 8.192 grid cells along a diameter of the disk).  
We find that this level of refinement provides adequately resolved computed solutions.  
Maintaining this ratio requires an extreme level of grid refinement, particularly for the smaller values of $\epsilon$.
In particular, for $\epsilon = 0.000125$ and a uniform grid of size $4096\times4096$, we utilize nested domains 
terminating at
\[
	\Omega_{10} = [-0.0625,0.0625]\times[-0.0625,0.0625] \ \ \Rightarrow \ \ \Delta{x} \approx 3.05 \cdot 10^{-5},
\]
yielding the desired ratio.  

All solutions in Fig.~\ref{fig:EpsTo0-t=1} below are presented at $t=1.0$.
Each row shows contour plots of solutions of Problem~1 (left) and Problem~2 (right), at the value of $\epsilon$ indicated  
to the row's left.  
The physical extent of the region shown is the same for all plots in the figure.  
The nested domains terminate at 
\[ 
	\Omega_9 = [-0.125,0.125]\times[-0.125,0.125]
\] 
for solutions depicted in plots (a)--(j) (i.e.\ $\epsilon \ge 0.00025$), and 
\[
	\Omega_{10} = [-0.0625,0.0625]\times[-0.0625,0.0625]
\]
for those in plots (k)--(l) ($\epsilon = 0.000125$).
Starting with the first row and moving downwards,
we see a progression in which the single spirals in the Problem~1 solutions become smaller and more compact, 
while the twin spirals in the Problem~2 solutions move closer together, as $\epsilon$ becomes smaller.  
By the time the fifth row is reached, $\epsilon = 0.00025$, 
the twin spirals have moved into close proximity and are beginning to wind together.
In the final row, $\epsilon = 0.000125$, the twin spirals have completely joined to become a single spiral.
Small regions 
of locally higher vorticity are visible in the tails of both the left and right solutions for $\epsilon = 0.00025$ and $\epsilon = 0.000125$.

In Fig.~\ref{fig:EpsTo0-t=2} all computational parameters are identical to those discussed in the
previous paragraph, except now solutions are integrated to $t=2.0$.
The physical extent displayed in all plots, however, is larger than presented in Fig.~\ref{fig:EpsTo0-t=1}
due to the additional motion of the double spirals from $t=1$ to $t=2$.   
Moving downwards in the figure a row at a time, we again see a progression
in which two-spiral solutions of Problem~2 merge into a one-spiral solution.
In the fourth row, $\epsilon = 0.0005$, the twin spirals have begun to wind onto each other; beginning with this row, the
Problem~1 solutions on the left and the Problem~2 solutions on the right of each row bear a mutual resemblance.
By the time we reach the fifth row, $\epsilon = 0.00025$, the twin spirals have essentially merged.   
In the last row, $\epsilon = 0.000125$, the Problem~1 and Problem~2 
solutions each contain a single spiral and appear similar, 
with small regions of locally higher vorticity present in the spiral tails of both.

It is interesting to note that in plots Fig.~\ref{fig:EpsTo0-t=2} (k) and (l), the displayed view 
overlaps the interface between nested domains $\Omega_9$ and $\Omega_{10}$.  
The vertical portion of this interface, located at $x = \pm 0.0625$, is visible. 
We encourage the reader to enlarge plots (k) and (l) near these $x$-values to confirm there are
no signs of undesirable artifacts coming from our numerical treatment of sub-domain boundary interfaces.

\begin{figure}[h]
  \centerline{\includegraphics[width=0.75\textwidth]{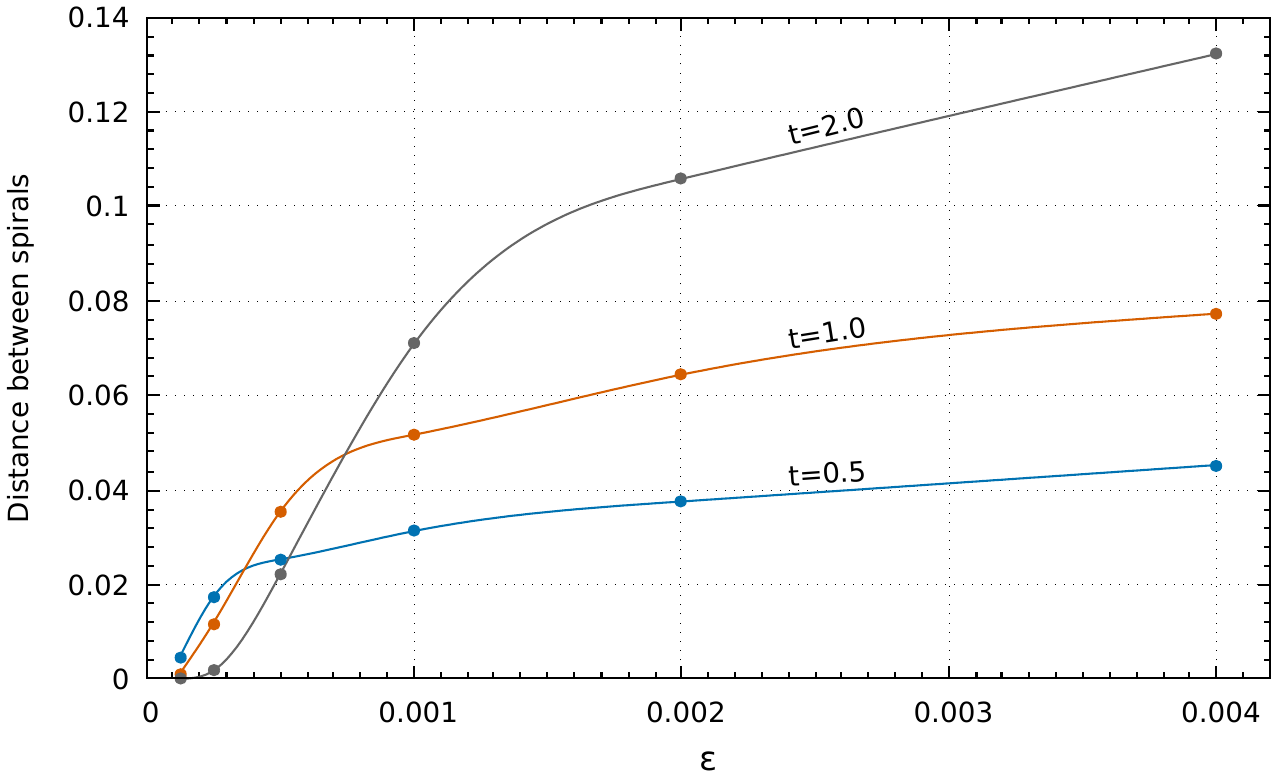}}
  \caption{An illustration of the transition of two-spiral solutions to one-spiral solutions 
  for fixed $t$ as $\epsilon$ tends to zero.  
  Marked data points in the figure above represent the distance between twin spiral centers
  depicted in Fig.~\ref{fig:EpsTo0-t=1} and Fig.~\ref{fig:EpsTo0-t=2}.  
  Earlier time data was used to generate the curved labeled t~=~0.5.  
  } \label{fig:Rate}
\end{figure}

Figures~\ref{fig:EpsTo0-t=1} and \ref{fig:EpsTo0-t=2} imply that solutions of Problem~2
for fixed $t>0$ converge to a one-spiral solution as $\epsilon$ tends to zero.
To illustrate the $\epsilon$-convergence graphically, in Fig.~\ref{fig:Rate} 
we plot the distance between spiral centers 
in solutions of Problem~2 as a function of $\epsilon$ for certain values of fixed $t$.  
Here and throughout, a spiral center is regarded as the spatial location of the maximum flow-field vorticity.
The computed data points in this figure are connected with smooth cubic spline curves.  
Marked data points on the curve labeled  t~=~1.0 correspond to solutions of Problem~2 depicted on the right in Fig.~\ref{fig:EpsTo0-t=1}; 
marked points on the curve labeled t~=~2.0 correspond to solutions depicted on the right in Fig.~\ref{fig:EpsTo0-t=2}.  
To give a more complete picture of how these curves depend on fixed values of $t$, 
an additional curve at t~=~0.5 is introduced.
Another way to visualize the $\epsilon$-dependent spiral center distance for fixed $t$, seen here, 
is provided in Fig.~\ref{fig:EpsTo0} below.  
There, an alternate interpolation technique is employed to display
spiral center {\it locations} as functions of $\epsilon$ for a variety of fixed $t$'s.

\begin{figure}[!h]
  \centerline{\includegraphics[width=0.7\textwidth]{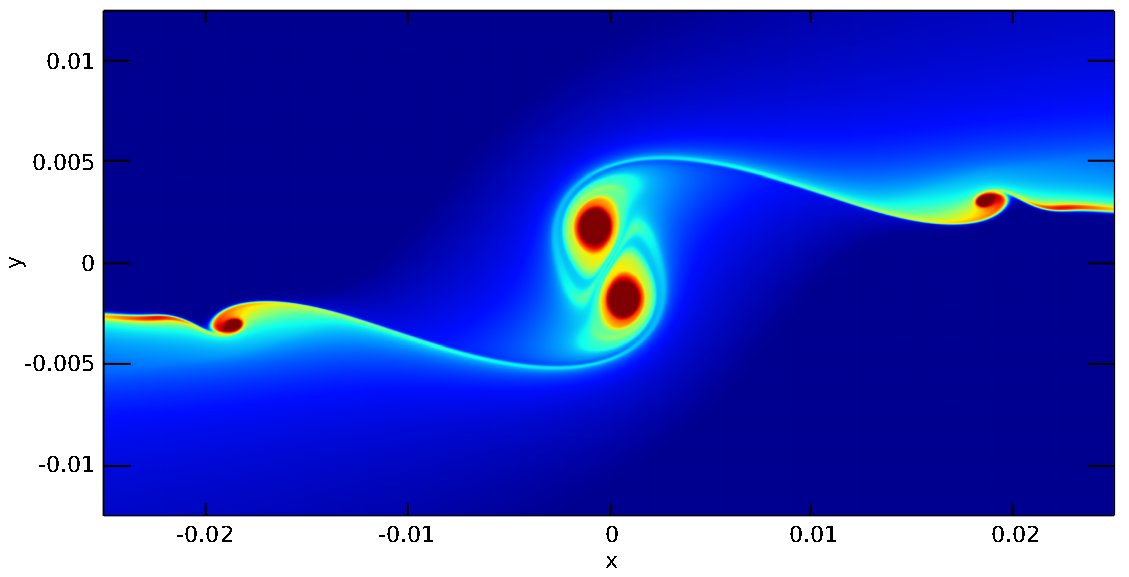}}
  \caption{A plot of vorticity contours for the solution of Problem 2 at 
  $t=0.5, \, \epsilon = 0.000125$,
  corresponding to the leftmost data point on the blue curve in Figure~\ref{fig:Rate} above.  
  Two spirals are visible, but they have almost completely merged to form one spiral.  
  } \label{fig:t=0.5_eps=0.001-Detail}
\end{figure}

The leftmost data point on the blue curve in Fig.~\ref{fig:Rate} indicates that when $t=0.5$ the spirals are still a 
small distance apart even in the solution with \hbox{$\epsilon=0.000125$}.  
We present a vorticity contour plot for the solution corresponding to this data point 
in Fig.~\ref{fig:t=0.5_eps=0.001-Detail}.  
As shown, the twin spirals have begun to intertwine and are close to merging.
Further reduction in $\epsilon$ would certainly  
result in the twin spirals fully merging into a single spiral;
see Fig.~\ref{fig:EpsTo0-t=1} (j,l) and Fig.~\ref{fig:EpsTo0-t=2} (h,j).

\begin{figure}[h]
  \centerline{\includegraphics[width=1.0\textwidth]{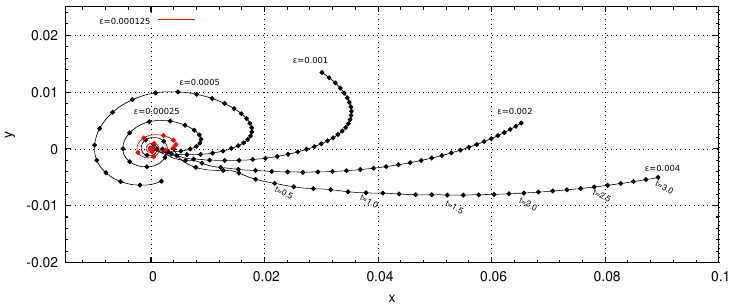}}
  \caption{The computed time-trajectories of the initially-right side spiral, for solutions of Problem~2 
  for the indicated values of $\epsilon$. 
  The spacing between marked points is $\Delta{t}=0.1$.
  For clarity, only points on the $\epsilon=0.004$ curve are labeled with their $t$-values, and 
  the trajectory for $\epsilon=0.000125$ is colored red.  
  Data points at $t=1$ and $t=2$ on each of the six curves correspond
  to contour plots depicted in Fig.~\ref{fig:EpsTo0-t=1} and Fig.~\ref{fig:EpsTo0-t=2} respectively.  
  }
  \label{fig:TimeTraj}
\end{figure}

We now examine the {\it time evolution} of spiral centers in solutions of Problem~2 for fixed selected 
values of $\epsilon$.  
Only one of these centers, the one that initially forms in the right half-plane, 
is tracked and displayed here and throughout.  
The other center is symmetric about the origin to the one that is tracked.
Our time-marching scheme outputs the spatial location of maximum vorticity at every time step, 
but only for the spiral center currently in the right half-plane.    
By post-processing we first detect and correct the actual location of the tracked spiral, 
accounting for its crossings between half-planes.  
The second post-processing step is to reduce the aliasing effect due to sub-cell motion of the spiral center;  
this effect is analogous to the aliasing error (``the jaggies'') encountered in computer graphics when
performing rasterization.  
This goal is accomplished by passing the raw time-dependent location data 
through a piecewise-cubic least squares filter.  
We refer to the resulting smoothed paths as {\it time-trajectories}.

Figure~\ref{fig:TimeTraj} is a plot of spiral center time-trajectories with data points marked at intervals 
of $\Delta{t}=0.1$ for each indicated value of $\epsilon$.  
Data points at $t=1$ and $t=2$ correspond to the solutions whose contour plots are depicted in 
Fig.~\ref{fig:EpsTo0-t=1} and Fig.~\ref{fig:EpsTo0-t=2}. 
Because points on trajectories when $\epsilon=0.00025$ and $\epsilon=0.000125$ lie so close to each other, 
for reasons of graphical clarity the $\epsilon=0.000125$ curve is plotted in red.
The figure indicates that the smaller the value of $\epsilon$, the more closely and rapidly the twin spirals wind around the origin, in accordance with our scaling arguments given in \S\ref{sec:scaling}.  
This observed phenomenon is demonstrated in greater detail in the next paragraph.
The parameter values and grid sizes used in the computations depicted in Fig.~\ref{fig:TimeTraj}
are identical to those used in the earlier figures.

\begin{figure}[h]
  \centerline{\includegraphics[width=0.85\textwidth]{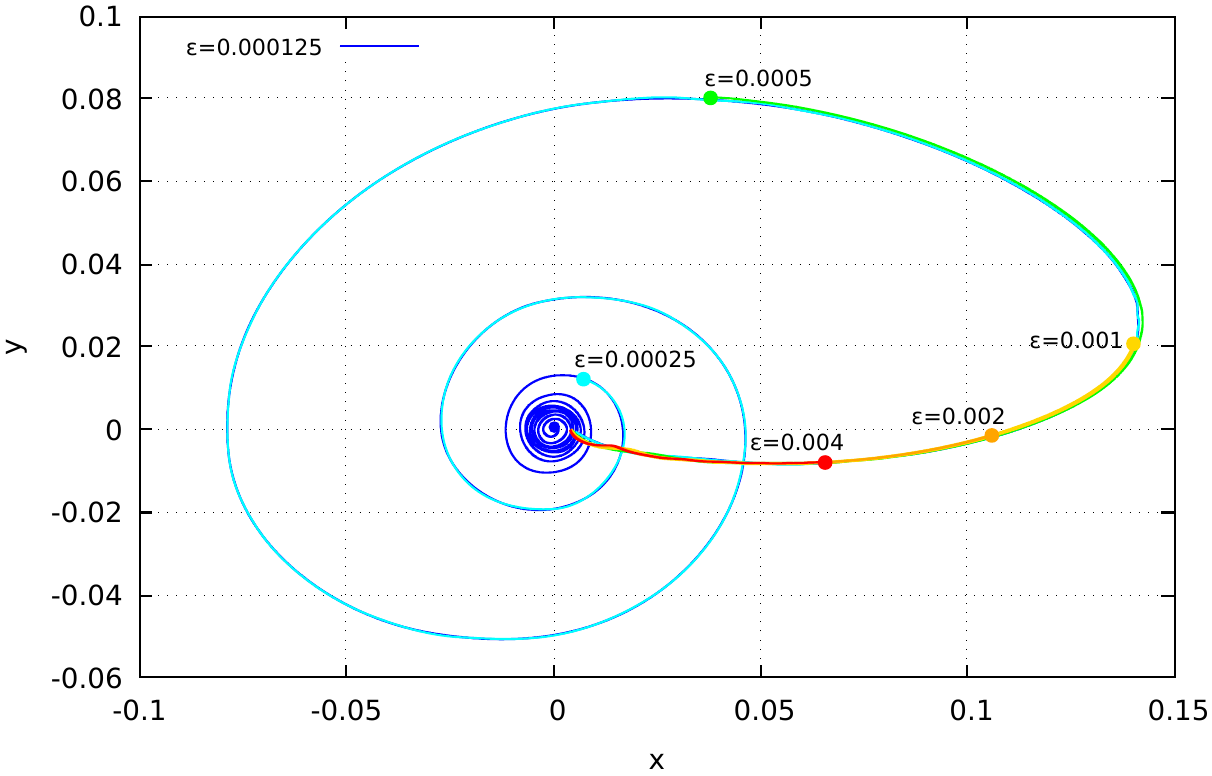}}
  \caption{Validation that our nested-grid calculations accurately obey the scale-invariance
  satisfied by the exact solution.  Given are six superimposed scaled time-trajectories from our 
  computed data displayed in Fig.~\ref{fig:TimeTraj} for different values of $\epsilon$. 
  Here time is restricted to $0 \le t \le 2$.
  The labeled point on each curve represents the scaled spiral position at $t=2.0$.
}
  \label{fig:overlay-mtraj}
\end{figure}

Let $\V{\bar{x}}_\epsilon(t)$ denote the six computed trajectories displayed in 
Fig.~\ref{fig:TimeTraj}, 
\[
   \V{\bar{x}}_\epsilon(t) \ \hbox{with $\epsilon = 0.004, \, 0.002, \ldots, 0.000125$}.
\]
From the scaling identity (\ref{traj-symmetry}) given in \S\ref{sec:scaling}, the trajectories for the 
exact solution to the full free-space problem, (\ref{euler-psiomega}), (\ref{approx-omega0}) 
and (\ref{phi-def2}), obey 
\begin{equation}
  (\epsilon_0/\epsilon)\, \V{x}_{\epsilon} ( t ) = \V{x}_{\epsilon_0}((\epsilon_0/\epsilon)^\alpha \, t).
  \label{eq:traj-symm-mod}
\end{equation}
While equality (\ref{eq:traj-symm-mod}) is true for the exact solution, numerically there is discretization error, 
and we are working on a very large but still finite computational domain.
We take the six curves, $\V{\bar{x}}_\epsilon(t)$, 
and scale each by 
$\epsilon_0/\epsilon$, where $\epsilon_0 = 0.004$ and $\epsilon$ ranges over the other five values.
In Fig.~\ref{fig:overlay-mtraj} we plot 
\[
  (\epsilon_0/\epsilon)\, \V{\bar{x}}_{\epsilon} ( t ) \ \ \hbox{for $0 \le t \le 2$,}
\]
and see they very accurately overlay each other, as expected from (\ref{eq:traj-symm-mod}).
The curves are rendered in the order of smallest to largest $\epsilon$, so that for each successively larger 
$\epsilon$ the corresponding curve is on top.  %
The labeled dots indicate each curve's endpoint.  
Observe that use of our nested but finite computational domains has essentially no adverse effects on the 
expected agreement with the free-space scaling.
The small but noticeable discrepancy between trajectory positions predicted by (\ref{eq:traj-symm-mod}) 
near $(0.004, 0)$ is almost certainly due exclusively to discretization error. 

\begin{figure}[!h]
  \centerline{\includegraphics[width=1.0\textwidth]{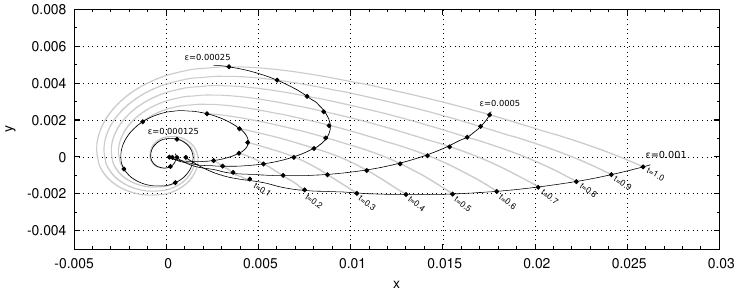}}
  \caption{The definite trend as $\epsilon$ tends to zero of spiral center locations 
  is illustrated here.  
  The gray curves depict spiral locations as functions of decreasing $\epsilon$ 
  at the indicated fixed $t$,
  while the black curves depict locations as functions of $t$ at indicated fixed 
  $\epsilon$, as seen in Fig.~\ref{fig:TimeTraj}.
  The full grid nesting covers a physical domain $[-64,64]\times[-64,64]$, 
  with each 10/11 subgrids individually containing $4096^2$ grid cells.
  }
  \label{fig:EpsTo0}
\end{figure}

We are confident that the scaling identity (\ref{eq:traj-symm-mod}) is
accurately obeyed when applied to our numerical solutions.
It can be exploited to derive a family of curves, which we call 
{\em epsilon-trajectories}, used to further illustrate the solution 
trends observed earlier in Fig.~\ref{fig:TimeTraj}.
Recall we currently have data needed to determine six
time-trajectories, 
$\V{\bar{x}}_{\epsilon_0}(t), \ldots, \V{\bar{x}}_{\epsilon_5}(t)$,
$\epsilon_0 = 0.004,\ldots,\epsilon_5 = 0.000125$.
Each are defined on $0\le t\le 3$.
These are mapped to $\epsilon = \epsilon_0$ via (\ref{eq:traj-symm-mod})
and merged into one trajectory, say $\V{\tilde{x}}_{\epsilon_0}(t)$.
The overlaps depicted in Fig.~\ref{fig:overlay-mtraj} are discarded in
a manner to help minimize the discretization error seen and discussed earlier.
The new merged curve $\V{\tilde{x}}_{\epsilon_0}(t)$ is defined on
$0 \le t \le 3.0 \, (\epsilon_0/\epsilon_5)^\alpha$.
A so-called epsilon-trajectory is denoted by $\V{\bar{y}}_{t_k}(\epsilon)$ 
for a given fixed $t_k$
and displays the locations of spiral centers parameterized with respect to a 
continuously defined variable $\epsilon$, 
$\epsilon_5 \le \epsilon \le \epsilon_0$.
Each is deduced from (\ref{eq:traj-symm-mod}) as follows.
\[
  \V{\bar{y}}_{t_k}(\epsilon) 
   \equiv \V{\bar{x}}_\epsilon(t_k) 
  =
  \left(\epsilon/\epsilon_0\right)
  \V{\tilde{x}}_{\epsilon_0}(  (\epsilon_0/\epsilon)^\alpha t_k  )
  \ \ \hbox{for $\epsilon_5 \le \epsilon \le \epsilon_0$.} 
\]
Regard this as a particular type of interpolation procedure
which simply takes advantage of the given problem's symmetry.
In Fig.~\ref{fig:EpsTo0}, we display four black time-trajectories 
which for visual clarity are restricted to 
$\epsilon = 0.001$, $0.0005$, $0.00025$, $0.000125$ and 
$0 \le t \le 1$. 
Data points are marked at intervals $\Delta t = 0.1$.
The ten epsilon-trajectories, with $t_1=0.1,\ldots,1.0=t_{10}$, 
are colored in light gray.
These ten curves tell us for a given fixed time
how spiral locations change when the continuously defined parameter 
$\epsilon$ is reduced in value.

The interested reader is encouraged to compare and contrast results we give in our
Figs.~\ref{fig:EpsTo0-t=1},~\ref{fig:EpsTo0-t=2},~\ref{fig:TimeTraj} to those
in \cite[Example~3]{bjl2021}.

\subsubsection{The effect of initial data parameter $\alpha$.}
\label{subsubsec:param-values-alpha}

The solutions presented up to this point have been computed using  
baseline values of $\alpha = 0.95$, $\theta_w = \pi/8$.  
Here we consider the effect of the radial exponent $\alpha$  
on the $\epsilon$-family of solutions; 
in \S\ref{subsubsec:param-values-theta} we consider the effect of the wedge angle $\theta_w$.

\begin{figure}[!h]
  \centerline{
  \vbox{\hbox{\includegraphics[width=0.725\textwidth]{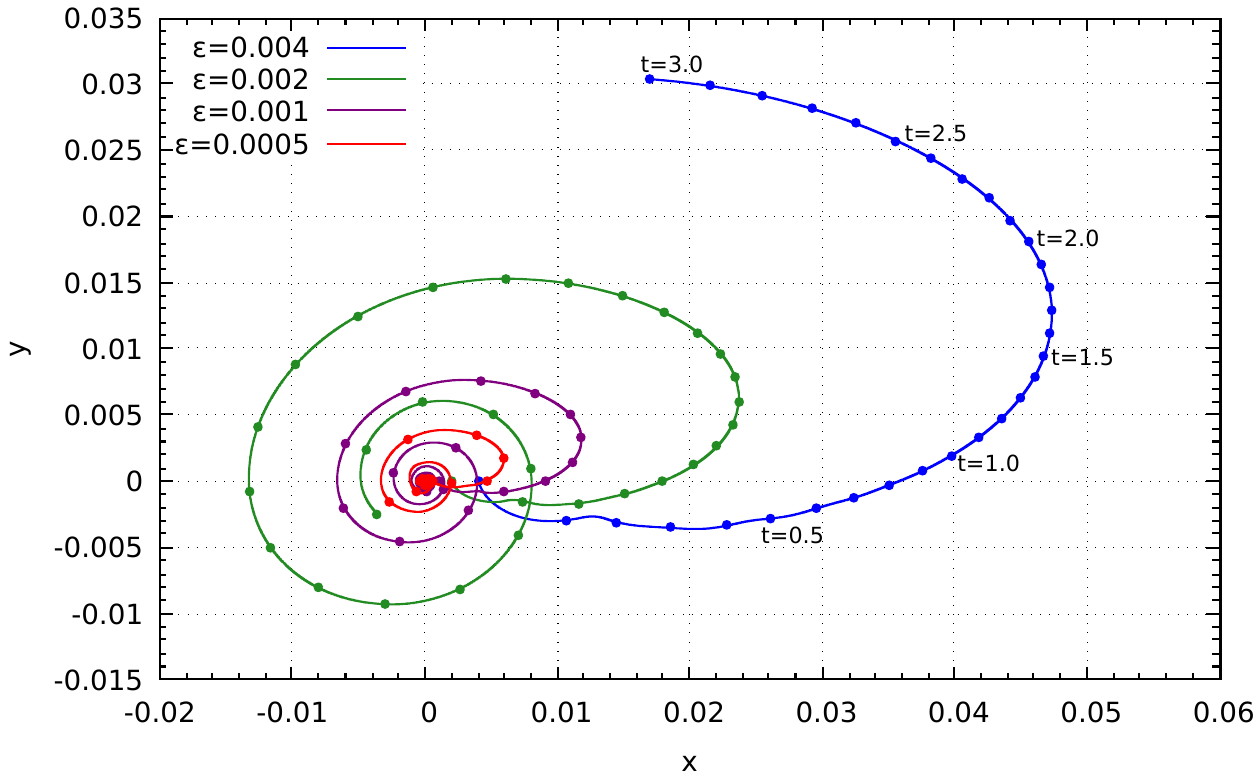}}
            \hbox{(a) $\alpha = 1.05$, trajectories for $0 \le t \le 3.$}}
            }  \vspace{3pt}
  \centerline{
  \vbox{\hbox{\includegraphics[width=0.725\textwidth]{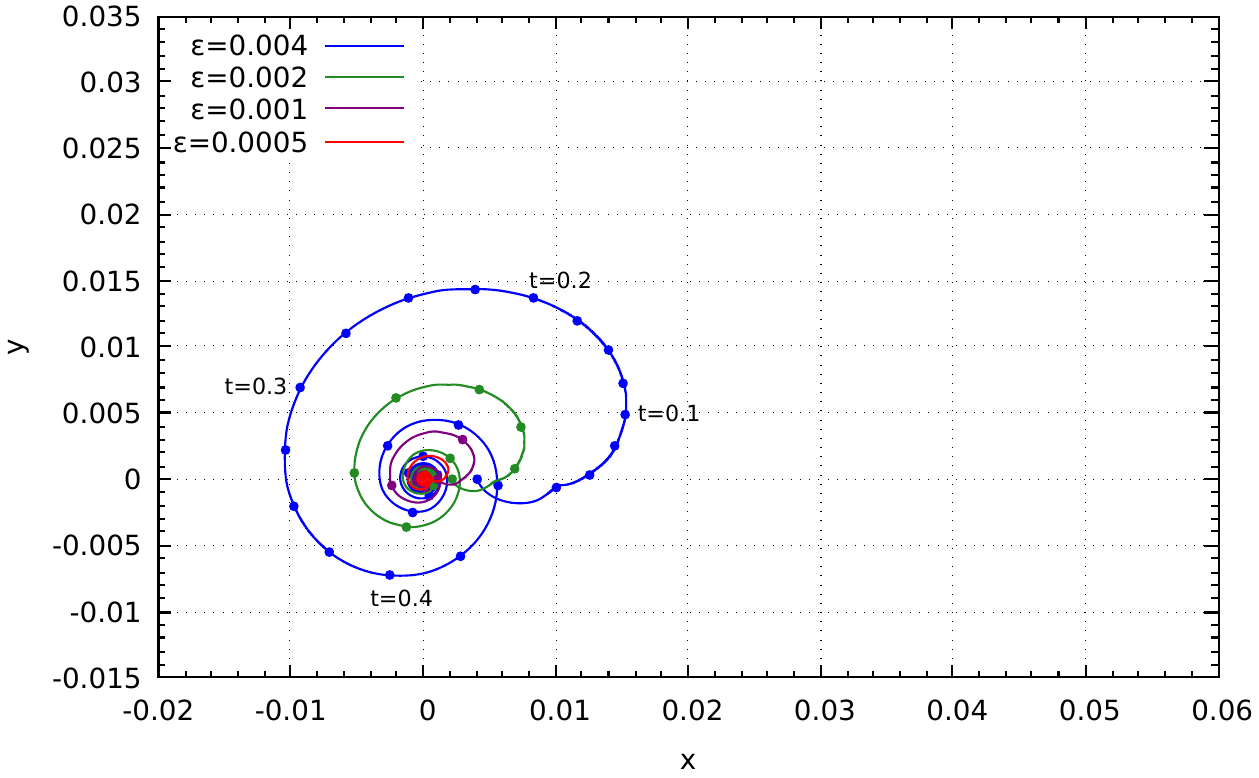}}
            \hbox{(b) $\alpha = 1.4$, trajectories for $0 \le t \le 0.8.$}}
            }   
  \caption{Trajectories of the right-side spiral showing the effect of $\alpha$. 
  The spiral motion takes place in a small portion of the full final grid which 
  has extent $[-0.125,0.125]\times[-0.125,0.125]$.  
  In (a) data points are plotted at intervals $\Delta{t}=0.1$.  
  In (b) data points are plotted at intervals $\Delta{t}=0.025$.
  Solutions were computed using grids containing $4096^2$ grid cells,
  and $\theta_w = \pi/8$.  
  Compare the time-trajectories here at identical $\epsilon$-values to 
  those in Fig.~\ref{fig:TimeTraj} where $\alpha= 0.95$.
   }
   \label{fig:TimeTraj-1.05+1.4}
\end{figure}

Figure~\ref{fig:TimeTraj-1.05+1.4} depicts time-trajectories for solutions of Problem~2 
with larger values of $\alpha=1.05$~(a) and $\alpha=1.4$~(b) 
which are both displayed on identically sized regions.  
Compare the time-trajectories in this figure to those in Fig.~\ref{fig:TimeTraj} 
for identical values of $\epsilon$, for example $\epsilon_0 = 0.004$.
Observe that the larger the value of $\alpha$, 
the more closely and rapidly the twin spirals orbit around the origin. 
In our calculations we always halve consecutive values of $\epsilon$,
so that the ratio $\epsilon_0 / \epsilon_{k} = 2^k$.
From (\ref{eq:traj-symm-mod}), for $\epsilon_k = 2^{-k} \epsilon_0$,
relative to $\epsilon_0$, the problem is expected to scale 
in space by $2^{-k}$ and in time by $(2^k)^\alpha$.  
Therefore, with decreasing $\epsilon$, for larger values of $\alpha$
this time scaling exacerbates the increase in spiral rotation rate and
collapse to the origin, as just observed for $\epsilon = 0.004$.
Moreover, for $\alpha = 1.05$ and $\alpha = 1.4$ 
this rotation and collapse is so rapid that the 
$\epsilon=0.00025,\ 0.000125$ trajectories are not given, 
and in (b) we only plot the spiral's motion up to time $t=0.8$.

\begin{figure}[h]
  \centerline{
    \vbox{\hbox{\includegraphics[width=0.45\textwidth]{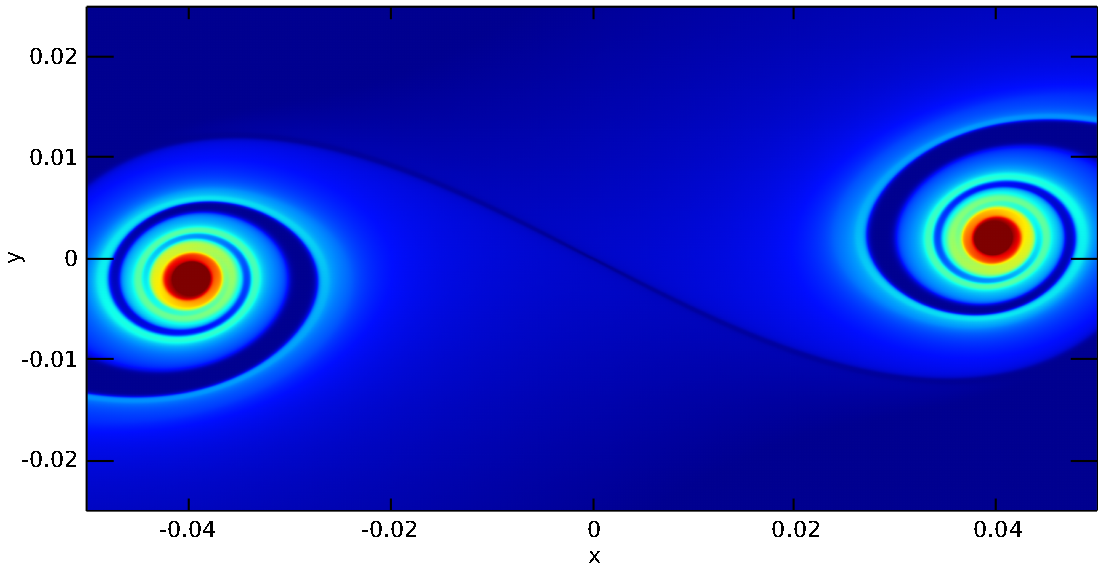}}
    \hbox{\quad (a) {$\epsilon=0.004, \, t=1.0$ }}}
    \hskip0ex
    \vbox{\hbox{\includegraphics[width=0.45\textwidth]{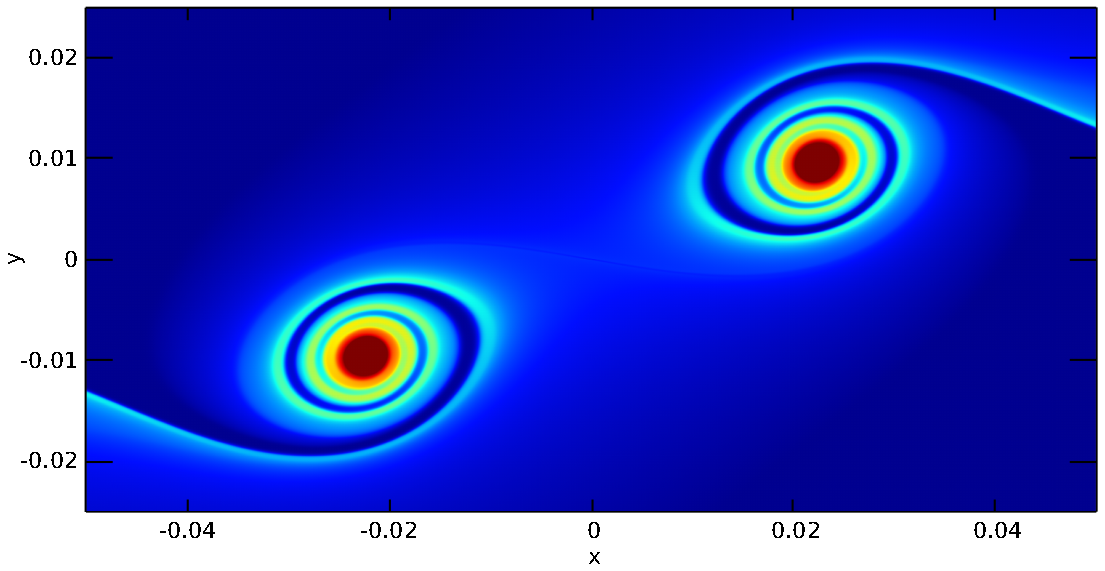}}
    \hbox{\quad (b) {$\epsilon=0.002, \, t=1.0$}}}}
  \vskip0.25cm
  \centerline{
    \vbox{\hbox{\includegraphics[width=0.45\textwidth]{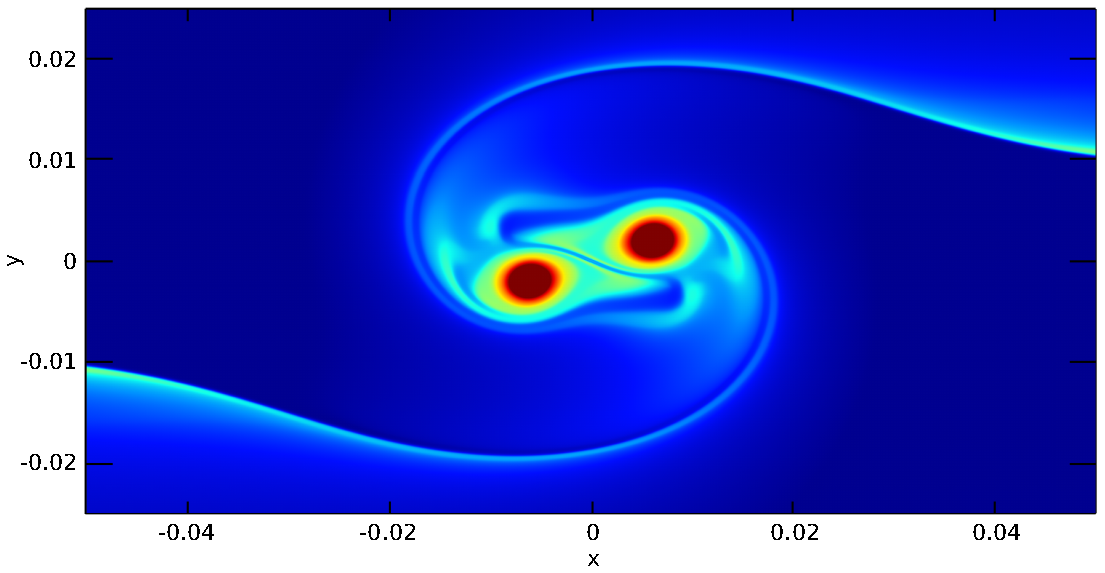}} 
    \hbox{\quad (c) {$\epsilon=0.001, \, t=1.0$}}}
    \hskip0ex
    \vbox{\hbox{\includegraphics[width=0.45\textwidth]{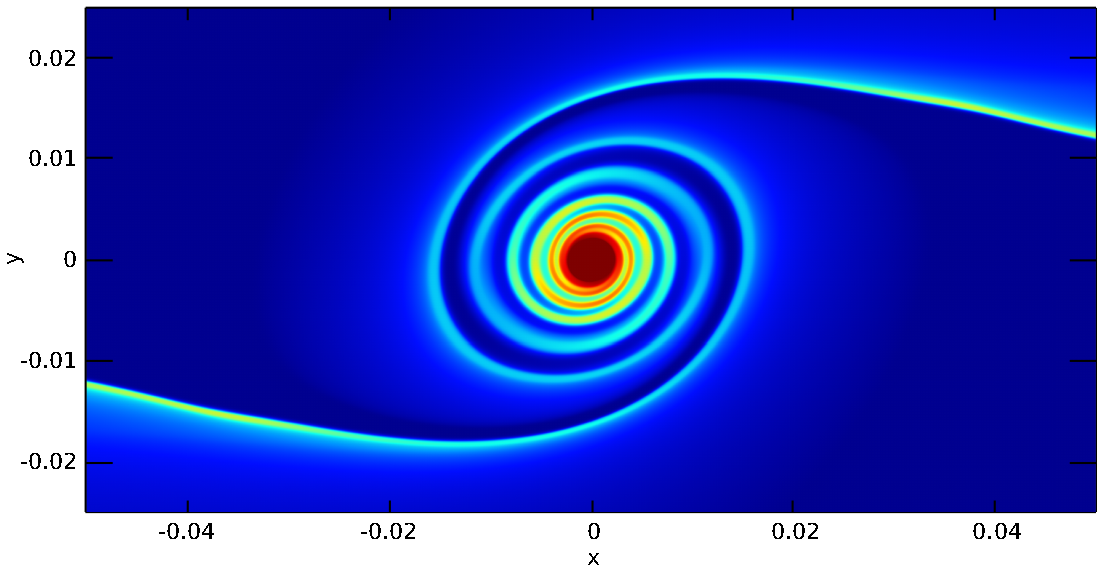}}
    \hbox{\quad (d) {$\epsilon=0.0005, \, t=1.0$}}}
  }
  \caption{Vorticity contour plots for a decreasing sequence of 
  $\epsilon$ with $\alpha=1.05$; 
  compare to the plots in Fig.~\ref{fig:EpsTo0-t=1} which use $\alpha=0.95$.
  The solutions depicted here correspond to data points in 
  Fig.~\ref{fig:TimeTraj-1.05+1.4} (a) at $t=1$.  
  For (a)--(d), $\theta_w = \pi/8$.  
  }
  \label{fig:Cont:Alpha1.05}
\end{figure}

In light of our numerical results presented in Fig.~\ref{fig:TimeTraj-1.05+1.4}
and the scaling relation (\ref{eq:traj-symm-mod}) visualized in Fig.~\ref{fig:EpsTo0},
it is evident that for any fixed $t > 0$ the twin spirals will collapse to
a single spiral as $\epsilon$ tends to zero.
To verify this fact for the case $\alpha=1.05$, in Fig.~\ref{fig:Cont:Alpha1.05} we show contour plots 
at $t=1$ of the four solutions whose time-trajectories are depicted in Fig.~\ref{fig:TimeTraj-1.05+1.4}~(a). 
As shown, the twin spirals have essentially merged even at the relatively large value $\epsilon = 0.0005$.  
Figure~\ref{fig:TimeTraj-1.05+1.4}~(b) indicates that for $\alpha = 1.4$ 
a similar set of contour plots would show the same behavior, and at much earlier times than $t=1$. 

\subsubsection{The effect of initial data parameter $\theta_w$.}
\label{subsubsec:param-values-theta}

\begin{figure}[!b]
  \centerline{
  \vbox{\hbox{\includegraphics[width=0.72\textwidth]{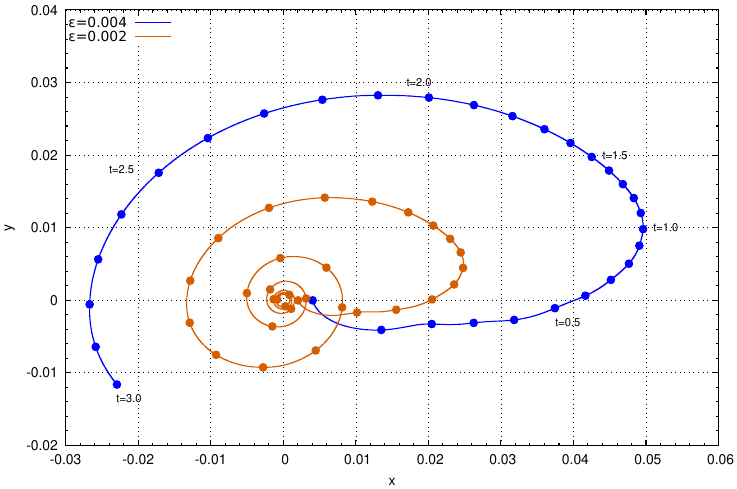}}
            \hbox{(a) $\theta_w = 3 \pi/16,$ trajectories for $0 \le t \le 3$.}}
            }  \vskip1em
  \centerline{
  \vbox{\hbox{\includegraphics[width=0.72\textwidth]{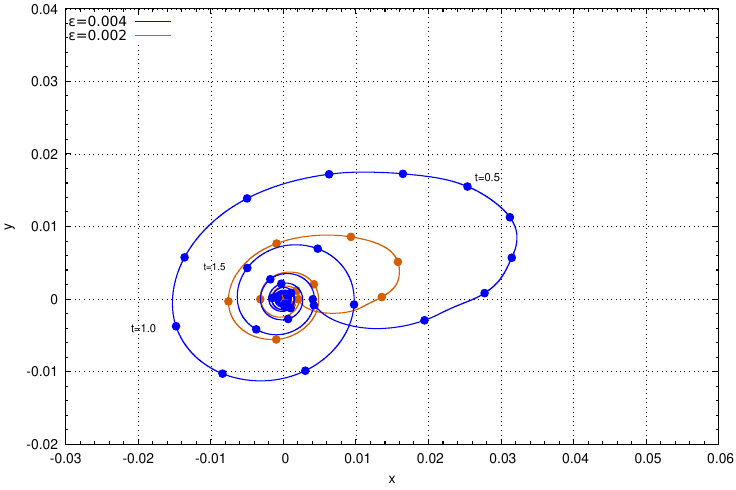}}
            \hbox{(b) $\theta_w = \pi/4,$ trajectories for $0 \le t \le 3$.}}
            }   
  \caption{Trajectories of the right-side spiral showing the effect of $\theta_w$. 
  The spiral motion takes place in a small portion of the full final grid 
  which has extent $[-0.125,0.125]\times[-0.125,0.125]$.  
  Data points are plotted at intervals $\Delta{t}=0.1$.  
  The solutions shown were computed with grids 
  of size $4096^2$ and $\alpha=0.95$. Compare the plots in (a) and (b) to 
  Fig.~\ref{fig:TimeTraj} which has 
  $\theta_w= \pi/8$.
  }
   \label{fig:TimeTraj-3pio16+pio4}
\end{figure}

Here we consider the effect of changing the wedge angle $\theta_w$ on the $\epsilon$-family of solutions.
We will use the results presented in Fig.~\ref{fig:TimeTraj}, 
where $\theta_w=\pi/8$, as a basis of comparison.  
Figure~\ref{fig:TimeTraj-3pio16+pio4} shows spiral time-trajectories 
in solutions of Problem~2 with $\alpha = 0.95$ and $\theta_w = 3 \pi/16$ in (a) and 
$\theta_w = \pi/4$ in (b).  
For visual clarity, only trajectories for solutions with $\epsilon=0.004$ and 
$\epsilon=0.002$ are depicted in these plots.  
Comparing time-trajectories here to those in Fig.~\ref{fig:TimeTraj} 
at the same $\epsilon$ values shows that the larger the value of $\theta_w$, 
the more rapidly and closely the twin spirals orbit the origin.  
In fact, the blue curve in Fig.~\ref{fig:TimeTraj-3pio16+pio4}~(b) shows that when 
$\theta_w = \pi/4$ and $\epsilon=0.004$, the tracked spiral circuits around the origin approximately 8 times in three units of time.  
Moreover, its distance from the origin is eventually decreasing.

\begin{figure}[!b]
  \centerline{
    \vbox{\hbox{\includegraphics[width=0.45\textwidth]{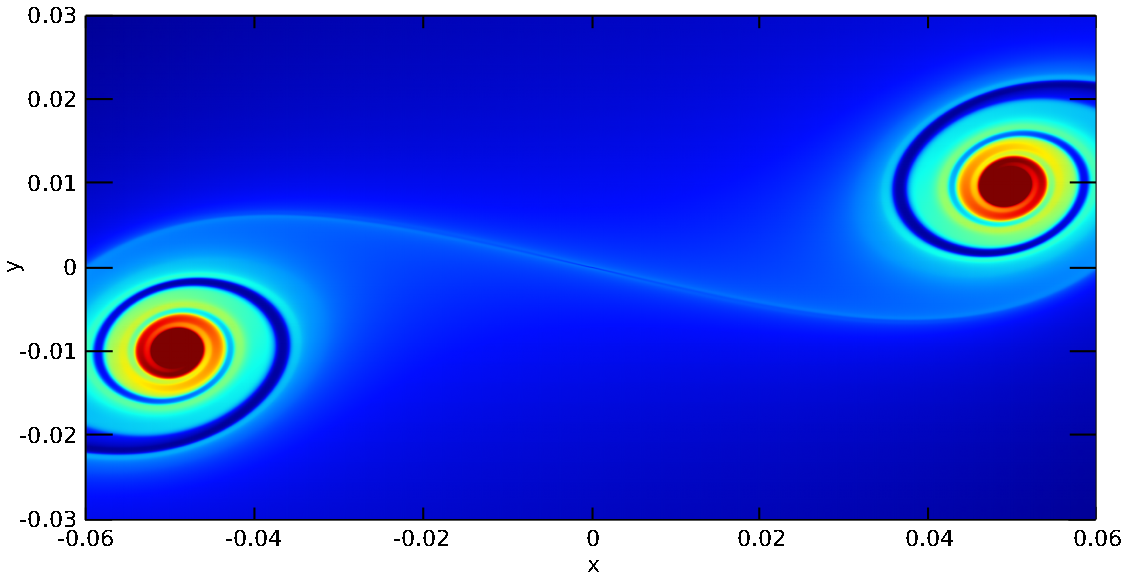}}
              \hbox{\hskip \dimexpr 0.225\textwidth - 2em \relax $\epsilon=0.004$}}
    \vbox{\hbox{\includegraphics[width=0.45\textwidth]{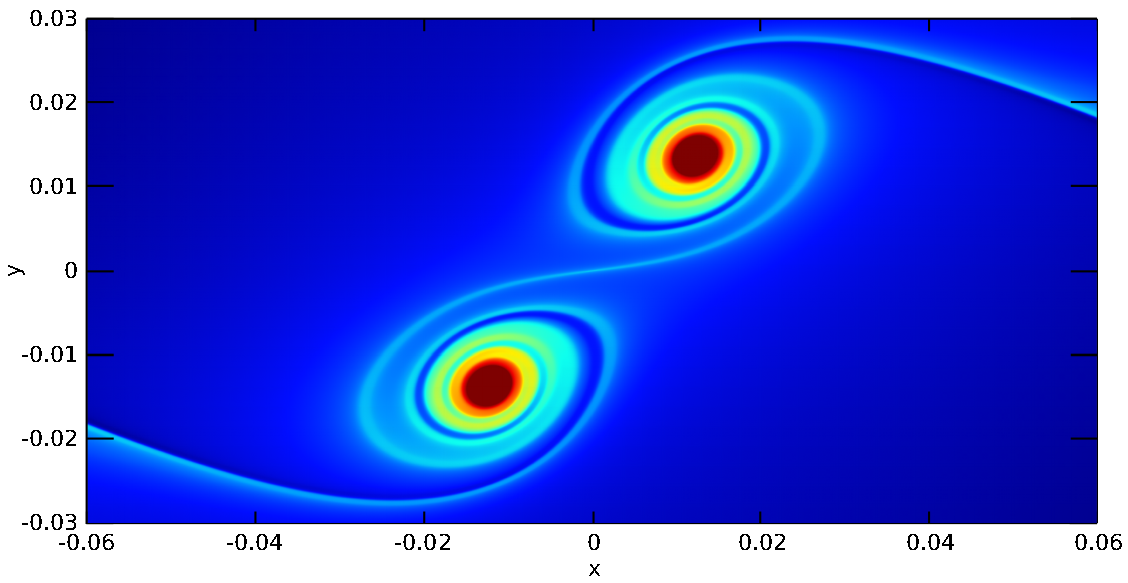}}
              \hbox{\hskip \dimexpr 0.225\textwidth - 2em \relax $\epsilon=0.002$}}
  } 
  \smallskip
  \centerline{  
    \vbox{\hbox{\includegraphics[width=0.45\textwidth]{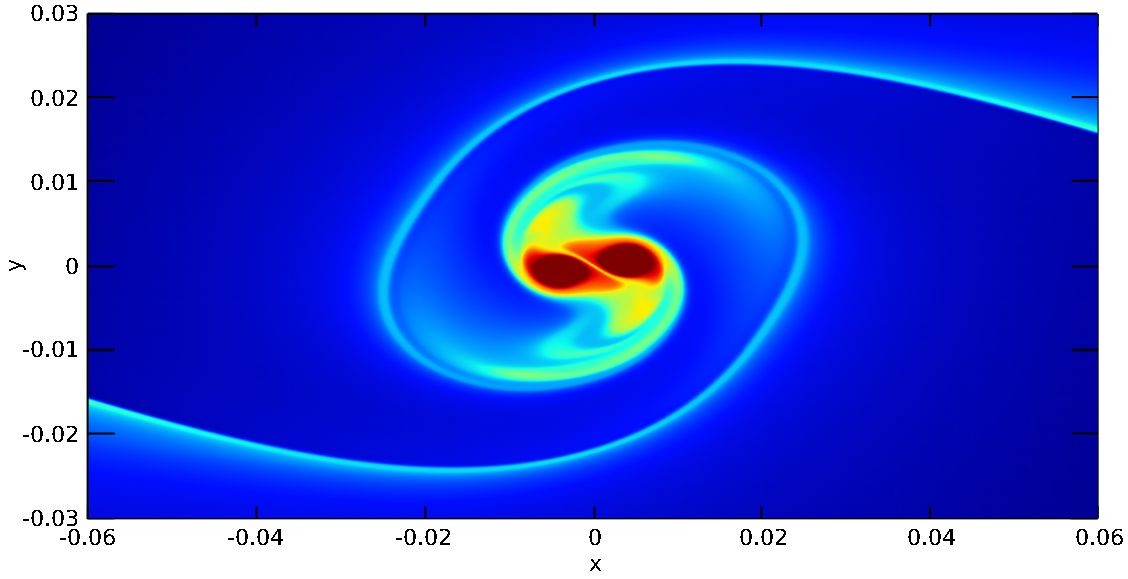}} 
              \hbox{\hskip \dimexpr 0.225\textwidth - 2em \relax $\epsilon=0.001$}}
    \vbox{\hbox{\includegraphics[width=0.45\textwidth]{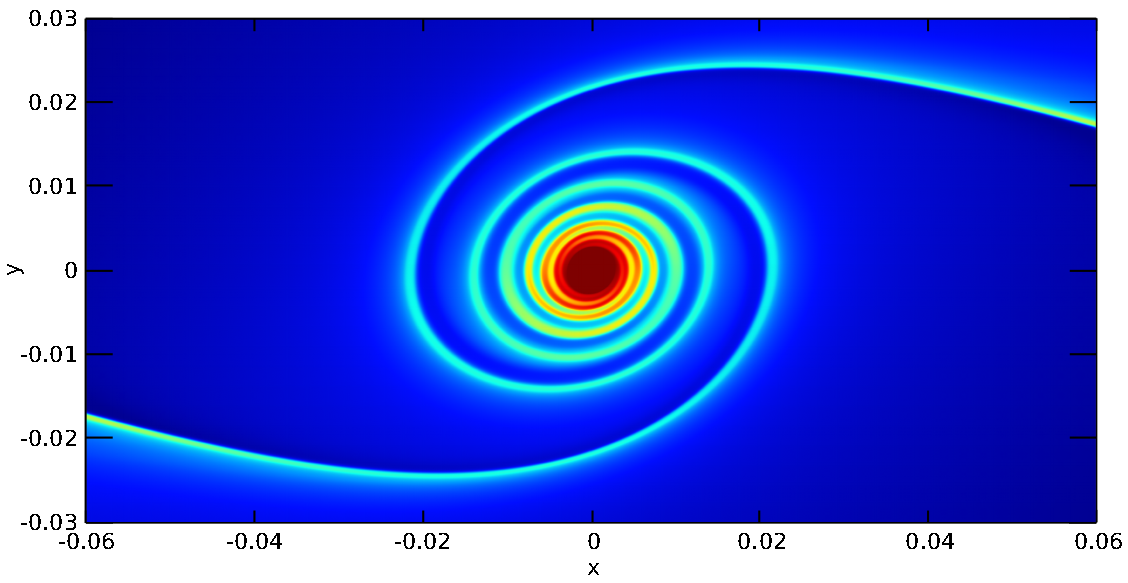}}
              \hbox{\hskip \dimexpr 0.225\textwidth - 2em \relax $\epsilon=0.0005$}}
  }
  \hbox{(a)  $\theta_w = 3 \pi/16$, $t=1$. }
  \vspace{2ex}
  \centerline{  
    \vbox{\hbox{\includegraphics[width=0.45\textwidth]{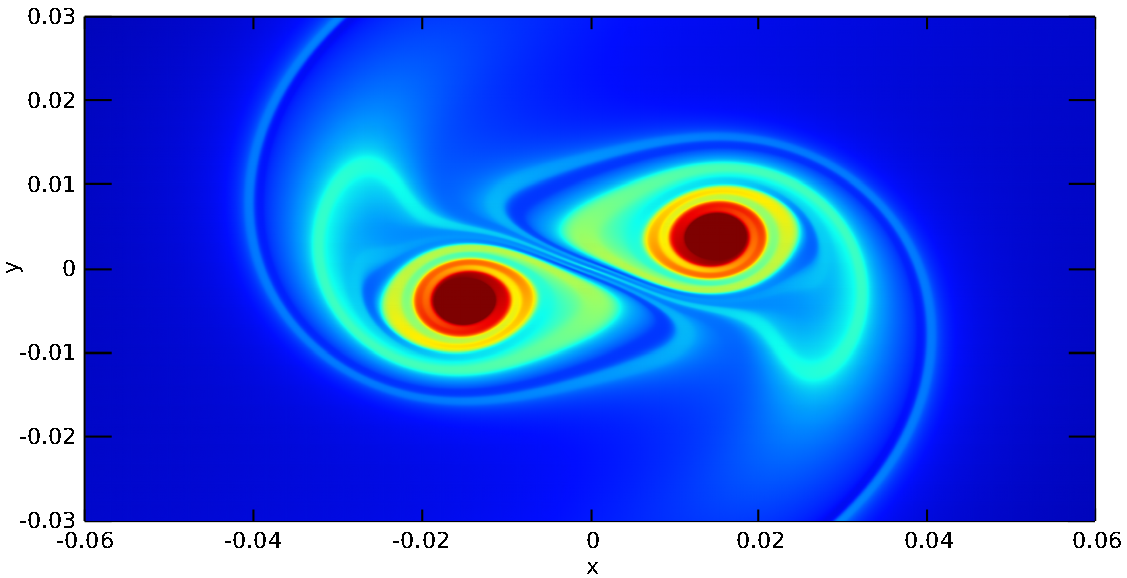}}
              \hbox{\hskip \dimexpr 0.225\textwidth - 2em \relax $\epsilon=0.004$}}
    \hskip0ex
    \vbox{\hbox{\includegraphics[width=0.45\textwidth]{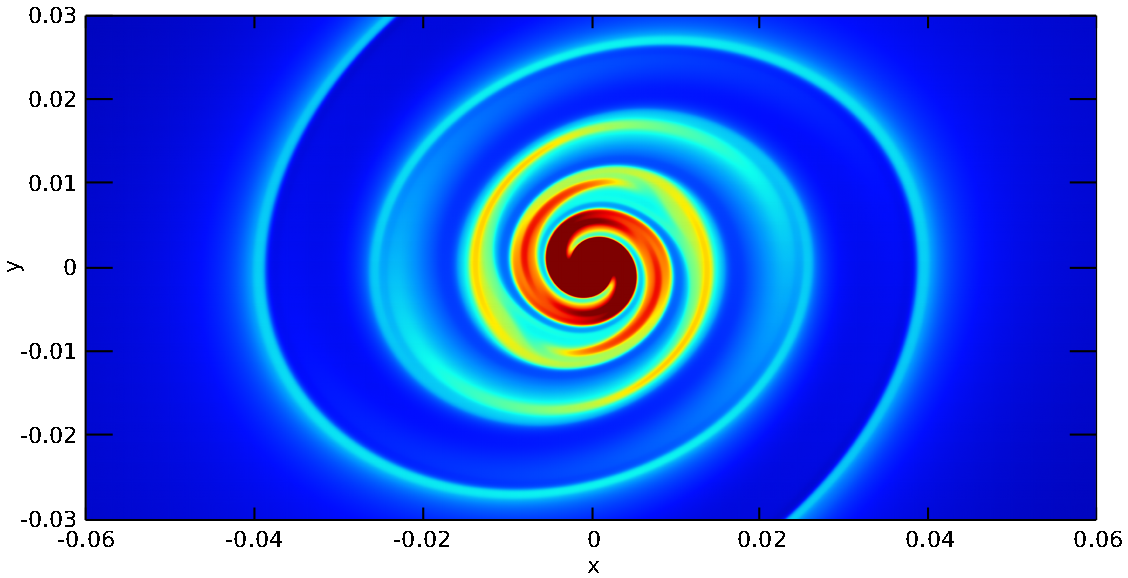}}
              \hbox{\hskip \dimexpr 0.225\textwidth - 2em \relax $\epsilon=0.002$}}
  }
  \hbox{(b)  $\theta_w = \pi/4$, $t=1$. }
  \caption{Vorticity contour plots for decreasing sequences of $\epsilon$, for two values of $\theta_w$.
  In (a) $\theta_w = 3 \pi/16$; contour plots displayed in the top row correspond to data points in 
  Fig.~\ref{fig:TimeTraj-3pio16+pio4}~(a) at $t=1$.
  In (b) $\theta_w = \pi/4$; contour plots depicted correspond to data points in 
  Fig.~\ref{fig:TimeTraj-3pio16+pio4}~(b) at $t=1$.  Compare the sequences of contours plots in 
  (a) and (b) to the sequence of plots displayed in Fig.~\ref{fig:EpsTo0-t=1} for which $\theta_w = \pi/8$.  
  All solutions depicted here use $\alpha=0.95$.  
  }
  \label{fig:Cont:3pio16+pio4}
\end{figure}

In light of our numerical results presented in Fig.~\ref{fig:TimeTraj-3pio16+pio4}
and the scaling relation (\ref{eq:traj-symm-mod}) visualized in Fig.~\ref{fig:EpsTo0},
it is once again evident that for any fixed $t > 0$ the twin spirals will collapse to
a single spiral as $\epsilon$ tends to zero.
To verify this fact for the case $\theta_w = 3\pi/16$, in Fig.~\ref{fig:Cont:3pio16+pio4}~(a) 
we show contour plots 
at $t=1$ of the two solutions whose time-trajectories are depicted in Fig.~\ref{fig:TimeTraj-3pio16+pio4}~(a)
and add, in addition to these, contour plots for two smaller values, $\epsilon=0.001$, $\epsilon=0.0005$. 
As shown, the twin spirals have for practical purposes merged at the relatively large 
value $\epsilon = 0.0005$.  
We see very similar phenomena with the slightly larger value of $\theta_w = \pi/4$.  
In Fig.~\ref{fig:Cont:3pio16+pio4}~(b) we display vorticity contour plots of the solutions 
whose trajectories were depicted in Fig.~\ref{fig:TimeTraj-3pio16+pio4}~(b).
As shown, the transition from two-spiral to one-spiral solution occurs at a value of $\epsilon$ between $0.004$ and $0.002$.

\begin{figure}[!h]
  \centerline{\includegraphics[width=0.9\textwidth]{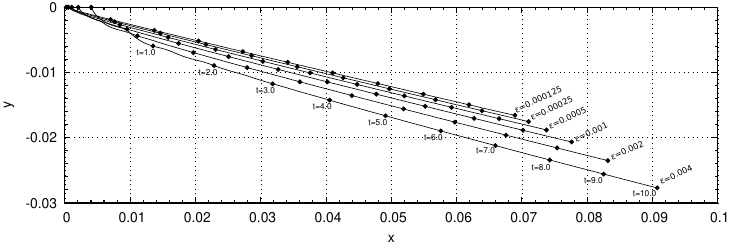}}
  \caption{Trajectories of the right-side spiral when $\theta_w = \pi/16$. 
  Times are marked at intervals of $\Delta{t}=1.0$.  All solutions were 
  computed with grids of size $4096^2$, and $\alpha = 0.95$.  
  Compare to Fig.~\ref{fig:TimeTraj} which has $\theta_w= \pi/8$.
  }
  \label{fig:TimeTraj_piO16}
\end{figure}

In all of the results presented to this point, our numerical calculations indicate that solutions of 
Problem~1 and Problem~2 both converge to a single spiral.  
Now we decrease $\theta_w$ from its baseline value of $\pi/8$ to $\pi/16$
and seek to determine whether we see this phenomenon again.  
First, we take $\alpha=0.95$.
Again Problem~1 resolves into a single spiral but is not displayed 
here for brevity. 
Recall that if the limit solution to Problem~2 converges to a twin spiral self-similar solution, 
(\ref{what-we-expect}) tells us that time-trajectories should
approximately behave like $\V{x}_\epsilon(t) \sim t^{1/\alpha} \V{y}_*$, for some $\V{y}_* \not= \V{0}$, i.e.\ a straight line.   
Figure~\ref{fig:TimeTraj_piO16} is a plot of spiral time-trajectories, 
for $0 \le t \le 10$ at the indicated values of $\epsilon$.
Note that here, with smaller wedge angle, 
we have integrated much further in time than was done earlier.
Consequently, the scale-invariance expressed in (\ref{eq:traj-symm-mod}) 
was again checked, with equally good agreement as depicted earlier in Fig.~\ref{fig:overlay-mtraj}.

The time-trajectories depicted in Fig.~\ref{fig:TimeTraj_piO16} 
do not obviously rule out a twin spiral $\epsilon$-limit solution to Problem~2.  
In order to check further, we double the final time to $t=20$ and 
use our smallest $\epsilon = 0.000125$ only; recall from (\ref{traj-symmetry}) in 
\S\ref{sec:scaling} that the smaller the value of $\epsilon$ the more rapidly a trajectory evolves
in time.  
The time-trajectory is depicted in Fig.~\ref{fig:TimeTraj_piO16_t=20}(a).  
Figure~\ref{fig:TimeTraj_piO16_t=20}(b) depicts the computed values of 
$\V{x}_\epsilon(t) / t^{1/\alpha}$ for $0 \le t \le 20$.  
Given that what may appear here to be horizontal asymptotes are indeed such, 
this provides the first indication in this work of a possible attracting self-similar solution 
with spiral center $\V{y}_* \not= \V{0}$; see ({\ref{what-we-expect}).

\begin{figure}[t]
   \centerline{ \raisebox{13.0ex}{(a)} \quad \hbox{\includegraphics[width=0.8\textwidth]{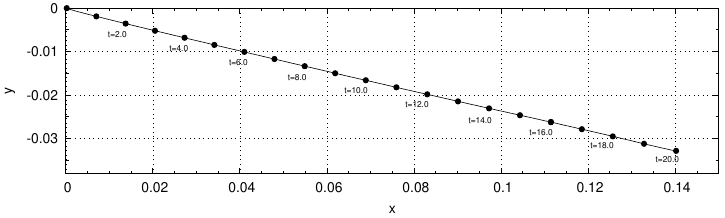}}}  
   \centerline{ \raisebox{16.5ex}{(b)} \quad \hbox{\includegraphics[width=0.8\textwidth]{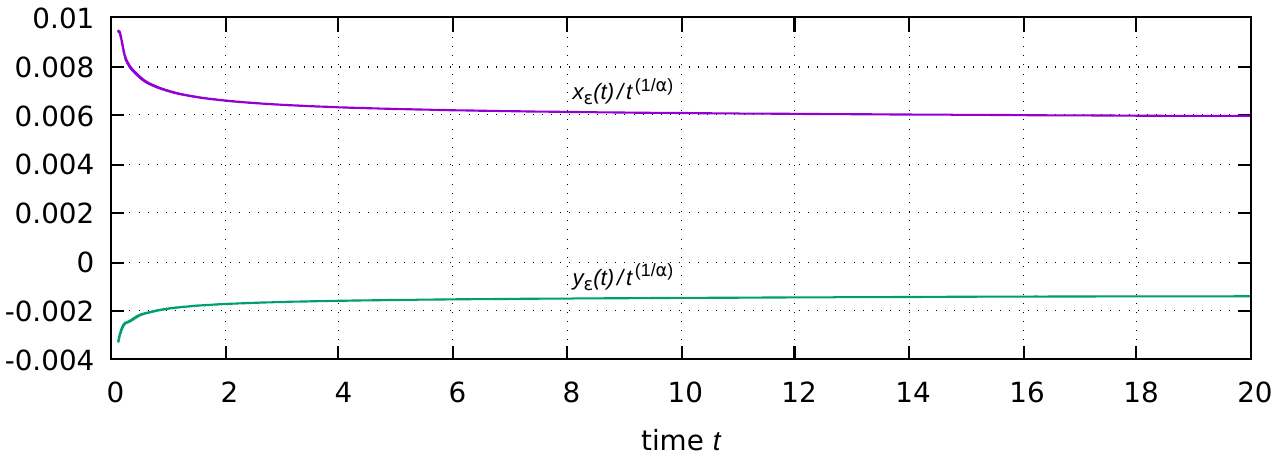}}}   
   \caption{In (a), the trajectory of the right-side spiral with $\epsilon=0.000125$, 
   displayed in Fig~\ref{fig:TimeTraj_piO16} for $0 \le t \le 10$, but here integrated further 
   in time to $t=20$.  In (b), a depiction of $\V{x}_\epsilon(t)/t^{1/\alpha}$ with 
   $\epsilon = 0.000125, \, \alpha=0.95$.
   }
   \label{fig:TimeTraj_piO16_t=20}
\end{figure}

Finally, we keep $\theta_w$ at the value $\pi/16$ and increase $\alpha$ from its baseline
value of $0.95$ to $1.4$.  Figure~\ref{fig:TimeTraj_piO16_1.4} shows spiral time-trajectories
for the $\epsilon$-family of solutions with 
$\epsilon = 0.004,\,0.002,\,0.001,\,0.0005$.  For $\epsilon = 0.00025, \, 0.000125$ the 
rate of rotation and collapse to the origin of spiral centers is so rapid 
that these trajectories are not displayed for reasons of graphical clarity.   
It is once again clear, from our numerical results presented in Fig.~\ref{fig:TimeTraj_piO16_1.4}
and the scaling relation (\ref{eq:traj-symm-mod}) visualized in Fig.~\ref{fig:EpsTo0}, that 
the twin spirals initially present in the Problem~2 data collapse to 
a single spiral as $\epsilon$ tends to zero for any fixed $t > 0$.

\begin{figure}[H]
  \centerline{\includegraphics[width=0.75\textwidth]{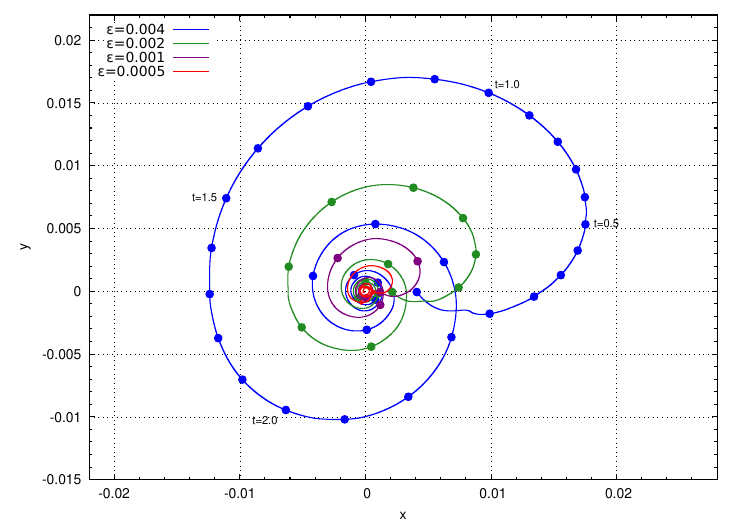}}
  \caption{Time-trajectories of the right-side spiral when $\theta_w = \pi/16$, $\alpha=1.4$,
  for the indicated values of $\epsilon$.  Trajectories with $\epsilon = 0.00025, 0.000125$
  are not displayed for graphical clarity.
  Times are marked at intervals of $\Delta{t}=0.1$.  Solutions computed on all nested grids are 
  of size $4096^2$. Note that here the region in which spiral motion occurs is tiny, only slightly
  larger than that in Fig.~\ref{fig:TimeTraj-1.05+1.4} (b).}
  \label{fig:TimeTraj_piO16_1.4}
\end{figure}


\section{Discussion}
\label{sec:discussion}

In all of our calculations, except one, we have presented clear 
numerical evidence that indicates as $\epsilon$ tends to zero, 
solutions of Problem~2 converge to a solution that contains only one spiral;
see Figs.~\ref{fig:EpsTo0-t=1}--\ref{fig:Cont:3pio16+pio4}, \ref{fig:TimeTraj_piO16_1.4}.
This contradicts statements made in works \cite{ab2021,ab2020,bjl2021},  
discussed in the next paragraph.  The exceptional case is $\alpha=0.95$, 
$\theta_w=\pi/16$; see Figs.~\ref{fig:TimeTraj_piO16},~\ref{fig:TimeTraj_piO16_t=20}.
Our numerical solution for this case does not obviously refute
the possibility of a two-spiral $\epsilon$-limit solution.  The trajectories are close to 
straight lines and move with approximately the right speed 
as would be expected here.  Nevertheless, this does 
not confirm with certainty that the limit solution, in this case, 
is in fact a twin-spiral self-similar solution.

Bressan and Shen in \cite{ab2021} and Bressan and Murray in \cite{ab2020} state  
that their numerical simulations indicate that, as $\epsilon \to 0$, 
two distinct limit solutions are obtained: solutions of Problem~1 contain
a single spiral, while solutions of Problem~2 contain two spirals.  
However, the numerical results presented in each of these works consist solely of two
vorticity contour plots, one depicting a one-spiral and the other a two-spiral solution
(similar to our plots in Fig.~\ref{fig:Our_SingleGrid} (a) and (b)).
The plots in \cite{ab2020} and \cite{ab2021} appear to be identical, 
although parameter values used are not specified.  
Evidence that the numerical solutions depicted in these plots
are in fact limit solutions, or self-similar, is not given in either paper.
However, both works reference the Matlab code \cite{ws2020}
they used to compute their solutions.
Accordingly, in order to compare our numerical results to theirs it was necessary to run their code, 
which we now describe.

\begin{figure}[!b]
  \centerline{\vbox{\hbox{\includegraphics[width=0.85\textwidth]{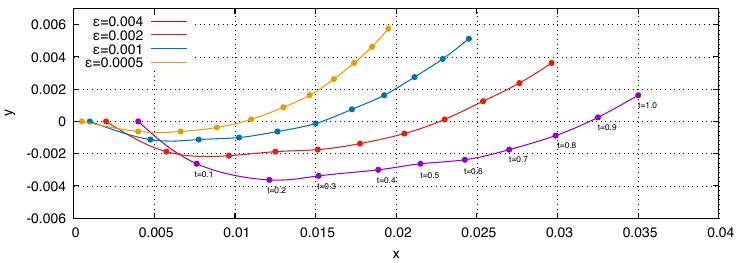}}
                             \vspace{-3ex} \hbox{(a) }}}
  \medskip
  \centerline{
    \vbox{\hbox{\includegraphics[width=0.5\textwidth]{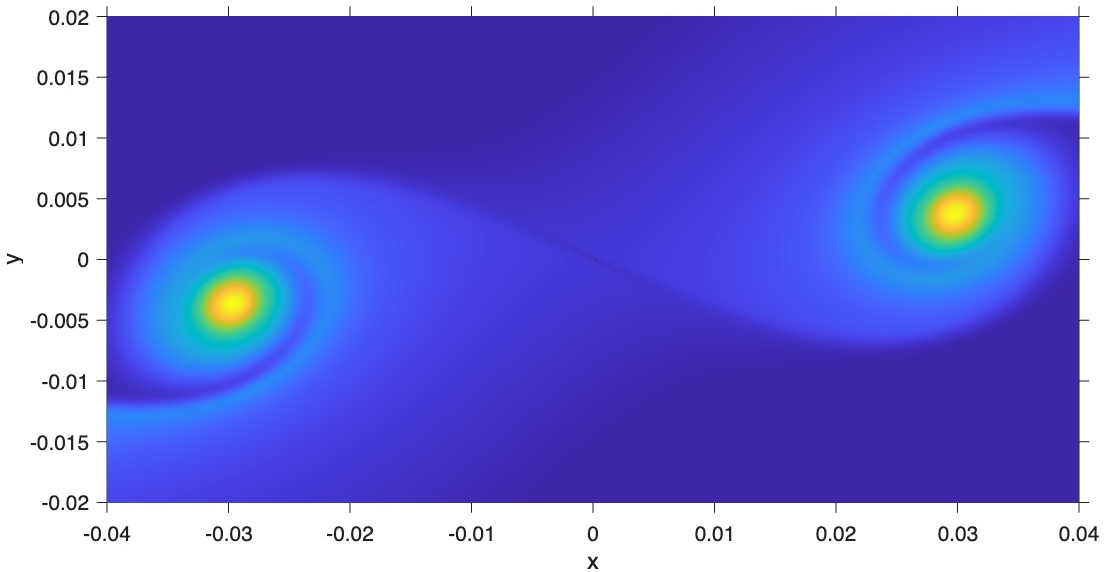}}
              \hbox{(b) $\epsilon=0.002$}}
    \vbox{\hbox{\includegraphics[width=0.5\textwidth]{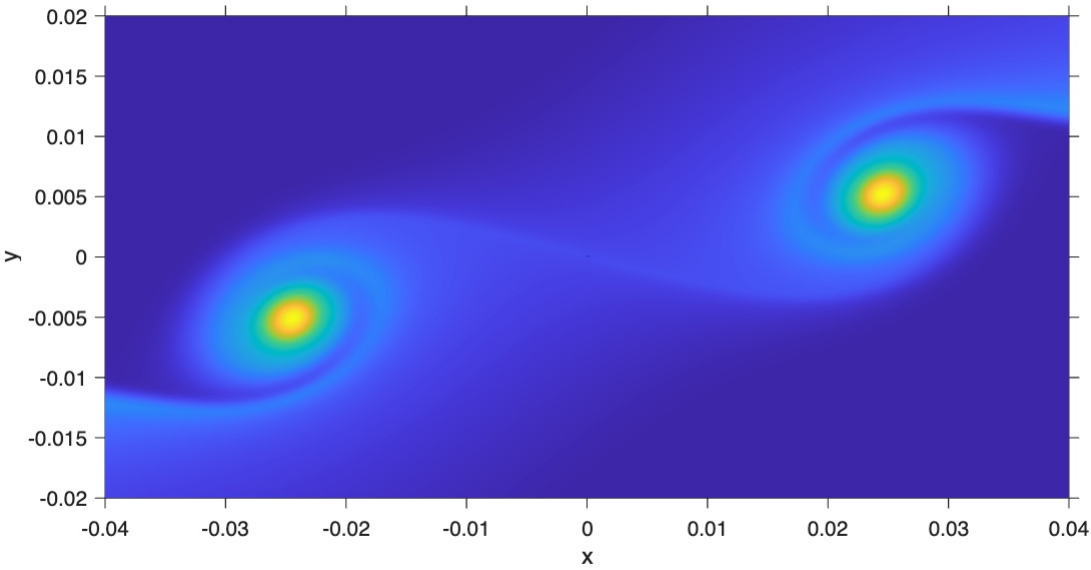}}
              \hbox{(c) $\epsilon=0.001$}}
  }
  \caption{Results using the Matlab code of references \cite{ab2021,ab2020} with 
  $\alpha=0.95$ and $\theta_w = \pi/8$.  Their code uses a single uniform grid,
  here with extent $[-0.2,0.2]\times[-0.2,0.2]$.  
  In (a), time-trajectories for $0 \le t \le 1$ of the initially right-side spiral 
  at the indicated values of $\epsilon$; data points are marked at intervals $\Delta{t}=0.1$.  
  In (b-c), vorticity contour plots for two consecutive values of $\epsilon$ at $t=1.0$.
  }
  \label{fig:Matlab}
\end{figure}

The referenced Matlab code uses a first-order finite-difference scheme on a single 
uniform grid of extent, as of this writing, $[-0.2,0.2]\times[-0.06,0.06]$,
and applies hard-wall boundary conditions $\psi = 0$ on the domain boundary.  
We increased the extent to $[-0.2,0.2]\times[-0.2,0.2]$  
to somewhat mitigate the
effect of the numerical boundaries, 
and set code parameters to solve our Problem 2 for the 
case $\alpha=0.95$, $\theta_w = \pi/8$, taking 
$\epsilon = 0.004,\,0.002,\,0.001$, and $0.0005$.  
To maintain an acceptable ratio $\epsilon/\Delta{x}$ for all of these   
values of $\epsilon$, we set $\Delta{x} = \Delta{y} = 0.000125$ (i.e.\ $3200^2$ grid points). 
All solutions depicted here were integrated in time to $t=1.0$, and   
Fig.~\ref{fig:Matlab} (a) displays four time-trajectories of the initially right-side spiral, 
$\V{\bar{x}}_{\epsilon_k}(t)$, with $\epsilon_0 = 0.004$, $\epsilon_1 = 0.002$, $\epsilon_2 = 0.001$,
$\epsilon_3 = 0.0005$.    
The marked spiral center locations were computed by us from the Matlab code's output 
every 0.1 time unit, and are interpolated here in time by cubic splines.  
The indicated location at $t=1.0$ on the $\epsilon=0.004$ curve appears to agree fairly 
well with our second-order method result displayed in Fig.~\ref{fig:Our_SingleGrid} (b), 
where we also used a single uniform grid discretizing $[-0.2,0.2]\times[-0.2,0.2]$.
Figure~\ref{fig:Matlab} (b-c) shows vorticity contours at $t=1.0$.
All data depicted in Fig.~\ref{fig:Matlab} illustrate that the distance between the twin spiral centers 
decreases with decreasing $\epsilon$.

\begin{figure}[!b]
  \centerline{\vbox{\hbox{\includegraphics[width=0.75\textwidth]{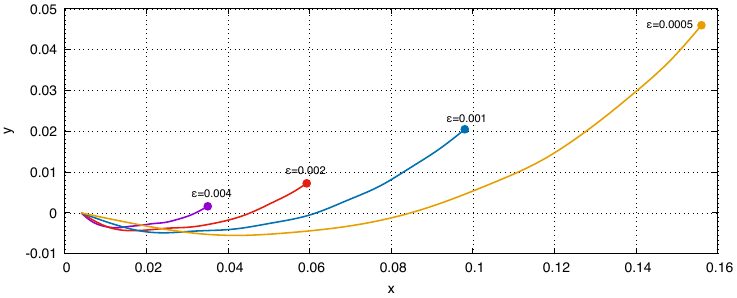}}
                             \hbox{(a) The referenced Matlab code.}}}
  \medskip
  \centerline{\vbox{\hbox{\includegraphics[width=0.75\textwidth]{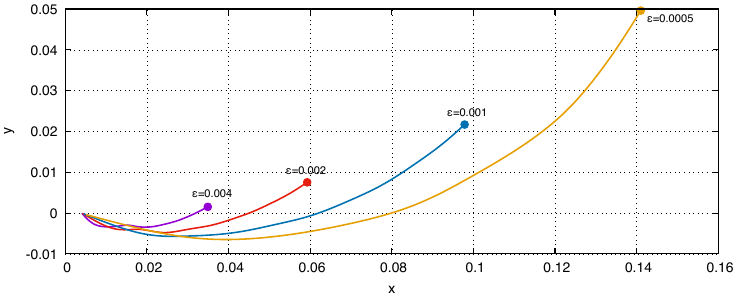}}
                             \hbox{(b) Our code.}}}
  \caption{{\em Scaled trajectories} (see (\ref{eq:traj-symm-mod2}) and its including paragraph)
  illustrating that the numerical solutions computed on the tiny domains, $[-0.2,0.2] \times [-0.2,0.2]$,
  fail to obey the scaling symmetry of the exact solution on $\mathbb{R}^2$. 
  If they did obey this symmetry, the trajectories marked $\epsilon=0.004, 0.002, \ldots$
  would overlay each other in both (a) and (b), as seen in Fig.~\ref{fig:overlay-mtraj}.
  The computations in (a) used $3200^2$ grid points, and in (b), $4096^2$ grid cells.
  }
  \label{fig:Matlab+Ours_overlay}
\end{figure}

The scaling symmetry discussed in \S\ref{sec:scaling} is a vital component of the analysis given in 
\cite{ab2021} and \cite{ab2020} and consequently here. 
Recall that from the scaling identity (\ref{traj-symmetry}) given in \S\ref{sec:scaling}, the trajectories for the 
exact solution to both problems posed on the entire plane $\mathbb{R}^2$ obey 
\begin{equation}
  (\epsilon_0/\epsilon)\, \V{x}_{\epsilon} ( t ) = \V{x}_{\epsilon_0}((\epsilon_0/\epsilon)^\alpha \, t).
  \label{eq:traj-symm-mod2}
\end{equation}
We take the four curves $\V{\bar{x}}_\epsilon(t)$ defined in the previous paragraph, 
generated from the Matlab code's results, and scale each by 
$\epsilon_0/\epsilon$, where $\epsilon_0 = 0.004$ and $\epsilon$ ranges over the other three values
$\epsilon=0.002$, $0.001$, $0.0005$.
In Fig.~\ref{fig:Matlab+Ours_overlay} (a) we plot 
\[
  (\epsilon_0/\epsilon)\, \V{\bar{x}}_{\epsilon} ( t ) \ \ \hbox{for $0 \le t \le 1$.}
\]
In Fig.~\ref{fig:Matlab+Ours_overlay} (b) we plot 
\[
  (\epsilon_0/\epsilon)\, \V{\bar{\bar{x}}}_{\epsilon} ( t ) \ \ \hbox{for $0 \le t \le 1$,}
\]
where the four curves denoted by $\V{\bar{\bar{x}}}_{\epsilon} ( t )$ are obtained using our second-order
method from \S\ref{sec:num-meth}, also using a single uniform grid of $[-0.2,0.2]\times[-0.2,0.2]$.  
The curves in (a) and (b) corresponding to $\epsilon = 0.004$, $0.002$, and $0.001$ appear to
agree fairly well, while the curves labeled $\epsilon = 0.0005$ differ markedly for larger $t$.  
This is almost certainly due to the order of the methods:  
our method displayed in (b) is second order in space 
and time, while the referenced Matlab code (a) is only 
first order in space and time.  
If the computed numerical solutions respected the scaling symmetry (\ref{eq:traj-symm-mod2}), 
the curves in both (a) and (b) would accurately overlay each other, as seen earlier in Fig.~\ref{fig:overlay-mtraj} for our nested-grid solutions computed on a much larger domain of size $[-64,64] \times [-64,64]$.
However, neither the curves in (a) nor (b) overlay at all.  
This is precisely why we have used the nested-grid approach. 
Recall equation (\ref{nested-domains}) and the paragraph including it in \S\ref{sec:num-meth}.  
Here we find that taking $L_0=64$ provides us a sufficiently large domain in order to adequately
obey the scaling symmetry (\ref{eq:traj-symm-mod2}).


\section{Conclusion}
\label{sec:conclusion}

We have presented numerical solutions of a Cauchy problem for the incompressible Euler equations
with initial data regularized in two different manners, subject to a parameter $\epsilon$,  
that converge strongly to the same limit as $\epsilon \to 0$. 
This problem was formulated in \cite{ab2021} and \cite{ab2020},  
and in both papers it was stated 
that numerical solutions show that as $\epsilon \to 0$  
two distinct solutions are obtained, one containing a single spiral and the other two spirals.  
In an analogous, small Mach number compressible study \cite{bjl2021}, 
numerical results were presented and similar statements are made there.
In the present work we have used a high-resolution finite volume scheme 
for the incompressible problem using extreme grid refinement. 
Since self-similar solutions were a critical component of the analysis given in
\cite{ab2021} and \cite{ab2020}, and self-similar solutions require that the Cauchy problem
be posed on all of $\mathbb{R}^2$, in our numerical computations we have effectively employed 
a very large spatial domain on which the requisite scale symmetry is accurately satisfied.
In almost every set of parameter values ($\alpha,\,\theta_w$) considered here, 
we do not find evidence of self-similar two-spiral solutions coming from either Problem~2 or Problem~1. 
We find that solutions of Problem~1, in every case, appear to tend in the $\epsilon$-limit to a single-spiral solution.  
All solutions of Problem~2 also appear to tend to a single-spiral self-similar solution, 
in all but one exceptional case:  $\alpha=0.95$, $\theta_w=\pi/16$.  
If what we see in this case is, in fact, approaching a two-spiral self-similar limit solution, 
which we can neither confirm nor refute with certainty, 
this indicates that whether or not the $\epsilon$-limit of Problem~2 evolves into a two-spiral 
solution depends on a subtle connection between problem parameters.

\section*{CRediT authorship contribution statement}
Allen Tesdall, Richard Sanders:  Equal contribution.

\section*{Declaration of competing interest}
The authors declare that they have no known competing financial interests or personal relationships that could have appeared to influence the work reported in this paper.

\section*{Data availability}
The data that support the findings of this study are available from the corresponding author upon reasonable request.

\section*{Acknowledgements}
A.~Tesdall was supported by NSF under grant DMS-1516131, and by 
CUNY Research Foundation under PSC-CUNY Awards 65263-00 53 and 67177-00 55.  
The authors wish to thank John K.~Hunter for his comments and advice on the current
manuscript. 
A.~Tesdall is grateful to the University of Houston, Department of Mathematics, for
hosting his sabbatical visit during academic year 2022-2023.

\bibliographystyle{elsarticle-num-names} 
\bibliography{sn-bibliography}





%
\end{document}